\documentclass[aap]{imsart}
 
 \RequirePackage{amsthm,amsmath,amsfonts,amssymb}
 \RequirePackage[numbers,sort&compress]{natbib}
 \RequirePackage[colorlinks,citecolor=blue,urlcolor=blue]{hyperref}
 \RequirePackage{graphicx}
  
 \usepackage{tikz} 
 \usepackage{mathrsfs,enumitem} 
 
 \usepackage{pgfplots}
 
 \startlocaldefs
 \theoremstyle{plain}

 \newtheorem{theorem}{Theorem}[section]
 \newtheorem{lemma}[theorem]{Lemma}
 \theoremstyle{remark}
 \newtheorem{definition}[theorem]{Definition}
 \newtheorem*{example}{Example}

 \newtheorem{corollary}[theorem]{Corollary}
 \newtheorem{proposition}[theorem]{Proposition}
 \newtheorem{remark}[theorem]{Remark}
 \newtheorem{condition}[theorem]{Condition}

  \makeatletter
 \def\namedlabel#1#2{\begingroup
 	#2%
 	\def\@currentlabel{#2}%
 	\phantomsection\label{#1}\endgroup
 }
 \makeatother

 \def\beqlb{\begin{eqnarray}}\def\eeqlb{\end{eqnarray}} 
 \def\beqnn{\begin{eqnarray*}}\def\eeqnn{\end{eqnarray*}} 
 \def\ar{&} 

 \def\proof{\noindent{\it Proof.~~}} 
 \def\qed{\hfill$\Box$\medskip}
 
 \endlocaldefs
 
\begin{document}
 
 \begin{frontmatter}
 \title{Diffusion Approximations for Self-Excited Systems with Applications to General Branching Processes}
 \runtitle{Self-Excited Systems and General Branching Processes}

 \begin{aug}
 \author{\fnms{Wei}~\snm{Xu}\ead[label=e1]{xuwei.math@gmail.com}\orcid{0000-0002-9370-8846}},
 \address{School of Mathematics and Statistics, Beijing Institute of Technology\printead[presep={,\ }]{e1}} 
 \end{aug}
 
 \begin{abstract}
 In this work, several convergence results are established for nearly critical self-excited systems in which event arrivals are described by multivariate marked Hawkes point processes. 
 Under some mild high-frequency assumptions, the rescaled density process behaves asymptotically like a multi-type continuous-state branching process with immigration, which is the unique solution to a multi-dimensional stochastic differential equation with dynamical mechanism similar to that of multivariate Hawkes processes.
 To illustrate the strength of these limit results, we further establish diffusion approximations for multi-type Crump-Mode-Jagers branching processes counted with various characteristics by linking them to marked Hawkes shot noise processes.
 In particular, an interesting phenomenon in queueing theory, well-known as state space collapse, is  observed in the behavior of the population structure at a large time scale. This phenomenon reveals that the rescaled complex biological system can be recovered from its population process by a lifting map.
 \end{abstract}

 \begin{keyword}[class=MSC2020]
 \kwd[Primary ]{60F17}
 \kwd{60G55}
 \kwd[; secondary ]{60J80}
 \kwd{62P10}
 \end{keyword}

 \begin{keyword}
 \kwd{scaling limit}
 \kwd{marked Hawkes point measure}
 \kwd{shot noise process}
 \kwd{multi-type continuous-sate branching process}
 \kwd{Crump-Mode-Jagers branching process}
 \kwd{population structure}
 \kwd{state space collapse}
 \end{keyword}

 \end{frontmatter}
 \tableofcontents

  \section{Introduction}
 
 This paper is concerned with developing a diffusion approximation for a stochastic dynamical system enjoying self-exciting property.
 In such a system, events are likely to not only occur in clusters but also mutually depend on each other.
 To capture both the self-exciting property and the clustering effect, we model the event arrivals with a \textit{multivariate marked Hawkes point measure with homogeneous immigration} (multivariate MHPI-measure) on $(0,\infty)\times\mathbb{U}$, denoted by $N_\mathcal{H}(dt,du):=(N_i(dt,du))_{i\in \mathcal{H}}$, where $\mathbb{U}$ is a measurable space and  $\mathcal{H}:=\{1,2,\dots , d  \}$ for some $d\in \mathbb{Z}_+$.
 To be precise, the random point measure $N_i(dt,du)$ has a predictable intensity $\Lambda_i(t-)\cdot dt\cdot \nu_i(du)$ in which $\nu_i(du)$ is a probability law on $\mathbb{U}$ and
 \beqlb\label{Intensity}
 \Lambda_i(t)
 \ar:=\ar \mu_i(t) + \int_0^t \int_{\mathbb{U}}\phi_i(t-s,u) N_I(ds,du) \cr
 \ar\ar +  \sum_{j\in\mathcal{H}} \int_0^t \int_{\mathbb{U}} \phi_i(t-s,u) N_j(ds,du), \quad  t\geq 0 
 \eeqlb
 for some non-negative functional-valued random variable $\{  \mu_\mathcal{H}(t):=(\mu_i(t))_{i\in\mathcal{H}} :t\geq 0 \}$, {kernel}\footnote{  \label{Footnote.1}
 	For different kernels $(\phi_{i,j})_{i,j\in\mathcal{H}}$ and $(\phi_{i,I})_{i\in\mathcal{H}}$, we can extend the mark space to $\bar{\mathbb{U}} :=  \mathbb U \times \{1,\cdots,d,I\}$ and take $\bar\nu_j(d\bar{u}):= \nu_j(d \bar{u}_1) \delta_{j} (d \bar{u}_2)$ and $\phi_i(t,\bar{u}):=\sum_{l\in\mathcal{H}}\phi_{i,l}(t, \bar{u}_1)\mathbf{1}_{\{\bar{u}_2=l\}} + \phi_{i,I}(t, \bar{u}_1)\mathbf{1}_{\{\bar{u}_2=I\}}$ for $i\in\mathcal{H}$ and $j\in\mathcal{H}\cup\{ I \}$.} 
 $\phi_\mathcal{H}:= (\phi_{i})_{i\in\mathcal{H}} : \mathbb{R}_+\times \mathbb{U} \to \mathbb{R}_+^d$, Poisson random measure $N_I(dt,du)$ on $(0,\infty)\times\mathbb{U}$ with intensity $\lambda_I\cdot dt\cdot \nu_I(du)$ for some constant $\lambda_I\geq 0$ and probability law $\nu_I(du)$ on $\mathbb{U}$; more accurate definitions can be found in Section~\ref{Section.1}.
 Usually, $\mu_\mathcal{H}$ is interpreted as the impact of all events prior to time $0$ on the arrivals of future events.
 This multivariate MHPI-measure includes both self/mutually-excited jumps ($N_\mathcal{H}$) and externally excited jumps ($N_I$), which respectively model the impact of endogenous and exogenous factors of the underlying system.
 It can be considered as an extension of marked Hawkes processes (i.e., $\mathbb{U}=\mathbb{R}_+$, $\mu_\mathcal{H}$ is a vector and $\lambda_I=0$) introduced by Ogata \cite{Ogata1998} for the study of different effects of earthquakes of different magnitudes on the arrivals of the future earthquakes; see also \cite{BremaudMassoulie2002, BremaudNappoTorrisi2002}  for the case of abstract-valued marks.
 Specially, when $\mathbb{U}=\mathbb{R}_+$ and $\phi_i(t,u)=ue^{-\beta_i t}$ for some $\beta_i>0$,
 the embedded point process $\{N_\mathcal{H}(t):= N_\mathcal{H}([0,t],\mathbb{R}_+):t\geq 0\}$ turns to be a multivariate version of dynamical contagion process given in \cite{DassiosZhao2011}. Moreover, when $\mu_\mathcal{H}$ is a vector, $\lambda_I=0$ and the kernel is mark-independent, the point process $N_\mathcal{H}$ reduces to a classical multivariate Hawkes process, which was firstly introduced by Hawkes \cite{Hawkes1971a,Hawkes1971b}.
 
 As Hawkes processes are always able to provide convincing interpretations of the cascade phenomenon and clustering effect that have been widely observed in various fields (e.g., financial contagion (see \cite{Ait-SahaliaCacho-DiazLaeven2015}) and
 credit contagion (see \cite{JorionZhang2009})), their applications have nowadays gone far beyond the original purpose of modeling earthquakes and their aftershocks; readers may refer to  \cite{BacryMastromatteoMuzy2015} for reviews on Hawkes processes and their applications.
 In particular, since they were firstly used in the estimation of value-at-risk (see \cite{Chavez-DemoulinDavisonMcNeil2005}) and modelling market events (see \cite{Bowsher2007}), various financial models have been established in the Hawkes framework to investigate the foreign exchange rates (see \cite{Hewlett2006}), mid-quote prices (see \cite{BacryDelattreHoffmannMuzy2013a,BacryDelattreHoffmannMuzy2013b}), limit order books (see \cite{HX2017,Large2007}), stochastic volatility (see \cite{ElEuchFukasawaRosenbaum2018,JaissonRosenbaum2015,JaissonRosenbaum2016}) and so on.
 Readers are suggested to refer to the seminal references of Bacry, Rosenbaum and their coauthors for various micro-structure models and macroscopic models.
 
 Different from Hawkes-based models, stochastic models driven by marked Hawkes processes/measures are \textit{individual-based models}, also called \textit{agent-based models}, in which a high degree of complexity and differences of events is allowed, and each event has a set of state variables or attributes and behaviors.
 An advantage of marked Hawkes-based models over Hawkes-based models is that they can incorporate any number of event-level mechanisms.
 Therefore, they are usually more effective in the modelling of complex dynamical systems.
 For instance, Horst and Xu \cite{HX2022} used a class of MHPI-measures with exponential kernel to study stochastic volatility models with self-exciting jump dynamics.
 In this case, each order is associated with a mark from the space $\mathbb{U}=\mathbb{Z}\times\mathbb{R}_+$ that describes the changes of the price in ticks caused by the order along with its impact on the arrival dynamics of future orders.
 Because of the significant impact of some orders on the arrivals of future orders, jumps occur in the high-frequency limit volatility models.
 This never happens in the high-frequency limits of Hawkes-based models; see \cite{ElEuchFukasawaRosenbaum2018,JaissonRosenbaum2015}.
 In another example that illustrates the advantage of marked Hawkes-based models, Xu \cite{Xu2021} generalized a classical second Ray-Knight theorem to a spectrally positive stable process by linking the intrinsic branching structure of its local time to a MHPI-measure with kernel being a unit step function.
 More precisely,  each individual in the population is endowed with a mark from the space $\mathbb{U}=\mathbb{R}_+$ to represent its life-length and its survival state is described by the kernel of the form  $\mathbf{1}_{\{u>t\}}$ (i.e., it is alive when its life-length $u$ is larger than its age $t$).
 Furthermore, to emphasize the necessity of setting $\mathbb{U}$ to be an abstract space, in Section~\ref{Section.CMJ} we develop a new way to study the general branching particle systems by linking them to MHPI-measures, in which the abstract-valued marks represent individuals' characteristics, e.g., life-length, reproduction process, impact on host and so on.
 
 Similar to Hawkes processes, the MHPI-measure $N_\mathcal{H}(dt,du)$ can be constructed in collaboration with a labeled birth-immigration particle system, in which the embedding multi-type Galton-Watson process with immigration (GWI-process) has \textsl{mean matrix} $\|\phi_{\mathcal{H}^2}\|_{L^1}$ that is the $L^1$ norm of $\big\{\phi_{\mathcal{H}^2}(t):= (\phi_{ij}(t))_{i,j\in\mathcal{H}}:t\geq 0 \big\}$ with 
 \beqnn
 \phi_{ij}(t):= \int_\mathbb{U} \phi_i(t,u)\nu_j(du).
 \eeqnn 
 By the elementary theory of branching processes; see Chapter~V in  \cite{AthreyaNey1972}, the GWI-process is \textsl{subcritical},  \textsl{critical} or \textsl{supercritical} if the mean matrix $\|\phi_{\mathcal{H}^2}\|_{L^1}$ has spectral radius $\varrho<1$,  $=1$ or  $>1$ respectively.
 Furthermore, the stationary distribution exits if it is subcritical or critical with sparse immigrants (in this case the stationary distribution does not have finite mean). 
 In the supercritical regime, it grows exponentially to infinity.
 Therefore, analogous to Hawkes processes; see \cite{BremaudMassoulie1996}, the condition $\varrho<1$ is necessary for $N_\mathcal{H}(dt,du)$ to own an asymptotically stationary intensity process with finite first moment.
 
 Nowadays, because a significant part of financial transactions is carried out through electronic order books, high-frequency trading has enjoyed a growing popularity.
 This has made high-frequency financial models including Hawkes-based models receive considerable attention in the probability and financial mathematics literature in recent years.
 Two types of typical and important limit theorems have been widely established to study the behavior of Hawkes-based models at a large time scale.
 The first one mainly consists of functional law of large numbers (FLLN) and functional central limit theorem (FCLT), which were firstly established by Bacry et al. \cite{BacryDelattreHoffmannMuzy2013b} for a multivariate Hawkes process whose kernel enjoys short-memory property and spectral radius $\varrho$ is strictly smaller than one; also see \cite{HX2019a} for the case of MHPI-measures. Recently, Horst and Xu \cite{HX2023,HX2024} established a full FLLN and FCLT for subcritical and critical uni-variate Hawkes processes. 
 The second kind, usually known as scaling limit theorem, was firstly investigated by
 Jaisson and Rosenbaum \cite{JaissonRosenbaum2015} in the study of asymptotic behavior of Hawkes-based price models in the context of high-frequency trading.
 Their results state that under some short-memory condition, the rescaled intensity of nearly unstable Hawkes process converges weakly to a Feller diffusion (also known as CIR-model in finance).
 Different from the deterministic limit in FLLN and the Brownian motion in FCLT, CIR-model  inherits not only the randomness but also the self-exciting property from Hawkes processes.
 Under a heavy-tailed condition, they also proved that  the rescaled Hawkes point process converges weakly to the integral of a rough fractional diffusion; see \cite{JaissonRosenbaum2016}. A more refined convergence result has recently been established by Horst et al. \cite{HXZ2023}, which proved the weak convergence of the rescaled intensities, instead of their integrals, to a rough fractional diffusion.
 An analogous scaling limit was later established by El Euch et al. \cite{ElEuchFukasawaRosenbaum2018} for multivariate Hawkes processes with positive and diagonalizable kernel, see \cite{RosenbaumTomas2021} for the case of trigonalizable kernel.
 In addition, a jump-diffusion limit was provided in \cite{HX2022} for  MHPI-measures with real-valued mark and exponential kernel.
 
 In the first part of this work, we mainly investigate the behavior at a large time scale of stochastic dynamical systems driven by asymptotically critical multivariate MHPI-measures, which corresponds to the assumption that  $\|\phi_{\mathcal{H}^2}\|_{L^1}$ is close to a limit matrix with unit spectral radius.
 In addition, compared to the uni-variate case the asymptotic criticality for multivariate Hawkes processes/measures is much more complicated, since there are infinite possibilities for the limit matrix and meanwhile the limit of rescaled Hawkes based models varies greatly for different limit matrix.
 For instance, in the case of limit matrix being positive and diagonalizable, El Euch et al. \cite{ElEuchFukasawaRosenbaum2018} proved that the rescaled intensity process of multivariate Hawkes process gradually concentrates in one direction of $\mathbb{R}^d$ and the limit process is the multiplication of one-dimensional CIR-model by a vector.
 By contrast, the mean matrix $\|\phi_{\mathcal{H}^2}\|_{L^1}$ in our setting is assumed to converge to an identity matrix. 
 Under some short-memory conditions, we show that the rescaled intensity process asymptotically behaves like a multi-type continuous-state branching process with immigration (CBI-process) that is defined as the unique strong solution of a $d$-dimensional stochastic differential equation (SDE) with linear drift and $1/2$-H\"older continuous diffusion.
 
 Our assumption of the convergence of $\|\phi_{\mathcal{H}^2}\|_{L^1}$ to an identity matrix stems from three main reasons.
 Firstly, under this assumption both proofs and statements can be dramatically simplified. 
 Moreover, by using the rotation method developed in \cite{ElEuchFukasawaRosenbaum2018,RosenbaumTomas2021}, our proofs can be generalized to the case of trigonalizable limit matrices and the corresponding limit theorems can be established similarly.
 Secondly, in practice, the self-excitation is generally much stronger than the mutual-excitation, e.g. market, limit and cancel orders in financial market are likely to effect themselves (see \cite{BacryJaissonMuzy2016}); the foregoing topic in social media communities is preferred to be discussed continuously (see \cite{PhillipsGorse2018}).
 These are consistent with our assumption, i.e., $\|\phi_{ij}\|_{L^1} <<\|\phi_{ii}\|_{L^1}$ for $i\neq j$.
 It seems that the mutual-excitation in our setting can be asymptotically ignored, but its impact on the underlying system still can be observed in limit process.
 Finally, the diagonal entries of drift matrix in the limit $d$-dimensional SDE represents the net self-excitation. Moreover, the off-diagonal entries are non-negative and can be interpreted as the mutual-excitation.
 Hence, compared with the limit model established in \cite{ElEuchFukasawaRosenbaum2018} as the multiplication of one-dimensional CIR-model by a vector, our limit model is a more natural continuous version analogous in form to the stochastic dynamical system (\ref{Intensity}).
 
 The main results in the first part of this paper are proved by extending the method developed in \cite{JaissonRosenbaum2015}, in which the rescaled intensity of Hawkes process is rewritten in the form of an It\^o's SDE and then the limit theorem is proved by using the convergence results established in \cite{KurtzProtter1991} for finite-dimensional stochastic integrals.
 However, in addition to the technical difficulties encountered in  \cite{ElEuchFukasawaRosenbaum2018}, three main and new challenges are induced by our setting. 
 Firstly, the high degree of complexity and differences of events gives raise to not only some additional perturbations in the dynamics of intensity processes but also severe fluctuation in error processes.
 These make it much more difficult to establish a rigorous connection between multivariate MHPI-measures and multi-type CBI-processes.
 Secondly, as two sufficient conditions for the convergence results in \cite{KurtzProtter1991}, the weak convergence and uniform tightness of driving noises were proved easily based on the fact that intensities of nearly unstable Hawkes processes in \cite{JaissonRosenbaum2015} are uniformly strictly larger than zero.
 By comparison, they are extremely difficult to be identified in our setting, since the exogenous intensity $\mu_\mathcal{H}$ varies as time goes and the intensity process $\Lambda_\mathcal{H}$ may hit zero in finite time.
 Thirdly, the stability condition, widely considered in the Hawkes literature (see \cite{BacryDelattreHoffmannMuzy2013b,ElEuchFukasawaRosenbaum2018,JaissonRosenbaum2015, JaissonRosenbaum2016}), is not assumed in this work and MHPI-measures are allowed to be unstable ($\varrho >1$).
 In this case, both the resolvent and the intensity process may grow exponentially to infinity, which make the error estimates more difficult.
 To overcome the first two difficulties, we start by reconstructing MHPI-measures in collaboration with Poisson random measures.
 The mutual-excitation is interpreted as the non-local branching mechanism in the corresponding birth-immigration particle system. 
 Under the assumption that $\|\phi_{\mathcal{H}^2}\|_{L^1}$ converges to an identity matrix, we further consider it as a state-dependent immigration and translate the mutually-excited jumps into another kind of externally excited jumps.
 Inspired by the computations and techniques applied in \cite{HX2019a,HX2022}, 
 associated to the resolvent we introduce a two-parameter function to describe the average impact of an event with some mark on the future intensity. 
 It enables us to write the stochastic equation (\ref{Intensity}) approximately as an It\^o's SDE driven by an infinite-dimensional semimartingale.
 In particular, this semimartingale mainly consists of several compensated Poisson random measures whose weak convergence and uniform tightness follow immediately from their orthogonal increments.
 Moreover, with the help of the foregoing two-parameter function, our error analysis is successfully carried out through investigating the exact perturbations of each event of various marks in the error processes.
 The desired limit theorem for intensity processes is finally obtained by using the convergence results of infinite-dimensional stochastic integrals established by Kurtz and Protter \cite{KurtzProtter1996}.
 For the third difficulty, the exponential growth of the resolvent and intensity processes of supercritical MHPI-measures encourages us to modify the self-excited dynamical system by an exponential function.
 Due to the multiplicative property of exponential functions, the preceding representations and asymptotic analysis remain valid with some slight modifications.
 
 As mentioned above, the main contribution in the second part of this work is to illustrate the strength of the foregoing limit results for MHPI-measures by applying them to study the behavior at a large time scale of \textit{multi-type Crump-Mode-Jagers branching processes with immigration} (CMJI-processes).
 In the realistic pattern, as Peter Jagers \cite{Jagers2010} pointed out,  population models ``must be ultimately stochastic [...] individual based [...] life span can have an arbitrary distribution [...] reproduction should be modelling as it actually occurs''.
 As a result,  CMJI-processes, as a class of continuous-time and discrete-state stochastic population models with age-dependent reproduction mechanism, have received considerable attention in the probability and mathematical biology literature since they were firstly introduced in \cite{CrumpMode1968,CrumpMode1969,Jagers1969}.
 However, because they are generally neither Markov nor semimartingales, the instruments provided by modern probability theory are almost out of work and hence researches concerning CMJI-processes are relatively less than those of Markovian population models; see \cite{AthreyaNey1972,Jagers1975}.
 To the best of our knowledge, only few asymptotic results have been established for  CMJ-processes up to now, e.g.,
 a scaling limit was established by Lambert et al. \cite{Lambert2013} for homogeneous, binary and  single-type CMJ-processes without immigration via connecting them to the local time of compound Poisson processes.
 However, this connection can not be generalized to CMJ-processes with general branching mechanism and complex individual characteristics. 
 
 Here we develop a new way to investigate multi-type CMJI-processes by linking them to multivariate MHPI-measures.
 More precisely, we translate random point measures $N_I(dt,du)$ and $N_\mathcal{H}(dt,du)$ into the arrivals of immigrants and offspring respectively, and marks from an abstract space into individuals' information, e.g., size, type, life-length, reproduction process and characteristic.   
 Different to the existing literature (e.g. \cite{Lambert2013,Sagitov1995}) in which the population size is usually studied in the first place, we start by considering the asymptotic behavior of total reproduction rate of alive individuals that coincides with the intensity process of $N_\mathcal{H}(dt,du)$. 
 Using the preceding limit results for multivariate MHPI-measures, we show that with a suitable scaling, the total reproduction rate process behaves asymptotically as a multi-type CBI-process. 
 Additionally, for a multi-type CMJI-process counted with various characteristics (e.g., population size and total progeny), we link it to a shot noise process driven by $N_\mathcal{H}(dt,du)$ and then show that with a suitable scaling, it can be well approximated by a functional of the multi-type CBI-process. 
 Furthermore, an interesting phenomenon, known as \textit{state space collapse} in queueing theory, is observed in the population structure at a large time scale.  
 In precise, the joint distribution of age and residual life of alive individuals can be recovered directly from the population size by an appropriate lifting map.
 This indicates that when the life-length and the reproduction process enjoy short-memory property, more detailed information about the population,  except the  life-length distribution and the mean/variance of offspring, is not necessary for the study of complex biological systems. 
 
 \smallskip
 \textbf{\textit{Organization of this paper.}} In Section~\ref{Section.1} we provide a branching representation as well as a stochastic Volterra integral representation for multivariate MHPI-measures. A criticality criterion for their stationarity is also given as a byproduct.
 In Section~\ref{SLMHP}, we establish several limit theorems for stochastic dynamical systems driven by multivariate MHPI-measure, including the weak convergence of rescaled intensity processes to a multi-type CBI-process and scaling limits for marked Hawkes shot noise processes.
 In Section~\ref{Section.CMJ} we apply these limit results for self-excited systems to establish diffusion approximations for multi-type CMJI-processes.
 Section~\ref{Proofs} is devoted to the proofs for all limit theorems given in this work.
  
 \smallskip
 \textbf{\textit{Notation.}} 
 Denote by $[x]$ the integer part of  the real number $x\in\mathbb{R}$ and $z_{\mathcal{H}}= (z_i)_{i\in\mathcal{H}}=(z_1,\cdots, z_d)$ with $|z_{\mathcal{H}}|:= |z_1|+\cdots+|z_d|$.
 Let $f*g$ be the convolution of two functions $f,g$ on $\mathbb{R}_+$ and $f^{*n}$ be the $n$-th convolution of $f$. 
 For $h > 0$, we write 
 \beqnn
 \Delta_-f(x):=f(x)-f(x-) 
 \quad \mbox{and}\quad 
 \Delta_hf(x):=f(x+h)-f(x).
 \eeqnn 
 Let $\|f\|_{\rm TV}$ denote the total variation of $f$. 
 For any $p,q\in(0,\infty]$, let $L^{p,q}(\mathbb{R}_+)= L^{p}(\mathbb{R}_+)\cap L^{q}(\mathbb{R}_+)$ with norm $\|\cdot\|_{L^{p,q}}:= \|\cdot\|_{L^p}+ \|\cdot\|_{L^q}$.
 For $T>0$, let 
 \beqnn
 \|f\|_{L_T^\infty}:= \sup_{t\in[0,T]}|f(t)|
  \quad \mbox{and}\quad 
 \|f\|^q_{L_T^q}:= \int_0^T |f(t)|^q dt
 \eeqnn
  
 Denote by $\overset{\rm u.c.}\to$ ,$\overset{\rm u.c.p.}\longrightarrow $, $\overset{\rm a.s.}\to$ ,$\overset{\rm d}\to$, $\overset{\rm p}\to$ and $\overset{\rm f.d.d.}\longrightarrow $ the compact convergence, uniform convergence on compacts in probability, almost sure convergence, convergence in distribution, convergence in probability and convergence in the sense of finite-dimensional distributions.
 We also use  $\overset{\rm a.s.}=$, $\overset{\rm d}=$ and $\overset{\rm p}=$ to denote almost sure equality, equality in distribution and equality in probability.
 
 Given a measurable space $(\mathbb{V},\mathscr{V})$, let $B(\mathbb{V})$, $C(\mathbb{V})$ and $C_0(\mathbb{V})$ be the spaces of measurable functions on $\mathbb{V}$ that are  bounded, continuous and  continuous as well as vanishing at infinity respectively.
 For $T\in[0,\infty]$, denote by $\mathbf{D}([0,T],\mathbb{V})$ the space of c\'adl\'ag functions from $[0, T ]$ to $\mathbb{V}$ furnished with the Skorokhod topology.
 Let $\mathcal{M}(\mathbb{V})$ be the space of finite Borel measures on $\mathbb{V}$ equipped with the weak convergence topology. 
 Let $\delta_a$ be the Dirac measure at point $a\in\mathbb{V}$.
 For any $\upsilon\in \mathcal{M}(\mathbb{V})$ and $f \in B(\mathbb{V})$, we write
 \beqnn
 \upsilon(f) := \int_\mathbb{V} f(x) \upsilon(dx)
 \quad\mbox{and}\quad 
 f*\upsilon := \int_\mathbb{V} f(x-y)\upsilon(dy).
 \eeqnn
 
 Throughout this paper, we assume the generic constant $C$ may vary from line to line.

 \section{Preliminaries}\label{Section.1}

 Let $(\Omega,\mathscr{F},\mathbf{P})$ be a complete probability space endowed with a filtration $\{\mathscr{F}_t:t\geq 0\}$ satisfying the usual hypotheses and $(\mathbb{U},\mathscr{U})$ be a measurable space. 
 Let $\mathcal{H}:=\{1,2,\cdots,d \}$ for some $d\in\mathbb{Z}_+$.
 For $i\in \mathcal{H}$,  let $\{\tau_{i,k}\}_{k\geq 1}$  be a sequence of increasing, $(\mathscr{F}_t)$-stopping times and $\{ \xi_{i,k} \}_{k\geq 1}$ be  a sequence of i.i.d. $\mathbb{U}$-valued random variables with distribution $\nu_i(du)$ satisfying that $\xi_{i,k}$ is independent of $\mathscr{F}_{\tau_{i,k}}$ for any $k\geq 1$. 
 Associated to these two sequences we define an $(\mathscr{F}_t)$-random point measure on $(0,\infty)\times \mathbb{U}$ 
 \beqlb\label{eqn1.01}
 N_i(dt,du):= \sum_{k=1}^\infty \mathbf{1}_{\{\tau_{i,k}\in dt,\, \xi_{i,k}\in du  \}}.
 \eeqlb
 We say the random point measure $N_{\mathcal{H}}(dt,du):=(N_i(dt,du))_{i\in\mathcal{H}}$ is a \textit{multivariate marked Hawkes point measure} (multivariate MHP-measure) on $(0,\infty)\times \mathbb{U}$ with embedded point process $\{N_\mathcal{H}(t):=(N_i((0,t],\mathbb{U}))_{i\in\mathcal{H}}:t\geq 0\}$, if $N_i(dt,du)$ has $(\mathscr{F}_t)$-intensity $\Lambda_i(t-)\cdot dt\cdot \nu_i(du)$ and the $(\mathscr{F}_t)$-intensity process $\Lambda_i$ is of the form 
 \beqlb\label{eqn1.03}
  \Lambda_i(t)\ar=\ar \Lambda_{0,i}(t) + \sum_{j\in\mathcal{H}} \sum_{k=1}^{N_{j}(t)} \phi_{i}(t-\tau_{j,k}, \xi_{j,k}), \quad  t\geq 0,\ i\in\mathcal{H},
 \eeqlb
 for some non-negative, locally integrable, $(\mathscr{F}_t)$-progressive  \textit{exogenous intensity} $\Lambda_{0,\mathcal{H}}:=  ( \Lambda_{0,i})_{i\in\mathcal{H}}$,  and some \textit{kernel} 
 $\phi_\mathcal{H}:= (\phi_{i})_{i\in\mathcal{H}} : \mathbb{R}_+\times \mathbb{U} \to \mathbb{R}_+^d$. 
 Specially, we call $N_\mathcal{H}(dt,du)$ a \textit{multivariate MHP-measure with homogeneous immigration} (multivariate MHPI-measure) if 
 $\Lambda_{0,\mathcal{H}}$ admits the representation:
 \beqlb\label{eqn1.04}
 \Lambda_{0,i}(t)\ar:=\ar \mu_{i}(t) + \sum_{k=1}^{N_I(t)} \phi_{i}(t-\tau_{I,k}, \xi_{I,k}),
 \quad t\geq 0,\, i\in\mathcal{H},
 \eeqlb
 where $N_I$ is a Poisson process  with rate $\lambda_I$ and arrival times $\{\tau_{I,k} \}_{k\geq 1}$, $\{ \xi_{I,k}\}_{k\geq 1}$ is a sequence of i.i.d. $\mathbb{U}$-valued random variables with distribution $\nu_I(du)$ and independent of $N_I$, and
 $\mu_\mathcal{H}:=(\mu_i)_{i\in\mathcal{H}} $ is an $\mathscr{F}_0$-measurable $\mathbf{D}([0,\infty),\mathbb{R}_+^d)$-valued random variable. 
 For  $i\in\mathcal{H}$, let 
 \beqnn
 \mathcal{D}:=\mathcal{H}\cup\{I\},\quad 
 \mathcal{H}_i:=\mathcal{H}\setminus\{i \}
 \quad \mbox{and}\quad
 \mathcal{D}_i:=\mathcal{D}\setminus\{i \}.
 \eeqnn 
 For simplicity, we assume that $\lambda_I=1$ and $ \tau_{i,k}\neq\tau_{j,l} $ a.s. for $(i,k), (j,l) \in \mathcal{D}\times \mathbb{Z}_+$ with $(i,k)\neq (j,l)$. 
 We also refer all externally excited jumps as type~$I$ events.   
 Let $\phi_{\mathcal{H}^2}:=(\phi_{ij})_{i,j\in\mathcal{H}}$ and $\phi_{\mathcal{H}I}:=(\phi_{iI})_{i\in\mathcal{H}}$ with
 \beqlb\label{phi.ij}
 \phi_{ij}(t):= \int_{\mathbb{U}} \phi_i(t,u)  \nu_j(du),\quad t\geq 0,\,i\in\mathcal{H},\,j\in\mathcal{D}
 \eeqlb
 be the \textit{mean impact functions} of a type-$j$ event on the future arrivals of type-$i$ events. 
 In the sequel, we always assume that
 \beqnn
 \big\|\phi_{\mathcal{H}^2}\big\|_{L^1}:=\big(\big\|\phi_{ij}\big\|_{L^1}\big)_{i,j\in\mathcal{H}}<\infty
 \quad \mbox{and}\quad
 \big\|\phi_{\mathcal{H}I}\big\|_{L^1}:= \big(\big\|\phi_{iI}\big\|_{L^1}\big)_{i\in\mathcal{H}}<\infty.
 \eeqnn

 \subsection{Branching representation}\label{BrachingRep}
 In this section we show that the foregoing construction of multivariate MHPI-measure $N_\mathcal{H}(dt,du)$ can be done in collaboration with a multi-type birth-immigration particle system defined on the probability basis $(\Omega,\mathscr{F},\mathscr{F}_t,\mathbf{P})$ by the following properties:
 \begin{enumerate}
 	\item[\bf(A1)] There is an ancestor at time $0$, whose successive ages arrive according to a Cox point process with intensity process  $\big\{|\mu_\mathcal{H}(t)| :t\geq 0\big\}$.  
 	Only one child is born at each successive age. 
 	Conditioned on the birth time $t$, the child is type-$i$ and has mark $u\in\mathbb{U}$  with probability $\mu_{i}(t)\cdot|\mu_\mathcal{H}(t)|^{-1}\cdot \nu_i(du)$;
 	\smallskip
 	
 	\item[\bf(A2)] Immigrants enter into the population according to a Poisson process with rate $1$ and are endowed with a mark randomly and independently according to the probability law $\nu_I(du)$;
 	\smallskip
 	
 	\item[\bf(A3)] For each individual (except the ancestor) with mark $u\in \mathbb{U}$, it gives birth to a child at the rate $|\phi_\mathcal{H}(t,u)|$ at age $t$ .
 	Moreover, the child has probability  $\phi_{i}(t,u)\cdot|\phi_\mathcal{H}(t,u)|^{-1}$ to be type-$i$ and it picks up a mark according to the law $\nu_i(du')$. Moreover, all individuals produce their offspring independently.
 	
 \end{enumerate}
 
 Denote by $\mathcal{A}$  the collection of all individuals except the ancestor.
 Associated with each individual $x\in\mathcal{A}$ is a random triple $(\mathrm{t}'_x,\tau'_x,  u'_x)$ that represents its type, birth time and mark respectively.
 Define an $(\mathscr{F}_t)$-random point measure $N'_\mathcal{H}(dt,du):=(N'_i(dt,du))_{i\in\mathcal{H}}$ on $(0,\infty)\times\mathbb{U}$ with
 \beqnn
 N'_i(dt,du):= \sum_{x\in\mathcal{A}}\mathbf{1}_{\{\mathrm{t}'_x=i, \tau'_x\in dt,\, u'_x\in du \}}.
 \eeqnn
 The intensity of its embedded random point process $\big\{N'_\mathcal{H}(t):=N'_\mathcal{H}([0,t],\mathbb{U}):t\geq 0\big\}$, denoted by $\big\{\Lambda_\mathcal{H}'(t):=(\Lambda_i'(t))_{i\in \mathcal{H}}:t\geq 0\big\}$,  equals to the total birth rate of children of various types.
 In addition, by the branching property,  
 it admits the following representation
 \beqlb\label{BranchingRepresenation02}
 \Lambda'_i(t) 
 = \mu_i(t)+ \sum_{x\in\mathcal{A}_I} \phi_i(t-\tau'_x, u'_x)+\sum_{j\in\mathcal{H}} \sum_{x\in\mathcal{A}_j} \phi_i(t-\tau'_x, u'_x), \quad t\geq 0, i\in\mathcal{H},
 \eeqlb
 where $\mathcal{A}_I$ and $\mathcal{A}_j$ are the collections of all immigrants and type-$j$ offspring respectively.
 The following result can be obtained immediately by comparing (\ref{BranchingRepresenation02}) with (\ref{eqn1.03})-(\ref{eqn1.04}).
  
 \begin{proposition}
  The random point measure $N'_\mathcal{H}(dt,du)$ is a realization of the multivariate MHPI-measure defined by (\ref{eqn1.01})-(\ref{eqn1.04}).
 \end{proposition}
 
 The embedded point process  $N_\mathcal{H}$ (or $N'_\mathcal{H}$), also can be considered as a cluster process in which the process of cluster centres is the random point process formed by the arrivals of immigrants and the successive ages of the ancestor.
 The cluster at each centre is formed by all the descendants of an immigrant or a child of the ancestor.
 These clusters are mutually independent and identically distributed.
 Denote by $\{ X_{n,\mathcal{H}}:=(X_{n,i})_{i\in\mathcal{H}}:n=1,2,\cdots \}$ the embedded multi-type Galton-Watson process (GW-process) of a cluster produced by an immigrant.
 It is easy to see that elements of $X_{1,\mathcal{H}}$ are mutually independent and $X_{1,i}$ is Poisson distributed with rate $\|\phi_{iI}\|_{L^1}$.
 For $n\geq 2$,  $X_{n,i}$ is the number of type-$i$ individuals in the $n$-th generation, which can be written as
 \beqnn
 X_{n,i} =\sum_{j\in\mathcal{H}} \sum_{k=1}^{X_{n-1,j}} \xi_{n,j,k,i},\quad i\in\mathcal{H},
 \eeqnn
 with $\xi_{n,j,k,i}$ being the number of type-$i$ children born by the $k$-th type-$j$ individual in the $(n-1)$-th generation, which is Poisson distributed with parameter $\|\phi_{ij}\|_{L^1} $.
 Let $\varrho$ be the spectral radius of the  matrix $\|\phi_{\mathcal{H}^2}\|_{L^1} $ and $\mathbf{I}$ be an $d$-dimensional identity matrix.
 The mean cluster size equals to the mean of total progeny
 \beqnn
 \sum_{n=1}^\infty \mathbf{E} \big[\big| X_{n,\mathcal{H}}\big| \big] =   \bigg|\sum_{n=1}^\infty \|\phi_{\mathcal{H}^2}\|^{n-1}_{L^1} \cdot \|\phi_{\mathcal{H}I}\|_{L^1} \bigg|  = \Big| (\mathbf{I}- \|\phi_{\mathcal{H}^2}\|_{L^1})^{-1}\cdot \|\phi_{\mathcal{H}I}\|_{L^1} \Big| ,
 \eeqnn
 which is finite  if and only if $\varrho<1$. 
 The next proposition follows immediately from Theorem~3 and Corollary~3.2 in \cite{Westcott1971}.
 \begin{proposition} \label{CriticalityMHP}
 	If $\mu_\mathcal{H}$ is a non-negative constant vector and $\varrho<1$, the embedded point process $N_\mathcal{H}$ is asymptotically stationary.
 \end{proposition}
 Drawing from the criticality criterion for multi-type GW-processes, we say the multivariate MHPI-measure $N_\mathcal{H}(dt,du)$ is \textit{subcritical}, \textit{critical} or \textit{supercritical} if $\varrho<1$, $=1$ or $>1$.
 These correspond to the three phases  of a classical Hawkes process: \textit{stationary}, \textit{quasi-stationary} or \textit{non-stationary}; see \cite{BacryMastromatteoMuzy2015}. 
  
 \subsection{Stochastic Volterra representation}\label{Sec.SVR}
 We now provide a stochastic Volterra representation for the intensity process $\Lambda_\mathcal{H}$, which will play a considerably important role in the following asymptotic analysis. 
 Associated to the sequence $\{ (\tau_{I,k},\xi_{I,k})  \}_{k\geq 1}$ we define an $(\mathscr{F}_t)$-Poisson random measure
 \beqnn
 N_I(ds,du):= \sum_{k=1}^\infty \mathbf{1}_{\{\tau_{I,k}\in ds, \xi_{I,k}\in du  \}}
 \eeqnn
 on $(0,\infty)\times \mathbb{U}$ with intensity $dt\cdot \nu_I(du)$ and then rewrite the intensity process $\Lambda_\mathcal{H}$ under the form (\ref{Intensity}).
 Moreover, following the argument in \cite[p.93]{IkedaWatanabe1989}, on an extension of the original probability space we can define $d$ mutually orthogonal Poisson point measures $N_{0,i}(dt,du,dz),$ $ i\in\mathcal{H} $ on $(0,\infty)\times\mathbb{U}\times \mathbb{R}_+$ independent of $N_I(ds,du)$ such that  $N_{0,i}(dt,du,dz)$ has intensity $dt \cdot\nu_i(du)\cdot dz$ and
 \beqlb\label{Ni2N0i}
 \int_0^t\int_\mathbb{U} f(u) N_i(ds,du) = \int_0^t \int_\mathbb{U} \int_0^{\Lambda_i(s-)} f(u) N_{0,i}(ds,du,dz), \quad t\geq 0, i\in\mathcal{H},
 \eeqlb
 for any $ f\in B(\mathbb{U})$. 
 We can thus rewrite the last stochastic integral in (\ref{Intensity}) as
 \beqnn
 \int_0^t \int_\mathbb{U} \int_0^{\Lambda_i(s-)} \phi_i(t-s,u) N_{0,j}(ds,du,dz), \quad j\in\mathcal{H}.
 \eeqnn
 Actually, we can always construct multivariate MHPI-measures in collaboration with some Poisson random measures on the probability basis $(\Omega,\mathscr{F},\mathscr{F}_t,\mathbf{P})$; see Section~2 in \cite{HX2019a}.
 
 \begin{remark}
  If $\mu_\mathcal{H}$ is a positive constant vector and $\nu_{I}(\mathbb{U})=0$,  then $N_\mathcal{H}(dt,du_\mathcal{H})$  reduces to a classical multivariate MHP-measure (without immigration) on $(0,\infty)\times \mathbb{U}$.
  We now link it to a special multivariate MHPI-measure. For each $i\in\mathcal{H}$, let $ \Lambda^\circ_i:=\Lambda_i- \mu_i$ and 
  \beqnn
 	N^\circ_{i}(dt,du) \ar:=\ar  N_{0,i}(dt,du,[0,\Lambda^\circ_i(t-)),\cr
 	N_{I,i}^\circ(dt,du) \ar:=\ar N_{0,i}(dt,du,[\Lambda^\circ_i(t-),\Lambda_i(t-)).
  \eeqnn 
  Let $N_I^\circ(dt,du):=\sum_{j=1}^d N_{I,j}^\circ(dt,du)$, which is a Poisson random measure on  $(0,\infty)\times\mathbb{U}$ with intensity $ds \cdot\sum_{j\in\mathcal{H}} \mu_j \cdot \nu_j(du)$. 
  It is obvious that 
  \beqnn
  N_\mathcal{H}(dt,du)=N^\circ_\mathcal{H}(dt,du)+ N_{I,\mathcal{H}}^\circ(dt,du)
  \eeqnn  
  and $N^\circ_\mathcal{H}(dt,du)$ is a multivariate MHPI-measure on $(0,\infty)\times\mathbb{U}$ with intensity process $\Lambda^\circ_\mathcal{H}$ being of the form 
 	\beqnn
 	\Lambda^\circ_i(t)
 	\ar=\ar \int_0^t \int_\mathbb{U} \phi_i(t-s,u) N^\circ_I(ds,du) \cr
 	\ar\ar + \sum_{j\in\mathcal{H}} \int_0^t \int_\mathbb{U}\phi_i(t-s,u) N^\circ_j(ds,du),\quad t\geq 0, i\in\mathcal{H}.
 	\eeqnn

 \end{remark}
  
 To get the desired stochastic Volterra representation, we need several important quantities associated to $\phi_i$ and $\phi_{ij}$ for $i\in\mathcal{H}$ and $j\in\mathcal{D}$.
 Let $R_{ii}$ be the \textit{resolvent} of $\phi_{ii}$ defined by
 \beqlb
 R_{ii}(t)\ar=\ar \phi_{ii}(t) +  R_{ii}*\phi_{ii}(t),\quad t\geq 0. \label{Resolvent001}
 \eeqlb
 The existence and uniqueness of the solution $R_{ii}$ follow directly from the assumption $\|\phi_{ii}\|_{L^1}<\infty$ and Theorem~3.1 in \cite[p.32]{GripenbergLondenStaffans1990}.
 It is easy to identify that $R_{ii}$ admits the   representation
 \beqnn
 R_{ii}(t)= \sum_{k=1}^\infty \phi_{ii}^{*k}(t),\quad t\geq 0.
 \eeqnn
 It is usual to interpret $R_{ii}$ as the mean impact of a type-$i$ event and its triggered events on the future arrivals of type-$i$ events.
 In addition, the mean impact of a type-$j\in\mathcal{D}_i$ event and its triggered events on the future arrivals of type-$i$ events also can be described as
 \beqlb
 R_{ij}(t)\ar :=\ar \phi_{ij}(t)+  R_{ii}* \phi_{ij}(t),\quad t\geq 0 .\label{Resolvent002}
 \eeqlb
 Similarly, associated to the kernel $\phi_i$ we define a two-parameter function
 \beqlb
 R_{i}(t, u)\ar :=\ar \phi_i(t,u) +  R_{ii}* \phi_i(t,u) , \quad (t,u)\in \mathbb{R}_+\times \mathbb{U}, \label{Resolvent003}
 \eeqlb
 to recount the mean impact of an event with mark $u$ and its triggered events on the future arrivals of type-$i$ events.
 An argument similar to the one in \cite[Section~2]{HX2019a} deduces the next proposition immediately.  
 
 \begin{proposition}[Martingale representation]\label{SVR}
 	The intensity process $\Lambda_\mathcal{H}$ is the unique solution to the following stochastic Volterra integral equation 
 	\beqlb\label{eqn.SVR}
 	\quad\Lambda_i(t)
 	\ar=\ar  \mu_{i}(t) +  R_{ii}*\mu_i(t) +\sum_{j\in\mathcal{H}_i}  R_{ij}*\Lambda_j(t) \cr
 	\ar\ar + \int_0^t R_{iI}(s)ds +  \sum_{j\in\mathcal{D}}\int_0^t \int_\mathbb{U} R_i(t-s, u)\tilde{N}_{j}(ds,du), 
 	\quad i\in\mathcal{H}, 
 	\eeqlb
 	where $\tilde{N}_I(ds,du):=N_I(ds,du)-ds\cdot \nu_I(du)$ and $\tilde{N}_j(ds,du):=N_j(ds,du)-\Lambda_j(s-)\cdot ds\cdot \nu_j(du)$ for $j\in\mathcal{H}$.
 	Moreover, the last stochastic integral  with $j\in\mathcal{H}$ can be replaced by
 	\beqnn
 	\int_0^t \int_\mathbb{U} \int_0^{\Lambda_j(s-)} R_i(t-s, u)\tilde{N}_{0,j}(ds,du,dz)
 	\eeqnn
 	with
 	$\tilde{N}_{0,j}(ds,du,dz):=N_{0,j}(ds,du,dz)-ds\cdot \nu_j(du)\cdot dz$. 
 \end{proposition} 
 
 \subsection{Examples} \label{Sec2.3}
 In this section, we consider three specific examples, which will be revisited when analyzing scaling limits.

 \begin{example}[Exponential type]\label{ExpMHPI}  
 For $i\in\mathcal{H}$ and $j\in\mathcal{D}$, let $\beta_i >0$, $\hat{u}_i\geq 0$ and $\nu_j(du_\mathcal{H})$ be a probability law on $\mathbb{R}_+^d$. 
 The multivariate MHPI-measure $N_\mathcal{H}(dt,du_\mathcal{H})$ on $(0,\infty)\times\mathbb{R}_+^d$ is said to be of {\rm  exponential type}  with parameter $(\hat{u}_\mathcal{H},\beta_\mathcal{H},\nu_\mathcal{H},\nu_I)$ if
 \beqnn
 \mu_{i}(t)=\hat{u}_{i}e^{-\beta_i t}
 \quad \mbox{and}\quad
 \phi_i(t,u_\mathcal{H})= u_i\beta_ie^{-\beta_i t},\quad i\in\mathcal{H},\, t\geq 0,\, u_\mathcal{H} \in \mathbb{R}_+^d. 
 \eeqnn
 \end{example}

 Let $\mathcal{M}(\mathbb{R}_+)$ be the space of finite measures on $\mathbb{R}_+$ equipped with the weak convergence topology and a $\sigma$-algebra $\mathscr{M}(\mathbb{R}_+)$.
 Let $\mathcal{M}_0(\mathbb{R}_+)$ be the subspace of $\upsilon(dx) \in \mathcal{M}(\mathbb{R}_+)$ with $x\upsilon(dx) \in \mathcal{M}(\mathbb{R}_+)$ and $\mathscr{M}_0(\mathbb{R}_+)$ be the corresponding $\sigma$-algebra.
 For any $\upsilon \in \mathcal{M}_0(\mathbb{R}_+)$, denote by $L_\upsilon$ the Laplace transform of $x\upsilon(dx)$
 \beqnn
 L_\upsilon(t):= \int_0^\infty e^{-tx}x \, \upsilon(dx),\quad t\geq 0.
 \eeqnn
 By Bernstein's theorem; see Theorem~1.4 in \cite[p.3]{SSV95}, the function  $L_\upsilon$ is completely monotone on $\mathbb{R}_+$. 
 \begin{example}[Completely monotone type] 
 For $i\in\mathcal{H}$ and $j\in\mathcal{D}$, let $\hat{u}_i\in \mathcal{M}_0(\mathbb{R}_+)$ and $\nu_j(du_\mathcal{H})$ be a probability measure on $\mathcal{M}_0(\mathbb{R}_+)^d$. 
 The multivariate MHPI-measure $N_\mathcal{H}(dt,du_\mathcal{H})$ on $(0,\infty)\times \mathcal{M}_0(\mathbb{R}_+)^d$ is said to be of {\rm completely monotone type} with parameter $(\hat{u}_\mathcal{H}, \nu_\mathcal{H},\nu_I)$ if 
 \beqnn
 \mu_i(t)=L_{\hat{u}_{i}}(t)
 \quad \mbox{and}\quad
 \phi_i(t,u_\mathcal{H})= L_{u_i}(t),\quad i\in\mathcal{H},\, t\geq 0,\, u_\mathcal{H} \in \mathcal{M}_0(\mathbb{R}_+)^d.
 \eeqnn
 \end{example}

 \begin{example}[Convolution type] 
 For $i\in\mathcal{H}$ and $j\in\mathcal{D}$, let $\hat{u}_i\in \mathcal{M}(\mathbb{R}_+)$, $\rho_i$ be a non-negative, bounded, integrable function on $\mathbb{R}_+$ and $\nu_j(du_\mathcal{H})$ be a probability measure on $\mathcal{M}(\mathbb{R}_+)^d$.	
 The multivariate MHPI-measure $N_\mathcal{H}(dt,du_\mathcal{H})$ on $(0,\infty)\times \mathcal{M}(\mathbb{R}_+)^d$ is said to be of {\rm convolution type} with  parameter $(\hat{u}_\mathcal{H},\rho_\mathcal{H},\nu_\mathcal{H},\nu_I)$ if  
 \beqnn
 \mu_i(t)= \rho_i*\hat{u}_{i}(t)
 \quad\mbox{and}\quad
 \phi_i(t,u_\mathcal{H})=\rho_i*u_i(t),\quad i\in\mathcal{H},\, t\geq 0,\, u_{\mathcal{H}}\in \mathcal{M}(\mathbb{R}_+)^d.
 \eeqnn
 \end{example}

 \section{Limit theorems for self-excited dynamical systems}\label{SLMHP}
 We consider in this section the weak convergence of stochastic dynamical systems driven by nearly critical multivariate MHPI-measures, which are defined on the common filtered probability space $(\Omega,\mathscr{F},\mathscr{F}_t,\mathbf{P})$. 
 We start by presenting some basic setting on the self/mutual-excitation.
 In the $n$-th model, we assume that the MHPI-measure $N_\mathcal{H}^{(n)}(dt,du)$ has intensity process $\Lambda_\mathcal{H}^{(n)}$ and parameter\footnote{Actually, the kernel $\phi_\mathcal{H}$ is allowed to be different in various models, similarly as in footnote~\ref{Footnote.1} we also can unify them by extending the mark space.} $ (\mu_\mathcal{H}^{(n)}, \phi_\mathcal{H}, \nu_\mathcal{D}^{(n)})$.
 For $i\in\mathcal{H}$ and $j\in\mathcal{D}$, the mean impact function $\phi^{(n)}_{ij}$ is defined as (\ref{phi.ij}).
 Here we are interested in the case in which the impact of each event on the future intensity enjoys short-memory property and does not fluctuate drastically. In precise,

 \begin{enumerate} 	
 \item[\namedlabel{H1}{(H1)}]  \textit{There exist a constant    $\alpha\in(1,2)$ and a function $\Phi$ on $\mathbb{U}$ such that for any $i\in\mathcal{H}$, $j\in\mathcal{D}$ and $u\in\mathbb{U}$,
  \beqnn
 		\int_0^\infty t\cdot \phi_i(t,u)dt +   \big\|\phi_i(u) \big\|_{\rm TV}\leq \Phi(u)
 		\quad\mbox{and}\quad
 		\sup_{n\geq 1}   \int_\mathbb{U} \big| \Phi(u)\big|^{2\alpha} \nu_j^{(n)}(du) <\infty.
 \eeqnn }
 \end{enumerate}
 Actually, this is a less restrictive hypothesis, because in practice, impacts of events on the arrivals of future events usually decease fast as time goes. By using Fubini's lemma,
 \beqlb\label{ShortMemory}
 \sup_{n\geq 1} \int_0^\infty t\cdot \phi_{ij}^{(n)}(t)dt + \sup_{n\geq 1}\big\|\phi_{ij}^{(n)}\big\|_{\rm TV} <\infty.
 \eeqlb
 
 For any $K>0$, let $\mathcal{R}_K$ be the collection of non-negative functions $g$ on $\mathbb{R}_+$ whose resolvent satisfies that
 \beqnn
 R_g(t):= \sum_{k=1}^\infty g^{*k} (t) \leq K,\quad t\geq 0 . 
 \eeqnn 
 It is obvious that $\mathcal{R}_K$ comprises exponential functions in the form of $\lambda_0e^{-\lambda_1 t}$ with  $0<\lambda_0\leq (\lambda_1 \wedge K)$.
 It also contains the following two kinds of non-negative functions in the form of
 \begin{enumerate}
 	\item[$\bullet$] $f_{\lambda}*\upsilon$ in which  $\upsilon\in\mathcal{M}(\mathbb{R}_+)$ with $\upsilon(\mathbb{R}_+)\leq 1$ and $f_{\lambda}$ is the probability density function of exponential distribution with rate $\lambda\leq K$; see Lemma 4.1 in \cite{Kalashnikov1997}; \vspace{5pt}
 	
 	\item[$\bullet$] $\sum_{k=1}^\infty(-1)^{k+1}h^{*k}$ in which  $h$ is a positive, continuous, non-increasing and log-convex function with $\|h\|_{L^1_1}<\infty$, e.g., $h$ is completely monotone; see Theorem~1 in \cite{Gripenberg1978}.
 	
 \end{enumerate}
 In order to simplify the following asymptotic analysis and error estimates, we also assume an additional technical hypothesis on the mean self-excitation.
 \begin{enumerate} \it
  \item[\namedlabel{H2}{(H2)}] There exist two constants $\beta\geq 0$, $K>0$ and a non-negative function $\bar\phi$ with
  \beqnn
  \int_0^\infty t\cdot \bar\phi(t)dt<\infty
  \eeqnn
  such that for any $n\geq 1$ and $t\geq 0$,
  \beqnn
 	\phi^{(n)}_{\beta,ii}(t):= e^{-\beta t/n}\phi^{(n)}_{ii}(t) \in \mathcal{R}_K
 	\quad \mbox{and}\quad 
 	\phi^{(n)}_{\beta,ii}(t) \leq \bar\phi (t).
  \eeqnn 
 \end{enumerate}

 \subsection{Scaling limit for intensity processes}\label{Sec.SLMHP}
 We provide in this section a limit theorem for the rescaled intensity processes, which plays a key role in studying the asymptotics of self-excited dynamical systems.
 Before giving the theorem, we offer an intuitive description on how to derive it under the following asymptotic assumptions. The detailed and accurate proof can be found in Section~\ref{Section6.1}.
 \subsubsection{Asymptotic assumptions}
 By the criticality for multivariate MHPI-measures; see Proposition~\ref{CriticalityMHP},  the sequence $\{ N_\mathcal{H}^{(n)}(dt,du)  \}_{n\geq1}$ is \textit{asymptotically critical} when the matrix $\|\phi^{(n)}_{\mathcal{H}^2}\|_{L^1}$ converges to a limit matrix with spectral radius equals to one. 
 However, compared to the uni-variate case, the asymptotic criticality for multivariate Hawkes processes/measures is much more complicated because of the infinite possibilities for the limit matrix.
 In this work we mainly consider a special case in which the limit matrix is an identity matrix $\mathbf{I}$. 
 Compared to classical Hawkes processes, the random marks make additional perturbations in the convergence of rescaled intensity process via the variances of total self-excitation
 \beqnn
 c_i^{(n)}:=\int_{\mathbb{U}}\big\|\phi_{i}(u)\big\|_{L^1}^2 \nu_i^{(n)}(du) ,\quad i\in\mathcal{H},\   n\geq 1.
 \eeqnn
 We now give the detailed asymptotic assumptions on the matrix $ \phi^{(n)}_{\mathcal{H}^2} $ and vector $ \phi^{(n)}_{\mathcal{H}I} $.
 
 \begin{condition}\label{Con.PhiH} \it
 There exist a matrix $b_{\mathcal{H}^2} :=(b_{ij})_{i,j\in\mathcal{H}}$ and three vectors $ a_\mathcal{H}:= (a_i)_{i\in\mathcal{H}}\in[0,\infty)^d$, $  \sigma_\mathcal{H}:=(\sigma_i)_{i\in\mathcal{H}}\in(0,\infty)^d$, $c_\mathcal{H}:=(c_i)_{i\in\mathcal{H}} \in (0,\infty)^d$  such that  as $n\to \infty$, 
 \beqnn
 n\big(\big\|\phi^{(n)}_{\mathcal{H}^2}\big\|_{L^1}-\mathbf{I}\big)\to b_{\mathcal{H}^2}, \quad \big\|\phi_{iI}^{(n)}\big\|_{L^1}\to   a_i 
 \eeqnn
 and 
 \beqnn
 \sigma_i^{(n)}:= \int_0^\infty t \phi^{(n)}_{ii}(t)dt \to \sigma_i,\quad
 c_i^{(n)}\to c_i^2.
 \eeqnn 	
 
 \end{condition}
 
 \begin{remark} 
 Under this condition, the self-excitation in the underlying system is much stronger than the mutual-excitation, which is consistent to many practical applications. 
 For instance, the main influence of market, limit and cancel orders in financial market is on themselves, which is linked to the well known persistence of order flows and to the splitting of meta-orders into sequences of orders; see \cite{BacryJaissonMuzy2016}.
 In social media communities, it is usual that the discussion of a topic is likely to prompt further discussion as people reply to each other; see Figure 2(a,b) in \cite{PhillipsGorse2018}. 
 Moreover,  the mutual-excitation seems to be asymptotically ignorable, but it would dominate the limit process by  $b_{ij}\neq0$ with $i\neq j$.

 \end{remark}

 \begin{remark}
 It is obvious that $b_{ii}\in\mathbb{R}$ and $b_{ij}\geq 0$ for $i,j\in\mathcal{H}$ with $i\neq j$.
 Moreover, the pre-limit models are allowed to be supercritical or unstable. Indeed, if $b_{11}>0$ we have $\|\phi^{(n)}_{11}\|_{L^1}>1$ for large $n$ and
 \beqnn
 \Lambda^{(n)}_1(t)\geq \mu_1^{(n)}(t)+ \int_0^t \int_\mathbb{U} \phi_1(t-s,u)N_I(ds,du) + \int_0^t \int_\mathbb{U} \phi_1(t-s,u)N_1(ds,du).
 \eeqnn
 Similarly as in the proof of \cite[Theorem~1]{HawkesOakes1974}, we have $\mathbf{P}(\Lambda^{(n)}_1(t)\to \infty)>0$ and  $N_1(dt,du)$ is unstable.
 \end{remark}
 
 \begin{remark}
 Our asymptotic setting is different to that in Rosenbaum et al.'s works \cite{ElEuchFukasawaRosenbaum2018,JaissonRosenbaum2015}. In their setting, the kernel matrix in the $n$-th multivariate Hawkes process has the form of $\{a_n\cdot {\it \Phi}(t):t\geq 0\}$, where $\{a_n\}_{n\geq 1}$ is a positive sequence increasing to one, ${\it \Phi}(t)$ is a diagonalizable, positive matrix for each $t\geq 0$ and $\|{\it \Phi}\|_{L^1}$ has spectral radius equal to one\footnote{The matrix $\|{\it \Phi}\|_{L^1}$ is also assumed to be asymmetric in \cite{JaissonRosenbaum2015}.}.
 Because the eigenvalue of largest absolute value of $\|{\it \Phi}\|_{L^1}$ is simple and equals to one, the pre-limit model can be understood as a multivariate Hawkes process with a common intensity. This gives rise to the weak convergence of the rescaled intensity to the multiplication of one-dimensional CIR-model by a vector.

 \end{remark}

 For $i\in\mathcal{H}$ and $j\in\mathcal{D}$, the resolvent $R^{(n)}_{ij}$ and $R_i^{(n)}$ associated to the mean impact function $\phi^{(n)}_{ij}$ are defined as (\ref{Resolvent001})-(\ref{Resolvent003}), i.e.,  for any $(t,u)\in\mathbb{R}_+\times \mathbb{U}$,
 \beqlb
 R^{(n)}_{ij}(t)\ar=\ar \phi^{(n)}_{ij}(t) +  R^{(n)}_{ii}* \phi^{(n)}_{ij}(t) ,\label{ResovlentEqn01}\\
 R^{(n)}_i(t,u) \ar=\ar \phi_i(t,u) +  R^{(n)}_{ii}*\phi_i(t,u). \label{ResovlentEqn02}
 \eeqlb
 An argument similar to that in \cite{ElEuchFukasawaRosenbaum2018,JaissonRosenbaum2015} induces that the expectation $\mathbf{E}[\Lambda_\mathcal{H}^{(n)}(nt)]$ is of the order of $n$ and hence it is natural to consider the weak convergence of rescaled intensity process $\{Z^{(n)}_\mathcal{H}(t):=(Z^{(n)}_i(t))_{i\in\mathcal{H}}:t\geq 0 \}$ with $Z_i^{(n)}(t):=\Lambda_i^{(n)}(n t)/n$.
 From Proposition~\ref{SVR}, we see that $Z_\mathcal{H}^{(n)}$  satisfies the following $d$-dimensional stochastic Volterra system
 \beqlb  \label{ScaledDensity}
 Z_i^{(n)}(t)
 \ar=\ar  \frac{ \mu_i^{(n)}(nt)}{n} +   R_{ii}^{(n)}* \frac{ \mu_i^{(n)} }{n}(nt) + \sum_{j\in\mathcal{H}_i}\int_0^t nR^{(n)}_{ij}(n(t-s))Z^{(n)}_j(s)ds  \cr
 \ar\ar +\int_0^t R^{(n)}_{iI}(ns)ds + \sum_{j\in \mathcal{D}} \int_0^t \int_{\mathbb{U}}  \frac{R^{(n)}_{i}( n(t-s),u)}{n}\tilde{N}_j^{(n)}(n\cdot ds,du), 
 \quad i\in\mathcal{H},
 \eeqlb
 where $\tilde{N}_I^{(n)}(n\cdot ds,du):=N_I^{(n)}(n\cdot ds,du)-n\cdot ds\cdot \nu_I^{(n)}(du)$ and $\tilde{N}_{j}^{(n)}(n\cdot ds,du):=N_{j}^{(n)}(n\cdot ds,du)- n^2 \cdot Z^{(n)}_j(s-) \cdot ds\cdot \nu_j^{(n)}(du)$ for $j\in\mathcal{H}$. 
 
 We now give some asymptotic assumptions on the impact of events prior to time $0$ on the arrivals of future events. 
 Based on our previous argument that the mutual-excitation usually can be asymptotically ignored, it is understandable to assume that type-$i$ events prior to time $0$ make the main contribution to $ \mu^{(n)}_i$. 
 Denote by $\tau_x\leq 0$ and $u_x$ the arrival time and the mark of a typical type-$i$ event $x$ prior to time $0$.
 Because of the lack of information, we may assume it arrives uniformly before time $0$. Then its mean impact function would has the form of
 \beqnn
 I_{\phi,ii}^{(n)}(t):=\mathbf{E}\big[\phi_i(t+\tau_x,u_x)\big]
 \ar=\ar \int_{-\infty}^0 ds \int_{\mathbb{U}}\phi_i(t-s,u) \nu_i^{(n)}(du)\cr
 \ar=\ar \int_t^\infty \phi_{ii}^{(n)}(s)ds,\quad t\geq 0.
 \eeqnn 
 As the number of events prior to time $0$ goes to infinity, by the law of large numbers it is natural to assume the following condition holds. 
 Recall the constant $\alpha\in(1,2)$ in the hypothesis \ref{H1}. 
 
 \begin{condition} \label{MomentConditionInitialState} \it
 Assume that
 \beqnn
 \sup_{n\geq 1}\mathbf{E} \Big[ \big\|\mu_\mathcal{H}^{(n)}/n\big\|_{L^{1,\infty}}^{2\alpha} \Big]<\infty
 \quad \mbox{and}\quad 
 \big\|\mu^{(n)}_\mathcal{H}/n-\hat\mu^{(n)}_\mathcal{H}\big\|_{L^{1,\infty}}\overset{\rm d}\to 0, 
 \eeqnn
  as $n\to\infty$ with $\hat\mu^{(n)}_\mathcal{H}:=\big(Z_i^{(n)}(0) I_{\phi,ii}^{(n)}\big)_{i\in\mathcal{H}}$ for some random variable $Z_\mathcal{H}^{(n)}(0) \in\mathbb{R}^d_+$. 
 \end{condition}
  
 \subsubsection{Asymptotic analysis in intuition}\label{Sec.AsymptoticAna}
 We begin this section with some asymptotic analysis for the time-scaled resolvents 
 \beqnn
 \big\{R_{ij}^{(n)}(nt):t\geq 0 \big\}_{i\in\mathcal{H},j\in\mathcal{D}}
 \quad\mbox{and}\quad 
 \big\{ R^{(n)}_{i}(nt,u):t\geq 0, u\in\mathbb{U} \big\}_{i\in\mathcal{H}}.
 \eeqnn 
 From (\ref{ResovlentEqn01})-(\ref{ResovlentEqn02}), it is necessary to study $\{R_{ii}^{(n)}(n\cdot)\}_{i\in\mathcal{H}}$ first.
 Integrating both sides of (\ref{ResovlentEqn01}) over $\mathbb{R}_+$ with $j=i$, we have
 \beqnn
 \big\|R_{ii}^{(n)}\big\|_{L^1}=  \big\|\phi_{ii}^{(n)}\big\|_{L^1}+  \big\|R_{ii}^{(n)}\big\|_{L^1} \cdot  \big\|\phi_{ii}^{(n)}\big\|_{L^1}
 \eeqnn
 and hence 
 \beqnn
 \int_0^\infty R_{ii}^{(n)}(nt)
 dt= \frac{\big\|\phi_{ii}^{(n)}\big\|_{L^1}}{n\big(1- \big\|\phi_{ii}^{(n)}\big\|_{L^1}\big)} ,
 \eeqnn 
 which is finite for large $n$ if and only if 
 \beqnn
 n\Big(1- \big\|\phi_{ii}^{(n)}\big\|_{L^1}\Big)\to -b_{ii}> 0.
 \eeqnn
 Otherwise, $R_{ii}^{(n)}(nt)$ may increase to infinity.
 To overcome this difficulty, we first adjust the kernel as follows. 
 Choosing the constant $\beta$ in the hypothesis \ref{H2} larger than $\lambda_b:=\max_{j\in\mathcal{H}}b_{jj}/\sigma_j$, we define
 \beqnn
 \phi^{(n)}_{\beta,i}(t,u):=e^{-\beta t/n}\phi_{i}(t,u)
 \quad \mbox{and}\quad
 \phi^{(n)}_{\beta,ij}(t):=e^{-\beta t/n}\phi_{ij}^{(n)}(t),
 \eeqnn
 for $ (t,u)\in \mathbb{R}_+\times \mathbb{U} $, $ i\in\mathcal{H}$ and $j\in\mathcal{D}$. 
 Their resolvents $R^{(n)}_{\beta,ij}$ and $R^{(n)}_{\beta,i}$ are defined as in (\ref{ResovlentEqn01}) and (\ref{ResovlentEqn02}) respectively. The following relationships are obvious
 \beqnn
 R^{(n)}_{\beta,ij}(t)=e^{-\beta t/n}R^{(n)}_{ij}(t)
 \quad \mbox{and}\quad
 R^{(n)}_{\beta,i}(t,u)= e^{-\beta t/n}R_i^{(n)}(t,u) ,\quad (t,u)\in \mathbb{R}_+\times \mathbb{U} .
 \eeqnn
 We first consider the modified process $\{ Z_{\beta,\mathcal{H}}^{(n)}(t):t\geq 0  \}$  with 
 \beqnn
 Z_{\beta,\mathcal{H}}^{(n)}(t):= e^{-\beta t }Z_\mathcal{H}^{(n)}(t). 
 \eeqnn
 Let $  \mu^{(n)}_{\beta,\mathcal{H}}(t) :=  e^{-\beta t/n} \mu_\mathcal{H}^{(n)}(t)$ for $t\geq 0$. From (\ref{ScaledDensity}) and the foregoing notation, it is easy to identify that $Z_{\beta,\mathcal{H}}^{(n)}$ satisfies the following stochastic system
 \beqlb\label{AdjustScaledDensity01}
 Z_{\beta,i}^{(n)}(t)
 \ar=\ar\frac{ \mu^{(n)}_{\beta,i}(nt)}{n} +    R^{(n)}_{\beta,ii}*  \frac{\mu_{\beta,i}^{(n)} }{n}(nt)  + \int_0^t R^{(n)}_{\beta,iI}(n(t-s))e^{-\beta s}ds \cr
 \ar\ar + \sum_{j\in\mathcal{H}_i}\int_0^t nR^{(n)}_{\beta,ij}(n(t-s))Z^{(n)}_{\beta,j}(s)ds  \cr
 \ar\ar
 + \sum_{j\in \mathcal{D}}  \int_0^t \int_{\mathbb{U}} R_{\beta,i}^{(n)}(n(t-s),u)\frac{e^{-\beta s}}{n} \tilde{N}_j^{(n)}(n\cdot ds,du ),
 \quad i\in\mathcal{H}. 
 \eeqlb
 
 We now start to consider the convergence of the sequence $\big\{R^{(n)}_{\beta,ii}(n\cdot)\big\}_{n\geq 1}$ for each $i\in\mathcal{H}$. Notice that
 \beqnn
 n\big(1-\big\|\phi^{(n)}_{\beta,ii}\big\|_{L^1}\big)= n\big(1-\big\|\phi^{(n)}_{ii}\big\|_{L^1}\big)+ n\big(\big\|\phi^{(n)}_{ii}\big\|_{L^1}-\big\|\phi^{(n)}_{\beta,ii}\big\|_{L^1}).
 \eeqnn
 Applying the dominated convergence theorem together with the hypothesis \ref{H2} and Condition~\ref{Con.PhiH}, we have $n\big(\big\|\phi^{(n)}_{ii}\big\|_{L^1}-\big\|\phi^{(n)}_{\beta,ii}\big\|_{L^1}\big)\to \sigma_i \beta$ and hence
 \beqnn
 n\Big(1-\big\|\phi^{(n)}_{\beta,ii}\big\|_{L^1}\Big) \to  \sigma_i\beta - b_{ii}>0
 \eeqnn 
 as $n\to\infty$, which immediately induces that for large $n$,
 \beqlb\label{L1NormRii} 
 \int_0^\infty R_{\beta,ii}^{(n)}(n t)
 dt=\frac{\|\phi^{(n)}_{\beta,ii}\|_{L^1}}{n \big(1-\|\phi^{(n)}_{\beta,ii}\|_{L^1}\big)}  <\infty.
 \eeqlb
 \textit{Without loss of generality, in the sequel we will always assume that 
 	\beqnn
 	\|\phi^{(n)}_{\beta,ii}\|_{L^1} <1,\quad 
 	n\geq 1,\,i\in\mathcal{H}.
 	\eeqnn 
 }
 
 Denote by $\hat{\phi}^{(n)}_{\beta,ii}$ and $\hat{R}^{(n)}_{\beta,ii}$ the Fourier transforms of $\phi_{\beta,ii}^{(n)}$ and $R_{\beta,ii}^{(n)}$ respectively.
 Taking the Fourier transform of both sides of (\ref{ResovlentEqn01}) and then using the convolution theorem,  we have $\hat{R}^{(n)}_{\beta,ii}(\lambda)= \hat{\phi}^{(n)}_{\beta,ii}(\lambda) \big(1+\hat{R}^{(n)}_{\beta,ii}(\lambda)\big)$ for any $\lambda\in\mathbb{R}$ and hence
 \beqnn
 \int_0^\infty e^{\mathrm{i} \lambda t} R_{\beta,ii}^{(n)}(nt) dt = \frac{1}{n}\hat{R}^{(n)}_{\beta,ii}(\lambda/n)= \frac{\hat{\phi}^{(n)}_{\beta,ii}(\lambda/n)}{n\big(1-\hat{\phi}^{(n)}_{\beta,ii}(\lambda/n )\big)}.
 \eeqnn
 By the hypothesis \ref{H2} and the dominated convergence theorem,  the numerator goes to $1$ as $n\to\infty$.
 Moreover,  the dominator can be written as
 \beqnn
   n\big(1-\big\|\phi^{(n)}_{\beta,ii}\big\|_{L^1}\big) -\mathrm{i} \lambda \int_0^t t \phi^{(n)}_{\beta,ii}(t)dt  + \int_0^\infty n \big(e^{\mathrm{i} \lambda t/n }-1-\mathrm{i} \lambda t/n \big)\phi^{(n)}_{\beta,ii}(t)dt.
 \eeqnn
 By the inequality $|e^{\mathrm{i} \frac{\lambda}{n} t}-1-\mathrm{i} \frac{\lambda}{n} t| \leq  \frac{|\lambda| t}{n}\wedge \frac{|\lambda t|^2}{n^2}$ and the dominated convergence theorem, the last integral vanishes as $n\to\infty$.
 By Condition~\ref{Con.PhiH}, we have $n(1-\hat{\phi}^{(n)}_{\beta,ii}(\lambda/n))
 \to \sigma_i\beta - b_{ii} -\mathrm{i} \sigma_i \lambda$, 
 and hence
 \beqlb \label{ConvergenceFourierRii}
 \int_0^\infty e^{\mathrm{i} \lambda t} R_{\beta,ii}^{(n)}(nt) dt
 \to
 \frac{1}{\sigma_i\beta - b_{ii} -\mathrm{i} \sigma_i \lambda}
 = \int_0^\infty e^{\mathrm{i} \lambda t}
 \frac{1}{\sigma_i}e^{-(\beta - b_{ii}/\sigma_i) t } dt,
 \eeqlb
 which shows that $ R_{\beta,ii}^{(n)}(nt)  $ can be approximated by $ \sigma_i^{-1}e^{-(\beta-b_{ii}/\sigma_i) t }$.
 
 To analyze the asymptotics of the sequence $\big\{R_{\beta,i}^{(n)}(n\cdot,u) \big\}_{n\geq 1}$ for any $u\in\mathbb{U}$, we take the Fourier transform of both sides of (\ref{ResovlentEqn02}) and  obtain with some simple calculations that
 \beqnn
 \int_0^\infty e^{\mathrm{i} \lambda t} R_{\beta,i}^{(n)}(nt,u) dt =  \frac{\hat\phi_{\beta,i} (\lambda/n,u)}{ n\big(1- \hat\phi_{\beta,ii}^{(n)}(\lambda/n)\big)}, \quad \lambda \in\mathbb{R},
 \eeqnn
 where $\hat\phi_{\beta,i} (\lambda ,u)$ is the Fourier transform of $\phi_{\beta,i}(\cdot,u)$. Like the previous argument, we have as $n\to\infty$,
 \beqnn
 \int_0^\infty e^{\mathrm{i} \lambda t} R_{\beta,i}^{(n)}(nt,u) dt  \to \int_0^\infty e^{\mathrm{i} \lambda t} \frac{\|\phi_i(u)\|_{L^1}}{\sigma_i}e^{-(\beta - b_{ii}/\sigma_i) t }  dt,
 \eeqnn
 which induces that $R_{\beta,i}^{(n)}(nt,u)$ can be well approximated by  $\frac{\|\phi_i(u)\|_{L^1}}{\sigma_i}e^{-(\beta - b_{ii}/\sigma_i) t }$.
 For the rescaled resolvents $nR^{(n)}_{\beta,ij}(n\cdot)$ with  $i,j\in\mathcal{H}$ and $i\neq j$, by the dominated convergence theorem we have
 \beqnn
 \int_0^\infty e^{\mathrm{i} \lambda t} nR_{\beta,ij}^{(n)}(nt) dt
 \ar=\ar \int_{\mathbb{U}}\nu_j^{(n)}(du)\int_0^\infty e^{\mathrm{i} \lambda t} nR_{\beta,i}^{(n)}(nt,u) dt\cr
 \ar=\ar  \frac{n\hat{\phi}^{(n)}_{\beta,ij}(\lambda/n)}{n(1-\hat{\phi}^{(n)}_{\beta,ii}(\lambda/n ))}\cr
 \ar\to\ar \frac{b_{ij}}{\sigma_i\beta - b_{ii} -\mathrm{i} \sigma_i \lambda}
 = \int_0^\infty e^{\mathrm{i} \lambda t}
 \frac{b_{ij}}{\sigma_i }e^{-(\beta - b_{ii}/\sigma_i ) t } dt ,
 \eeqnn
 as $n\to\infty$ and hence it can be approximated by $\frac{b_{ij}}{\sigma_i }e^{-(\beta - b_{ii}/\sigma_i ) t }$.
 The same argument also induces that $R_{\beta,iI}^{(n)}(nt)$ can be approximated by $\frac{a_{i}}{\sigma_i }e^{-(\beta-b_{ii}/\sigma_i) t }$.

 We now turn to analyze the asymptotics of the impact of events prior to time $0$ on the future intensity.
 A simple calculation together with the hypothesis \ref{H2} shows that as $n\to\infty$,
 \beqnn
 \sup_{t\geq 0} \int_{t}^\infty 
 \big|1-e^{-\frac{\beta}{n}(s-t)} \big| \phi_{ii}^{(n)}(s)ds  \to 0,\quad i\in\mathcal{H}
 \eeqnn
 and the first two terms on the right side of (\ref{AdjustScaledDensity01}) can be approximated by $Z_i^{(n)}(0) \bar{\mu}_{\beta,i}^{(n)}(nt) $ with
 \beqnn
 \bar{\mu}_{\beta,i}^{(n)}(t):=  I^{(n)}_{\phi,\beta,ii}(t)    + R^{(n)}_{\beta,ii}*  I^{(n)}_{\phi,\beta,ii}(t)
 \quad \mbox{and}\quad
 I^{(n)}_{\phi,\beta,ii}(t):= \int_t^\infty \phi^{(n)}_{\beta,ii}(s) ds,\quad t\geq 0.
 \eeqnn
 Integrating both sides of (\ref{ResovlentEqn01}) on $[t,\infty)$ and then using Fubini's lemma, we have
 \beqnn
 \int_{t}^\infty R^{(n)}_{\beta,ii}(s) ds
 \ar=\ar   I^{(n)}_{\phi,\beta,ii}(t) +  \int_{t}^\infty R^{(n)}_{\beta,ii}(s) ds\cdot \big\|\phi^{(n)}_{\beta,ii}\big\|_{L^1}+R^{(n)}_{\beta,ii}*I^{(n)}_{\phi,\beta,ii}(t),
 \eeqnn
 which  induces that
 \beqlb\label{bar_mu}
 \bar{\mu}_{\beta,i}^{(n)}(nt) =  n\big(1-\big\|\phi^{(n)}_{\beta,ii}\big\|_{L^1}\big)   \int_t^\infty R^{(n)}_{\beta,ii}(ns)ds.
 \eeqlb
 From Condition~\ref{Con.PhiH} and (\ref{ConvergenceFourierRii}) we can approximate $ \bar{\mu}_{\beta,i}^{(n)}(nt)  $ with $ e^{ -(\beta - b_{ii}/\sigma_i)t  }$ and hence the sum of first two terms on the right side of (\ref{AdjustScaledDensity01}) is asymptotically equivalent to $Z_i^{(n)}(0) e^{ -(\beta - b_{ii}/\sigma_i)t  } $.

 Plugging all approximations above back into (\ref{AdjustScaledDensity01}),  we may have the following asymptotic equivalence for the process $ Z_{\beta,\mathcal{H}}^{(n)}$: for $i\in\mathcal{H}$,
 \beqnn
 Z_{\beta,i}^{(n)}(t)
 \ar\sim\ar Z_i^{(n)}(0)  e^{ -(\beta - b_{ii}/\sigma_i)t  }
 + \int_0^t e^{-(\beta - b_{ii}/\sigma_i) (t-s)} \cdot \frac{a_i}{\sigma_i} e^{-\beta s}ds \cr
 \ar\ar + \sum_{j\in\mathcal{H}_i}\int_0^t \frac{b_{ij}}{\sigma_i} e^{-(\beta - b_{ii}/\sigma_i) (t-s)} Z_{\beta,j}^{(n)}(s)  ds  \cr
 \ar\ar +   \sum_{j\in\mathcal{D}} \int_0^t \int_{\mathbb{U}} e^{-(\beta - b_{ii}/\sigma_i) (t-s)} \cdot \frac{\|\phi_i(u)\|_{L^1}}{n}\frac{e^{-\beta s}}{\sigma_i}\tilde{N}_j^{(n)}(n\cdot ds,du). 
 \eeqnn
 Using the fact that $e^{-(\beta-b_{ii}/\sigma_i) (t-s)}= 1-(\beta-b_{ii}/\sigma_i) \int_s^t e^{-(\beta-b_{ii}/\sigma_i) (r-s)}dr$ and Fubini's theorem, we can rewrite it into the following convenient form:
 \beqnn
 Z_{\beta,i}^{(n)}(t)
 \ar\sim \ar Z^{(n)}_i(0)  + \int_0^t  \Big( \frac{a_i}{\sigma_i}e^{-\beta s}- \beta  Z_{\beta,i}^{(n)}(s)    + \sum_{j\in\mathcal{H}}   \frac{b_{ij}}{\sigma_i} Z_{\beta,j}^{(n)}(s)\Big) ds + \sum_{j\in\mathcal{D}} M_{\beta,ij}^{(n)}(t),
 \eeqnn
 where  $M_{\beta,ij}^{(n)}$ is an $(\mathscr{F}_{nt})$-local martingale with representation
 \beqlb
 M_{\beta,ij}^{(n)}(t) \ar:=\ar \int_0^t \int_{\mathbb{U}} \frac{\|\phi_i(u)\|_{L^1}}{n}\frac{e^{-\beta s}}{\sigma_i}\tilde{N}_j^{(n)}(n\cdot ds,du),\quad t\geq 0.\label{MartMiI}
 \eeqlb
 By Condition~\ref{Con.PhiH}, for any $i\in\mathcal{H}$ and $j\in\mathcal{D}_i$ we will show the quadratic variation of $M^{(n)}_{\beta,ij}$ goes to $0$ as $n\to\infty$ and hence the sequence $\{ M^{(n)}_{\beta,ij} \}_{n\geq 1}$ converges to $0$.
 On the other hand, the quadratic variation of $M^{(n)}_{\beta,ii}$ admits the form of
 \beqnn
 [M^{(n)}_{\beta,ii}]_t= \int_0^t   \frac{c_i^{(n)}}{\sigma_i^2} e^{-\beta s} Z^{(n)}_{\beta,i}(s)ds + \int_0^t \int_{\mathbb{U}}  \frac{\|\phi_i(u)\|_{L^1}^2}{n^2}\frac{e^{-2\beta s}}{\sigma_i^2}\tilde{N}_i^{(n)}(n\cdot ds,du),\quad t\geq 0.
 \eeqnn
 Applying Doob's martingale inequality to the last stochastic integral, we see that it converges to $0$ uniformly on compacts in probability as $n\to\infty$.
 If $ Z_{\beta,\mathcal{H}}$ is a possible cluster point of the sequence $\{Z_{\beta,\mathcal{H}}^{(n)}\}_{n\geq 1}$, by Condition~\ref{Con.PhiH} we will show that
 \beqnn
 [M^{(n)}_{\beta,ii}]_t \to \int_0^t   \frac{c_i^2}{\sigma_i^2} e^{-\beta s}Z_{\beta,i}(s)ds, \quad t\geq 0, i\in\mathcal{H}.
 \eeqnn
 By Theorem III-7 in \cite{ElKarouiMeleard1990}, we can find an $(\mathscr{F}_t)$-Gaussian white noise $W_i(ds,dz)$ on $(0,\infty)^2$ with intensity $dsdz$ such that
 \beqnn
 M^{(n)}_{\beta,ii}(t) \overset{\rm d}\to \int_0^t  \int_0^{e^{\beta s}Z_{\beta,i}(s)} \frac{c_i}{\sigma_i} e^{-\beta s}W_i(ds,dz),
 \eeqnn
 in $\mathbf{D}([0,\infty),\mathbb{R})$. Additionally, the conditional orthogonality of $ \big\{\tilde{N}_{i}^{(n)}(ds,dz):i\in\mathcal{H}\big\}$ induces the mutual independence among the Gaussian white noises $\big\{ W_i(ds,dz): i\in\mathcal{H}\big\}$ .

 \subsubsection{Weak convergence of rescaled intensity processes}\label{Sec.SL}
 With all preparations above, we are ready to consider the weak convergence of the sequence $\{ Z_{\beta,\mathcal{H}}^{(n)} \}_{n\geq 1}$.   
 Letting $n\to\infty$, we may expect the limit process $Z_{\beta,\mathcal{H}}$ to be the unique  solution to the following stochastic system: for $i\in\mathcal{H}$,
 \beqnn
 Z_{\beta,i}(t)\ar=\ar  Z_i(0)+ \int_0^t \Big( \frac{a_i}{\sigma_i} e^{-\beta s}
 - \beta Z_{\beta,i}(s)+\sum_{j\in\mathcal{H}}\frac{b_{ij}}{\sigma_i} Z_{\beta,j}(s)\Big)ds\cr
 \ar\ar + \int_0^t \int_0^{e^{\beta s}Z_{\beta,i}(s)}\frac{c_i}{\sigma_i} e^{-\beta s} W_i(ds,dz).
 \eeqnn 
 By the fact that $Z_{\mathcal{H}}^{(n)}(t)=e^{\beta t}Z_{\beta,\mathcal{H}}^{(n)}(t)$ for any $t\geq 0$ and using It\^o's formula to $e^{\beta t}Z_{\beta,\mathcal{H}}(t)$, we can get the following main theorem immediately.
 \begin{theorem}\label{MainThm01}
 Under Condition~\ref{Con.PhiH} and \ref{MomentConditionInitialState},  if $Z_\mathcal{H}^{(n)}(0)\overset{\rm d}\to Z_\mathcal{H}(0) \in\mathbb{R}_+^d$, we have
 \beqnn
 Z_\mathcal{H}^{(n)} \overset{\rm d}\to Z_\mathcal{H},
 \eeqnn
  in $\mathbf{D}([0,\infty),\mathbb{R}_+^d)$ as $n\to\infty$ with the limit process $Z_\mathcal{H}  $ being the unique strong solution to 
 \beqlb \label{LimitDensity}
 Z_i(t)\ar=\ar Z_i(0)+ \int_0^t \Big(\frac{a_i}{\sigma_i} +\sum_{j\in\mathcal{H}}\frac{b_{ij}}{\sigma_i}  Z_j(s) \Big) ds  +\int_0^t \int_0^{Z_{i}(s)}\frac{c_i}{\sigma_i} W_i(ds,dz),\quad i\in\mathcal{H}.
 \eeqlb	
 \end{theorem}

 \begin{remark}
 By the martingale representation theorem in \cite[p.84]{IkedaWatanabe1989}, Theorem~\ref{MainThm01} remains valid with the stochastic integral in the limit model (\ref{LimitDensity}) replaced by
 \beqnn
 \int_0^t \frac{c_i}{\sigma_i} \sqrt{Z_{i}(s)} dB_i(s),\quad i\in\mathcal{H},
 \eeqnn
 where $B_{\mathcal{H}}$ is a standard $d$-dimensional Brownian motion.
 \end{remark}

 \begin{remark}
 By comparing (\ref{LimitDensity}) with (\ref{Intensity}) or (\ref{eqn.SVR}), we see that the limit model $Z_\mathcal{H}$ is a natural high-frequency version analogous to the intensity process of multivariate MHPI-measure.
 More precisely, we can translate the term $\int_0^\cdot \frac{b_{ii}}{\sigma_i} Z_i(s)ds$ into the net impact of type-$i$ events on themselves, i.e., as time goes, the impact of past events deceases while new impact is added.
 On the other hand, notice that $b_{ij}\geq 0$ for $j\neq i$, the term $\int_0^\cdot \frac{b_{ij}}{\sigma_i} Z_i(s)ds$ can be interpreted as the mutual-excitation of type-$j$ events on the future arrivals of type-$i$ events. 
 Clearly, in the short-memory setting, the intensity process $\Lambda_\mathcal{H}$ can be successfully recovered from the limit process $Z_\mathcal{H}$. 
 The evolution dynamic of $Z_\mathcal{H}$ is much simpler than that of $\Lambda_\mathcal{H}$. 
 Moreover, compared with the non-parametric estimation for the exogenous density and kernel in  (\ref{Intensity}), parameters in (\ref{LimitDensity}) are much easier to be estimated from the data; see \cite{Xu2014}.  
 \end{remark}

 \begin{remark}  
 From Condition~\ref{MomentConditionInitialState}, the direct impact of events prior to time $0$ in the scaling limit $Z_\mathcal{H}$ vanishes immediately after time $0$, i.e., $\mu^{(n)}_i(nt)/n \sim Z^{(n)}_i(0) I^{(n)}_{\phi,{ii}}(nt) \sim Z_i(0) \cdot \mathbf{1}_{\{t=0\}}$. 
 If the exogenous density decays slowly in the pre-limit model, then  events prior to time $0$ may continue dominating the scaling limit $Z_\mathcal{H}$ after time $0$.
 For instance, for $i\in\mathcal{H}$ and a non-negative, integrable function $g_i$ on $\mathbb{R}_+$, let $g^{(n)}_i(t):=g_i(t/n)$ for $t\geq 0$ and $n\geq 1$. If Condition~\ref{MomentConditionInitialState} holds with $\hat{\mu}^{(n)}_i = Z^{(n)}_i(0)\cdot  I^{(n)}_{\phi,{ii}}*g_i^{(n)} $, then Theorem~\ref{MainThm01} holds with  the first term on the right side of (\ref{LimitDensity}) replaced by $ Z_i(0) \int_0^t g_i(s)ds$. 
 \end{remark}

 \begin{remark}
 When the identity matrix $\mathbf{I}$ in Condition~\ref{Con.PhiH} is replaced by a diagonal matrix $ {\rm diag}(\lambda_\mathcal{H})$ with $\lambda_i\in[0,1]$, our previous asymptotic analysis remains valid with $R^{(n)}_{ij}\to 0$ for $i\in\mathcal{H}_{<1}:=\{l\in\mathcal{H}:\lambda_l<1  \}$ and $j\in\mathcal{D}$.
 Moreover, if $Z_i(0)=0$ for $i\in\mathcal{H}_{<1}$ then the weak convergence in Theorem~\ref{MainThm01} still holds with $Z_i \equiv 0$ for $i\in\mathcal{H}_{<1}$. 
 \end{remark}

 \begin{remark} 
 Jaisson and Rosenbaum \cite{JaissonRosenbaum2015} established a scaling limit for nearly unstable uni-variate Hawkes processes with a common and constant exogenous density $\mu_0$. They  approximated the rescaled intensity $Z^{(n)}$ with the solution of an It\^o's SDE with driving noise of the form  $\int_0^t|n\cdot Z^{(n)}(s-)|^{-1/2}\tilde{N}^{(n)}(ds)$, and then obtained the scaling limit by using Theorem~5.4 in \cite{KurtzProtter1991}. 
 As the key condition in \cite[Theorem~5.4]{KurtzProtter1991}, the weak convergence and uniform tightness of driving noises were identified easily by the fact that their jumps are uniformly bounded\,\footnote{Jumps of driving noise are proportional to $1/\sqrt{n\cdot Z^{(n)}}$ and uniformly bounded by $1/\sqrt{  \mu_0}$, since $Z^{(n)}\geq \mu_0 /n$ uniformly. \label{footnote_5}}; see the second paragraph in \cite[p.623]{JaissonRosenbaum2015}.
 By contrast, the exogenous intensity in our setting may vanish as time goes and the intensity process may hit $0$ in finite time\,\footnote{If $ \mu_\mathcal{H} =0 $,  we have $ \mathbf{P}( N_\mathcal{H}([0,1])+ N_I([0,1])=0)>0$ and hence $\mathbf{P}(\Lambda_i(t)=0, t\in[0,1])>0$.}. 
 In this case, driving noises may have unbounded jumps and infinite moments, which makes the proof of their weak convergence and uniform tightness difficult. 
 To get around these problems, we first approximate the rescaled intensity processes with a sequence of It\^o's SDEs driven by Poisson random measures, which are then rewritten under the form of It\^o's SDEs driven by infinite-dimensional semimartingales that have been studied in Kurtz and Protter \cite{KurtzProtter1996}.
 Notably, the infinite-dimensional semimartingales are mainly determined by a sequence of compensated Poisson random measures whose  weak convergence and uniform tightness can be identified by their orthogonal increments; see Section~\ref{ConvergenceW}.
 Finally, the weak convergence of rescaled intensity processes follows directly from \cite[Theorem 7.5]{KurtzProtter1996}; see Section~\ref{Section6.1}.
 \end{remark}

 From the argument in \cite[p.163-166]{IkedaWatanabe1989}, the unique strong solution $\big\{Z_\mathcal{H}(t):t\geq 0\big\}$ is a $d$-dimensional non-negative strong Markov process with infinitesimal generator $\mathscr{L}$ given by
 \beqnn
 \mathscr{L}f(x) := \sum_{i\in\mathcal{H}}\Big(\frac{a_i}{\sigma_i} + \sum_{j\in\mathcal{H}}\frac{b_{ij}}{\sigma_i} x_j  \Big)\frac{\partial f(x)}{\partial x_i}
 + \sum_{i\in\mathcal{H}} \frac{c_i^2}{2\sigma_i^2} x_i \frac{\partial^2 f(x)}{\partial x_i^2}, 
 \eeqnn
 for any $f\in C^2(\mathbb{R}_+^d)$.
 Here $C^2(\mathbb{R}_+^d)$ is the space of all twice differentiable functions on $\mathbb{R}_+^d$ with the first two derivatives being continuous. 
 Define a mapping $\varphi_\mathcal{H}:= (\varphi _i)_{i\in\mathcal{H}}$ from $ \mathbb{R}_+^d $ to $ \mathbb{R}^d$ by
 \beqnn
 \varphi _i(z_{\mathcal{H}}):= - \sum_{j\in\mathcal{H}} \frac{b_{ij}}{\sigma_i}  z_j + \frac{c_i^2}{2\sigma_i^2}  z_i^2,\quad z_{\mathcal{H}}\in\mathbb{R}_+^d.
 \eeqnn
 From Theorem~2.7 in \cite{DuffieFilipovicSchachermayer2003}, $Z_\mathcal{H}$ is a regular affine process with Feller transition semigroup $(Q_t)_{t\geq 0}$ on $\mathbb{R}_+^d$ defined by
  \beqnn
 \lefteqn{\int_{\mathbb{R}_+^d} e^{-\langle z_{\mathcal{H}},y_{\mathcal{H}}\rangle} Q_t(x_{\mathcal{H}},dy_{\mathcal{H}})}\ar\ar\cr
 \ar\ar = \exp \bigg\{ -\langle x_{\mathcal{H}}, v_\mathcal{H}(t,z_{\mathcal{H}})\rangle - \int_0^t \langle (a_i/\sigma_i)_{i\in\mathcal{H}}, v_\mathcal{H}(s,z_{\mathcal{H}})\rangle ds \bigg\}, 
 \eeqnn
 where  $x_{\mathcal{H}},z_{\mathcal{H}}\in\mathbb{R}_+^d$ and $ v_\mathcal{H} :=(v_i )_{i\in\mathcal{H}} $ is the unique solution to the Riccati equation
 \beqnn
 \frac{\partial}{\partial t}v_\mathcal{H}(t,z_{\mathcal{H}})=  -\varphi_\mathcal{H}( v_\mathcal{H}(t,z_{\mathcal{H}}))
 \quad \mbox{with}\quad v_\mathcal{H}(0,z_{\mathcal{H}})=z_{\mathcal{H}}  .
 \eeqnn
 Moreover, the conservative Markov process $Z_\mathcal{H}$ is also  known as a \textit{multi-type continuous-state branching process with immigration}, which has branching mechanism $\varphi_\mathcal{H}$ and immigration rate $a_\mathcal{H}$; see \cite{Watanabe1969,Xu2014}.

 \subsubsection{Examples}
 In this section, we provide scaling limits for the self-excited dynamical systems driven by multivariate MHPI-measures considered in Section~\ref{Sec2.3}. 
 For $n\geq 1$, define 
 \beqnn
 c^{(n)}_\mathcal{H}:=c^2_\mathcal{H}+\frac{1}{n}
 \quad \mbox{and}\quad
 b_{\mathcal{H}^2}^{(n)}:= \frac{b_{\mathcal{H}^2}}{n}+\mathbf{I},
 \eeqnn
 which is a positive matrix for large $n$. 
 For simplicity, we assume $b_{\mathcal{H}^2}^{(n)}$ is positive  for any $n\geq 1$.

 \begin{example}[Exponential type]\label{SLExponential}
 For each $n\geq1$, let $N^{(n)}_\mathcal{H}(dt,du_\mathcal{H})$ be a multivariate MHPI-measure of exponential type on $(0,\infty)\times\mathbb{R}_+^d$ with parameter $(\hat{u}^{(n)}_\mathcal{H}, \beta_\mathcal{H}, \nu^{(n)}_\mathcal{H},\nu_I)$ defined by: for $i\in\mathcal{H}$,
 \beqnn
 \hat{u}^{(n)}_i= Z_i(0)\cdot n,\quad \sup_{n\geq 1} \int_{\mathbb{R}_+^d} \big|u_\mathcal{H}\big|^{2\alpha} \nu_i^{(n)}(du_\mathcal{H})
 + \int_{\mathbb{R}_+^d} \big|u_\mathcal{H}\big|^{2\alpha} \nu_I (du_\mathcal{H})< \infty, \quad\\
 \int_{\mathbb{R}_+^d} u_\mathcal{H} \nu_I(du_\mathcal{H}) = a_\mathcal{H},\quad  \int_{\mathbb{R}_+^d} u_\mathcal{H} \nu_i^{(n)} (du_\mathcal{H})= b^{(n)}_{\mathcal{H}i},  \quad  \int_{\mathbb{R}_+^d} u_i^2 \nu_i^{(n)} (du_\mathcal{H}) =c_i^{(n)}.
 \eeqnn
 In this case, we have $\|\phi_i(u_\mathcal{H})\|_{L^1}= u_i$, $\|\phi_i(u_\mathcal{H})\|_{\rm TV}= u_i\beta_i$, $\int_0^\infty t\cdot\phi_i(t,u_\mathcal{H})dt=u_i/\beta_i$ and $\phi_{ij}^{(n)}(t)=b_{ij}^{(n)}\beta_{i}e^{-\beta_i t}$ for $i,j\in\mathcal{H}$.
 It is easy to identify that the two hypotheses \ref{H1}-\ref{H2} and Condition~\ref{Con.PhiH}-\ref{MomentConditionInitialState} hold. Hence Theorem~\ref{MainThm01} holds with $\sigma_\mathcal{H}= (b_{ii}/\beta_i)_{i\in\mathcal{H}}$.
 \end{example}

 \begin{remark} 
 The state $0$ is a polar for the intensity processes of multivariate MHPI-measures of exponential type, i.e., $\mathbf{P}(\Lambda^{(n)}_i(t)>0, \forall t\geq 0)=1$ for each $n\geq 1$ and $i\in \mathcal{H}$, but may be not for the limit process $Z_\mathcal{H}$, i.e. $\mathbf{P}(Z_i(t)=0, \exists t\geq 0)>0$ for some $i\in\mathcal{H}$. 
 For instance,  when $b_{ij}=0$ for $j\neq i$, then $Z_i$ is a classic CIR-model. In this case, we have $\mathbf{P}(Z_i(t)=0, \exists t\geq 0)>0$ if and only if the Feller Condition holds ($2a_i\sigma_i<c_i^2$); see \cite{Feller1951}. 
 For the general multi-type CBI-processes, several sufficient conditions are given in \cite{FriesenJinRudiger2020b} for them to not hit zero in finite time, but their polarity still remains unclear up to now.  
 \end{remark}

 Let $\chi(x):=x+1/x$ for $x>0$ and $\hat\nu(du)$ be a probability measure on $\mathcal{M}_0(\mathbb{R}_+)$ satisfying
 \beqnn
 \hat{u}(dx):= \int_{\mathcal{M}_0(\mathbb{R}_+)}\frac{ u(dx)}{x} \hat{\nu}(du)  \in \mathcal{M}(\mathbb{R}_+),
 \eeqnn
 \beqnn
 \int_{\mathcal{M}_0(\mathbb{R}_+)} |u(\chi )|^{2\alpha} \hat\nu(du)<\infty
 \quad \mbox{and}\quad
 \int_{\mathcal{M}_0(\mathbb{R}_+)} u(\mathbb{R}_+)\hat\nu(du)=1.
 \eeqnn
 Let $\hat\nu_{L}$ be a positive function on $\mathbb{R}_+$ defined by
 \beqnn
 \hat\nu_{L}(t):= \int_{\mathcal{M}_0(\mathbb{R}_+)} L_{u}(t)  \hat\nu(du) = \int_0^\infty e^{-tx} \int_{\mathcal{M}_0(\mathbb{R}_+)}  xu(dx)  \hat\nu(du),\quad t\geq 0,
 \eeqnn
 which is completely monotone with $\| \hat\nu_{L}\|_{L^1}=1$ and  $\int_0^\infty t \hat\nu_{L}(t)dt =\hat{u}(\mathbb{R}_+)<\infty$.

 \begin{example}[Completely monotone type]\label{SLCompletelyMonotone}
 For each $n\geq 1$, let $N^{(n)}_\mathcal{H}(dt,du_\mathcal{H})$ be a multivariate MHPI-measure of completely monotone type on $(0,\infty)\times \mathcal{M}_0(\mathbb{R}_+)^d$ with parameter $(\hat{u}^{(n)}_\mathcal{H},  \nu^{(n)}_\mathcal{H},\nu_I)$, where $ \hat{u}^{(n)}_\mathcal{H}(dx) = Z_\mathcal{H}(0) \cdot n\cdot \hat{u}(dx)$ and
 \beqnn
 \nu_I(du_\mathcal{H})
 \ar=\ar \int_{\mathcal{M}_0(\mathbb{R}_+)} \delta_{a_{\mathcal{H}}\cdot u}(du_\mathcal{H})\hat\nu( du),\cr
 \nu^{(n)}_j(du_\mathcal{H})\ar=\ar   \int_{\mathcal{M}_0(\mathbb{R}_+)} \delta_{b^{(n)}_{\mathcal{H}j}\cdot u}(du_\mathcal{H}) \hat\nu(du),\quad
 j\in \mathcal{H}.
 \eeqnn
 For each $i,j\in\mathcal{H}$, we have $\|\phi_i(u_\mathcal{H})\|_{\rm TV}+ \int_0^\infty t \phi_i(t,u_\mathcal{H})dt \leq C u_i(\chi)$, $\|\phi_i(u_\mathcal{H})\|_{L^1}= u_i(\mathbb{R}_+)$,
 $ \phi^{(n)}_{ij}(t)= b_{ij}^{(n)} \hat\nu_{L}(t) $, $ \|\phi_{ij}^{(n)}\|_{L^1}= b^{(n)}_{ij} $, $\|\phi^{(n)}_{iI}\|_{L^1}=a_i$ and
 \beqnn
 \mu^{(n)}_{\mathcal{H}}(t)  
 =  Z_\mathcal{H}(0)\cdot n\cdot \int_{\mathcal{M}_0(\mathbb{R}_+)} \int_0^\infty e^{-tx}u(dx)  \hat{\nu}(du)= Z_\mathcal{H}(0)\cdot n\cdot \int_t^\infty\hat\nu_{L}(s)ds. 
 \eeqnn
 These imply that the two hypotheses \ref{H1}-\ref{H2} and Condition~\ref{Con.PhiH}-\ref{MomentConditionInitialState} hold. Hence  Theorem~\ref{MainThm01} holds with
 \beqnn
 c_i^2=  \int_{\mathcal{M}_0(\mathbb{R}_+)} \big|u(\mathbb{R}_+)\big|^{2} \hat\nu(du),\quad i\in\mathcal{H}.
 \eeqnn

 \end{example}

 \begin{example}[Convolution type]\label{SLConvolution}
 For each $i\in\mathcal{H}$, let $\rho_i$ be the probability density function of exponential distribution with rate $\beta_i>0$.
 Let
 \beqnn
 \bar{\phi}_i(t)= \int_{\mathcal{M}_0(\mathbb{R}_+)} \rho_i*u(t)\hat\nu(du)
 \quad \mbox{and}\quad
 \bar{\phi}(t)= \sum_{i\in\mathcal{H}} \bar{\phi}_i(t),\quad t\geq 0.
 \eeqnn
 For each $n\geq 1$, let $N^{(n)}_\mathcal{H}(dt,du_\mathcal{H})$ be a multivariate MHPI-measure of convolution type on $(0,\infty)\times \mathcal{M}_0(\mathbb{R}_+)^d$ with parameter $(\hat{u}^{(n)}_\mathcal{H}, \rho_\mathcal{H}, \nu^{(n)}_\mathcal{H},\nu_I)$, where $\nu_I$ and $ \nu^{(n)}_\mathcal{H} $ are defined as in Example~\ref{SLCompletelyMonotone},
 \beqnn
 \hat{u}^{(n)}_i(ds)=Z_i(0)\cdot n\cdot \int_{\mathcal{M}_0(\mathbb{R}_+)} \big(\beta_i^{-1}u(\mathbb{R}_+)\delta_0(ds) + u(s,\infty)ds \big)\hat{\nu}(du) , \quad i\in\mathcal{H}.
 \eeqnn
 For $i,j\in\mathcal{H}$, we have $\|\phi_i(u_\mathcal{H})\|_{\rm TV}+ \int_0^\infty t \cdot \phi_i(t,u_\mathcal{H})dt \leq C u_i(\chi)$, $\|\phi_i(u_\mathcal{H})\|_{L^1}= u_i(\mathbb{R}_+)$, $\phi^{(n)}_{ij}(t)=b^{(n)}_{ij}  \bar{\phi}_i(t)$, $\|\phi^{(n)}_{ij}\|_{L^1}= b^{(n)}_{ij} $, $\int_0^\infty t \cdot \bar{\phi}_{i}(t)dt = \hat\nu_{L}(0) + \beta_i^{-1}$  and $\mu_i(t)= Z_i(0)\cdot n\cdot \int_t^\infty  \bar{\phi}_i(s)ds$.
 These imply that the two hypotheses \ref{H1}-\ref{H2} and Condition~\ref{Con.PhiH}-\ref{MomentConditionInitialState} hold. Thus  Theorem~\ref{MainThm01} holds with
 \beqnn
 \sigma_i=\hat\nu_{L}(0) + \beta_i^{-1}
 \quad \mbox{and}\quad
 c_i^2 =  \int_{\mathcal{M}_0(\mathbb{R}_+)} \big|u(\mathbb{R}_+)\big|^2 \hat\nu(du) ,\quad i\in\mathcal{H}.
 \eeqnn

 \end{example}

 \subsection{Scaling limits for marked Hawkes shot noise processes}\label{Section.ShotNoise}
 In this section we provide several limit theorems for shot noise processes driven by multivariate MHPI-measures, which are widely used to model the impact of events of various types on the underlying dynamical system, e.g.,  price models \cite{HX2022}, risk reserve models \cite{KluppelbergMikosch1995b},  workload input models \cite{MaulikResnick2003} and so on. 
 They also play an important role in establishing diffusion approximations for the general branching particle systems in the next section.
 In the $n$-th model, we denote by $S_\mathcal{D}^{(n)}(t):= (S^{(n)}_i(t))_{i\in\mathcal{D}}$ the total impact of all events of various types  at time $t$ with
 \beqlb\label{ShotNoiseProcess}
 S_{i}^{(n)}(t) :=\int_0^{t}\int_{\mathbb{U}} \zeta_i(t-s,u)N^{(n)}_i(ds,du),
 \eeqlb
 where $\zeta_i:\mathbb{R}_+\times \mathbb{U}\mapsto \mathbb{R}$, usually known as \textit{shape function} or \textit{response function}, is c\'adl\'ag in time and can be interpreted as the impact of each type-$i$ event on the dynamical system.  
 Specially, if $\zeta_i(t,u):= \mathbf{1}_{\{t\geq 0\}}$ for $u\in\mathbb{U}$, the shot noise process $S_\mathcal{D}^{(n)}$ reduces to the embedded point process $N_\mathcal{D}^{(n)}$. 

 As a typical application, the shot noise process (\ref{ShotNoiseProcess}) is usually considered as a natural model for the delay in claim settlement.
 Indeed, the process $S_i^{(n)}$ can be interpreted as the amount process of type-$i$ claims, in which the response function, being the form of $\{u_{i,k}(t-\tau_{i,k}):t\geq 0\}$, represents the pay-off process of the $k$-th type-$i$ claim with arrival time $\tau_{i,k}$.
 In particular, if $u_{i,k}$ is a random non-null, finite measure on $\mathbb{R}_+$  with $u_{i,k}(t)=u_{i,k}([0,t])$ for $t\geq 0$, then $S_i^{(n)}$ turns to be the total amount process of type-$i$ claims.
 Moreover, when $u_{i,k}$ is differentiable, it is usual to translate the derivative process of $S_i^{(n)}$ into the total rate at which the insurance company pays to the type-$i$ claims.
 In conclusion, here we are mainly interested in the following two kinds of response functions:
 \begin{enumerate}
 \item[$\bullet$] \textit{Cumulative response function}: $\zeta_i$ is non-negative and non-decreasing in time $t$;  \vspace{5pt}

 \item[$\bullet$] \textit{Instantaneous response function}: $\zeta_i$ is non-negative and integrable in time $t$. 
 \end{enumerate} 

 By Condition~\ref{Con.PhiH} and Theorem~\ref{MainThm01}, the arrival rates of external, mutually-triggered and self-triggered events in the $n$-th model are of the order of $1$, $1$ and $n$ respectively.
 Thus compared to that of self-triggered events, the impact of external events and mutually triggered events on the dynamical system  can be asymptotically ignored.
 Hence we mainly consider the \textit{marked Hawkes shot noise process} $S_\mathcal{H}^{(n)}$.
 Denote by $\{\zeta_{ii}^{(n)}(t):t\geq 0\} $ the \textit{mean response function} of a type-$i$ event in the $n$-th system with
 \beqnn
 \zeta_{ii}^{(n)}(t):= \int_\mathbb{U} \zeta_i(t,u)\nu_i^{(n)}(du),\quad   i\in\mathcal{H}.
 \eeqnn
 In this section we always assume that the two hypotheses \ref{H1}-\ref{H2} and Condition~\ref{Con.PhiH}-\ref{MomentConditionInitialState} hold.

 \subsubsection{Cumulative response function}
 In this section, we establish a limit theorem for the cumulative impact of events of various types on the dynamical system.
 Recall the constant $\alpha\in(1,2)$ in the hypothesis \ref{H1}.
 For each $i\in\mathcal{H}$, we assume that the total impact of a type-$i$ event with mark $u\in\mathbb{U}$  is finite, i.e., $\zeta_i(\infty,u):=\lim_{t\to\infty}\zeta_i(t,u)<\infty$, and satisfies the following condition.
 \begin{condition}\label{MomentConCumulativeEffect}
 For each $i\in\mathcal{H}$, assume that 
 \beqnn
  \sup_{n\geq 1}\int_\mathbb{U} |\zeta_i(\infty,u)|^\alpha \nu_i^{(n)}(du)<\infty
  \quad\mbox{and}\quad
  \lim_{n\to\infty}\zeta_{ii}^{(n)}(\infty) = b_{\mathtt{C}, i}\geq 0 . 
 \eeqnn 
 \end{condition}
 
 Taking expectations on both sides of (\ref{ShotNoiseProcess}) and then letting $n\to\infty$, we may have
 \beqnn
 \mathbf{E}[S_{i}^{(n)}(nt)] 
 \ar=\ar n^2
 \int_0^{t} \zeta_{ii}^{(n)}(n(t-s)) \mathbf{E}\big[Z_i^{(n)}(s) \big] ds \cr
 \ar\sim\ar n^2 \cdot \zeta_{ii}^{(n)}(\infty)\cdot \int_0^{t} \mathbf{E} \big[Z_i^{(n)}(s) \big] ds, \quad i\in\mathcal{H},
 \eeqnn
 which is  of the order of $n^2$. Thus a natural scaling in time and space leads us to consider the rescaled process $\{ S_{\mathtt{C},\mathcal{H}}^{(n)}(t):t\geq 0  \}$ with 
 \beqnn
 S_{\mathtt{C},\mathcal{H}}^{(n)}(t):= \frac{1}{n^2}\cdot S_\mathcal{H}^{(n)}(nt).
 \eeqnn 
 By the locally stochastic boundedness of the sequence $\{ Z^{(n)} \}_{n\geq 1}$, we see that the foregoing asymptotic equivalence holds if the \textit{mean residual impact} $\{\zeta_{ii}^{(n)\rm c}(t):=\zeta_{ii}^{(n)}(\infty)-\zeta_{ii}^{(n)}(t):t\geq 0,i \in\mathcal{H}\}$  satisfies the next condition.
 \begin{condition}\label{ConditionCumulativeEffect} \it
 For each $i\in\mathcal{H}$, assume that $\sup_{n\geq 1} \zeta_{ii}^{(n)\rm c}(t) \to0$ as $t\to\infty$.
 \end{condition}

 \begin{theorem}\label{ConvergenceCumulativeEffect}
 Under Condition~\ref{MomentConCumulativeEffect} and \ref{ConditionCumulativeEffect},  we have 
 \beqnn
  S^{(n)}_{\mathtt{C},\mathcal{H}} \overset{\rm d}\to S_{\mathtt{C},\mathcal{H}}
 \eeqnn
 in  $\mathbf{D}([0,\infty), \mathbb{R}_+^d)$  as $n\to\infty$ with the limit process $S_{\mathtt{C},\mathcal{H}}$ given by
 \beqnn
 S_{\mathtt{C},i}(t)
 = b_{\mathtt{C},i} \int_0^t Z_i(s)ds, \quad t\geq 0, i\in\mathcal{H}.
 \eeqnn
 \end{theorem}

 \begin{remark}
 When $\zeta_i(t,u)= \mathbf{1}_{\{t\geq 0\}}$, the shot noise process $ S_{\mathtt{C},\mathcal{H}}^{(n)}$ reduces to the rescaled embedded point process $\{  N^{(n)}_\mathcal{H}(nt)/n^2:t\geq 0\}$. 
 In this case, we have $\zeta_{ii}^{(n)}(\infty)= b_{\mathtt{C}, i}\equiv 1$ for $i\in\mathcal{H}$ and   $  N^{(n)}_\mathcal{H}(nt)/n^2 \overset{\rm d}\to \int_0^t Z_\mathcal{H}(s)ds$  in  $\mathbf{D}([0,\infty), \mathbb{R}_+^d)$  as $n\to\infty$.  
 \end{remark}

\begin{example}
For each $i\in\mathcal{H}$, let $\mathcal{P}_{\mathtt{C},i}$ be a probability law on $\mathcal{M}(\mathbb{R}_+)$ satisfying that
\beqnn
\int_{\mathcal{M}(\mathbb{R}_+)} \big|u(\mathbb{R}_+)\big|^\alpha \mathcal{P}_{\mathtt{C},i}(du) <\infty.
\eeqnn
Suppose that the pay-off processes of claims of various types in the $n$-th insurance model are mutually independent and distributed as
$\mathcal{P}_{\mathtt{C},\mathcal{H}}$.
It is easy to identify that Condition~\ref{MomentConCumulativeEffect} and \ref{ConditionCumulativeEffect} are satisfied. Then the rescaled total claim amount process converges weakly to $S_{\mathtt{C},\mathcal{H}}$ with
\beqnn
b_{\mathtt{C},i}= \int_{\mathcal{M}(\mathbb{R}_+)} u(\mathbb{R}_+) \mathcal{P}_{\mathtt{C},i}(du),\quad i\in\mathcal{H}.
\eeqnn

\end{example}

\subsubsection{Instantaneous response function}\label{Section4.1}
In this section, we consider the convergence of $\big\{S_{\mathcal{H}}^{(n)}\big\}_{n\geq1}$ with instantaneous response function that has low volatility and enjoys short-memory property, i.e.,

 \begin{condition}\label{ConInstantaneousEffects02} \it
 For each $i\in\mathcal{H}$, assume that
 \beqnn
  \sup_{n\geq 1} \int_\mathbb{U} \big( \big\|\zeta_i(u)\big\|^{2\alpha}_{\rm TV}+ \big\|\zeta_i(u)\big\|^{\alpha}_{L^1} \big) \nu_i^{(n)}(du)<\infty.
  \eeqnn
 \end{condition}

 Taking expectations on both sides of (\ref{ShotNoiseProcess}) and then letting $n\to\infty$, we may have
 \beqnn
 \mathbf{E}\big[S_{i}^{(n)}(nt)\big]
 \ar=\ar n \int_0^{nt} \zeta_{ii}^{(n)}( s) \mathbf{E}\big[Z_i^{(n)}(t-s/n)\big] ds  \cr
 \ar\sim\ar n\cdot \mathbf{E}\big[Z_i^{(n)}(t)\big]\cdot  \big\|\zeta_{ii}^{(n)}\big\|_{L^1},\quad t> 0, i\in\mathcal{H}.
 \eeqnn
 It is reasonable to consider the rescaled shot noise process
 $ \big\{S_{\mathtt{I},\mathcal{H}}^{(n)}(t) :t\geq 0\big\}$ with 
 \beqnn
 S_{\mathtt{I},\mathcal{H}}^{(n)}(t) :=  \frac{1}{n}\cdot S_\mathcal{H}^{(n)}(nt) 
 \eeqnn
 and also the mean response functions satisfying the next condition.

 \begin{condition}\label{ConInstantaneousEffects}
 For each $i\in\mathcal{H}$, there exist a non-negative, integrable function $\bar\zeta$ on $\mathbb{R}_+$ and a constant $b_{\mathtt{I},i}\geq 0$ such that for any $t\geq 0$,
 \beqnn
 \sup_{n\geq 1}\zeta_{ii}^{(n)}(t) \leq \bar\zeta(t)
 \quad \mbox{and}\quad 
 \lim_{n\to\infty}  \big\| \zeta_{ii}^{(n)} \big\|_{L^1}\to b_{\mathtt{I},i}.
 \eeqnn 
 \end{condition}

 From the intuitive analysis above, we may conjecture that  $S_{\mathtt{I},\mathcal{H}}^{(n)} $ can be approximated by $\hat{S}_{\mathtt{I},\mathcal{H}}^{(n)}$ in which
 \beqnn
 \hat{S}_{\mathtt{I},i}^{(n)}(t):=\int_0^{nt} \zeta_{ii}^{(n)}(s) Z_i^{(n)}(t-s/n) ds,\quad t\geq 0, i\in \mathcal{H}.
 \eeqnn 
 Since $Z_\mathcal{H}^{(n)} \overset{\rm d}\to Z_\mathcal{H}$ in $\mathbf{D}([0,\infty),\mathbb{R}_+^d)$; see Theorem~\ref{MainThm01},  it is natural to expect that
 \beqnn
 S_{\mathtt{I},i}^{(n)}(t)\overset{\rm d}\to S_{\mathtt{I},i}(t):= b_{\mathtt{I},i} \cdot Z_i(t),\quad t\geq 0, i\in\mathcal{H}.
 \eeqnn
 Unfortunately, this convergence may fail around time $0$, because $S_{\mathtt{I},\mathcal{H}}^{(n)}(0)\overset{\rm a.s.}=0$  but  $Z_\mathcal{H}(0)$  may be positive.

 \begin{theorem}\label{ConvergenceInstantanuousEffect}
 Under Condition~\ref{ConInstantaneousEffects02} and \ref{ConInstantaneousEffects}, we have for any $\delta>0$,
 \beqnn
 S_{\mathtt{I},\mathcal{H}}^{(n)} \overset{\rm d}\to  S_{\mathtt{I},\mathcal{H}}
 \eeqnn
 in $\mathbf{D}([\delta,\infty),\mathbb{R}_+^d)$  as $n\to\infty$.
 Moreover, if $Z_\mathcal{H}(0)\overset{\rm a.s.}=0$, this convergence also holds for $\delta=0$.
 \end{theorem}

 This convergence result fails around time $0$ mainly because the shot noise process (\ref{ShotNoiseProcess}) excludes the impact of events prior to time $0$ on the dynamical system.
 In the $n$-th model, denote by $\psi_{\mathtt{I},\mathcal{H}}^{(n)}(t)$ the total instantaneous impact of events of various types prior to time $0$ at time $t\geq 0$.
 An argument similar to that before Condition~\ref{MomentConditionInitialState} deduces that in the $n$-th model, the mean instantaneous response function of each event prior to time $0$ admits the form of
 \beqnn
 I_{\zeta,ii}^{(n)}(t):= \int_t^\infty \zeta_{ii}^{(n)}(s)ds,\quad t\geq 0, i\in\mathcal{H}.
 \eeqnn
 Applying the law of large numbers again, it is natural to assume the next condition holds for $\{\psi_{\mathtt{I},\mathcal{H}}^{(n)}\}_{n\geq 0}$.
 \begin{condition}\label{Con.InitialState01}
 Assume that $|\psi_{\mathtt{I},\mathcal{H}}^{(n)}/n-\hat{\psi}_{\mathtt{I},\mathcal{H}}^{(n)}|\overset{\rm d}\to 0$ in  $\mathbf{D}([0,\infty),\mathbb{R}^d)$ as $n\to\infty$ with $\hat{\psi}_{\mathtt{I},\mathcal{H}}^{(n)}:= (Z^{(n)}_i(0)\cdot I_{\zeta,ii}^{(n)})_{i\in\mathcal{H}}$.
 \end{condition} 

 \begin{theorem}\label{ConvergenceInstantanuousEffectAncestor}
 Under Condition~\ref{ConInstantaneousEffects02}, \ref{ConInstantaneousEffects} and \ref{Con.InitialState01},
 we have  $  \psi_{\mathtt{I},\mathcal{H} }^{(n)}(n\cdot)/n +S_{\mathtt{I},\mathcal{H} }^{(n)} \overset{\rm d}\to   S_{\mathtt{I},\mathcal{H}}  $ in  $\mathbf{D}([0,\infty),\mathbb{R}_+^d)$ as $n\to\infty$.
 \end{theorem}

 Let $L^1_{\rm TV}(\mathbb{R}_+)$ be the space of non-negative, integrable and c\'adl\'ag functions on $\mathbb{R}_+$ with bounded variation. 
 It is endowed with the norm $\|\cdot\|_{\rm TV}+ \|\cdot\|_{L^1}$. 
 For each $i\in\mathcal{H}$, let $\mathcal{P}_{\mathtt{I},i}$ be a probability measure on $L^1_{\rm TV}(\mathbb{R}_+)$ satisfying that
 \beqnn
 \int_{L^1_{\rm TV}(\mathbb{R}_+)} \big(\big\|u\big\|_{\rm TV}^{2\alpha}+ \big\|u\big\|_{L^1}^\alpha\big) \mathcal{P}_{\mathtt{I},i}(du)<\infty.
 \eeqnn

 \begin{example}
 In the $n$-th insurance model, suppose that there are $[Z_\mathcal{H}(0)\cdot n]$ claims at time $0$ and the pay-off rate of claims of various types is distributed as $\mathcal{P}_{\mathtt{I},\mathcal{H}}$.
 Here we are interested in the total rate at which the insurance company pays to claims of various types.
 It is described as $S^{(n)}_\mathcal{H}$ with $\zeta_i(t,u)=u(t)$ for $u\in L^1_{\rm TV}(\mathbb{R}_+)$ and $t\geq 0$.
 In this case, we see that Condition~\ref{ConInstantaneousEffects02} and \ref{ConInstantaneousEffects} hold  with $\zeta_{ii}^{(n)}(t)= \int_{L^1_{\rm TV}(\mathbb{R}_+)} u(t) \mathcal{P}_{\mathtt{I},i}(du)$.
 Additionally, the pay-off rate of a typical claim $x$ prior to time $0$ with arrival time $\tau_x<0$ is $ u_x(t-\tau_x) $ at time $t\geq 0$.  
 By the law of large numbers and an argument similar to that before Condition~\ref{MomentConditionInitialState}, the total pay-off rate process of type-$i$ claims prior to time $0$ can be approximated by
 \beqnn
 Z_i(0)\int_t^\infty \int_{L^1_{\rm TV}(\mathbb{R}_+)} u(s) \mathcal{P}_{\mathtt{I},i}(du) = Z_i(0)I_{\zeta,ii}^{(n)}(t),\quad t\geq 0.
 \eeqnn
 Hence Condition~\ref{Con.InitialState01} holds and the rescaled total pay-off rate process converges weakly to $ S_{\mathtt{I},\mathcal{H}} $ with
 \beqnn
 b_{\mathtt{I},i} = \int_{L^1_{\rm TV}(\mathbb{R}_+)}\|u\|_{L^1} \mathcal{P}_{\mathtt{I},i}(du),\quad i\in\mathcal{H}.
 \eeqnn
 \end{example}

 \section{Limits theorems for multi-type CMJI-processes}\label{Section.CMJ}
 
 In this section, we apply our limit theorems for self-excited dynamical systems to establish diffusion approximations for \textit{multi-type Crump-Mode-Jagers branching processes} (CMJI-processes).
 In order to clarify the connection between multi-type CMJI-processes and multivariate MHPI-measures, we try to use the same notation to represent quantities of population that play the similar roles in MHPI-measures.  
 
 \subsection{Multi-type CMJI-processes} A $d$-type CMJI-process is a general branching process with $d$ kinds of distinguishable individuals. 
 These are usually to be called type-$1,2,\cdots,d$.
 In order to illustrate its considerable importance in biology, we give its definition with budding microbes as a typical example.
 It is usual to assume that the observation on the population starts from the appearance of symptoms on the host. The microbes alive at time $0$ are considered as ancestors. 
 \begin{enumerate}
 	\item[\namedlabel{P1}{(P1)}] ({\it Ancestors}) The population starts with $\Xi_\mathcal{H}(0) := \big(\Xi_i(0)\big)_{i\in\mathcal{H}} \in \mathbb{N}^d$ ancestors at time $0$. 
 \end{enumerate}
 
 Compared to binary fission microbes, budding microbes usually live much longer before dying,  being killed or spreading out of the host.
 Moreover, their life-lengths are rarely exponentially distributed; see \cite[Table~4]{HolbrookMenninger2002} and \cite[Figure 2-4]{WoodRogina2004}. 
 \begin{enumerate}
 	\item[\namedlabel{P2}{(P2)}] ({\it Life-length}) Individuals of type-$i$ have a common life-length distribution $\mathcal{P}_{\mathtt{L},i}(dy)$ on $\mathbb{R}_+$ with finite first and second moments 
 	\beqnn
 	\mathrm{m}_{\mathtt{L},i}:= \int_0^\infty y \mathcal{P}_{\mathtt{L},i}(dy).
 	\quad \mbox{and}\quad
 	\mathrm{v}_{\mathtt{L},i}:= \int_0^\infty y^2 \mathcal{P}_{\mathtt{L},i}(dy) .
 	\eeqnn  
 \end{enumerate}
 
 Moreover, different to binary fission in which the fully grown parent cell either dies or splits into equally sized daughter cells, the mother budding microbes usually produce buds several times during their lifetime.
 As so often, the budding rate is low during the growth stage and then increases to the highest level after separating from the mother cell.
 As the bud scars accumulate on the surface, the microbe enters into the senescence state and the budding rate starts to decrease; see \cite[Figure 2]{JiangJaruga2000}.
 We collect the possible budding rate functions in
 \beqnn
 \mathbb{B}:=\Big\{\mathtt{B}: \mathbb{R}_+^2\mapsto \mathbb{R}_+: \mathtt{B}(t,y)=0  \mbox{ if }t\geq y \mbox{ and } \big\|\mathtt{B}(y)\big\|_{\rm TV}+    \int_0^\infty t \cdot \mathtt{B}(t,y)dt <\infty  \Big\}
 \eeqnn
 and describe the reproduction process of each budding microbe by a Cox process with intensity process selected randomly in $\mathbb{B}$; see the following properties.
 \begin{enumerate}
 	\item[\namedlabel{P3}{(P3)}] ({\it Budding rate}) Each type-$i$ individual is endowed with a budding rate function randomly according to the probability law $\mathcal{P}_{\mathtt{B},i}(d \mathtt{B})$  on  $\mathbb{B}$. The mean budding rate is finite and light-tailed, i.e.,
 	\beqnn
 	\mathtt{B}_i(t) := \int_0^\infty \mathcal{P}_{\mathtt{L},i} (dy) \int_\mathbb{B}  \mathtt{B}(t,y)  \mathcal{P}_{\mathtt{B},i} (d\mathtt{B})<\infty
 	\quad \mbox{and}\quad
 	\mathrm{d}_{\mathtt{B},i} := \int_0^\infty t \cdot \mathtt{B}_i(t) dt<\infty .
 	\eeqnn
 	\vspace{3pt}

 \item[\namedlabel{P4}{(P4)}] ({\it Successive ages}) Conditioned on the life-length $y$ and budding rate function $\mathtt{B}$, the successive ages $0<t_1<t_2<\cdots<y$ at which the individual gives birth to offspring are described by an in-homogeneous Poisson  process on $(0,y)$ with intensity $\mathtt{B}(\cdot,y)$ and the mean number of successive ages is $ \big\|\mathtt{B}(y)\big\|_{L^1}:=\int_0^\infty  \mathtt{B}(t,y)dt$.
 The first and second moments of successive ages are finite, i.e., for $i\in\mathcal{H}$,
 \beqnn
  \mathrm{m}_{\mathtt{B},i} 
  \ar:=\ar \int_0^\infty \mathcal{P}_{\mathtt{L},i} (dy) \int_\mathbb{B}  \big\|\mathtt{B}(y)\big\|_{L^1}  \mathcal{P}_{\mathtt{B},i} (d\mathtt{B}),\cr
  \mathrm{v}_{\mathtt{B},i}
  \ar:=\ar \int_0^\infty \mathcal{P}_{\mathtt{L},i}(dy) \int_\mathbb{B}  \big\|\mathtt{B}(y)\big\|_{L^1}^2  \mathcal{P}_{\mathtt{B},i} (d\mathtt{B}).
 \eeqnn

 \end{enumerate}

 Usually, only one bud forms on the mother cell at each successive age.
 But multiple-budding is also widely observed in enveloped virus such as HIV and COVID-19; see \cite[p.384]{TortoraFunke2018}. 
 
 \begin{enumerate}
 	\item[\namedlabel{P5}{(P5)}] ({\it Branching mechanism}) At each successive age, a type-$i$ individual gives birth to a random number of offspring of various types according to a probability law $p_i:=\{ p_i(k_\mathcal{H}): k_\mathcal{H} \in\mathbb{N}^d \}$, where $p_i( k_\mathcal{H})$ is the probability to produce $k_1$ children of type-$1$, $k_2$ of type-$2$, ..., $k_d$ of type-$d$. The first and second moments of offspring of various types are finite
 	\beqnn
 	\mathrm{m}_{ij} := \sum_{k_\mathcal{H}\in\mathbb{N}^d} k_i \cdot p_j(k_\mathcal{H}) \quad\mbox{and}\quad
 	\mathrm{v}_{ij} := \sum_{k_\mathcal{H}\in\mathbb{N}^d} k_i^2 \cdot p_j (k_\mathcal{H}) , \quad j\in\mathcal{H}.
 	\eeqnn
 \end{enumerate}

 In addition to budding, microbes may enter into the host from the external environment or the neighboring hosts.
 For simplicity, we assume that 
 \begin{enumerate} 	
 	\item[\namedlabel{P6}{(P6)}] ({\it Immigration rate}) The arrivals of immigrants follow a Poisson point process with a unit rate ; \vspace{3pt}
 	
 	\item[\namedlabel{P7}{(P7)}] ({\it Immigration mechanism}) The number of invading microbes of various types in each immigration is distributed as a probability law $p_I:=\{p_I(k_\mathcal{H}) : k_\mathcal{H} \in \mathbb{N}^d \}$, where $p_I(k_\mathcal{H})$ is the probability that $k_1$ immigrants of type-$1$, $k_2$ of type-$2$, ..., $k_d$ of type-$d$ enter into the population. The mean number of immigrants is finite
 	\beqnn
 	\mathrm{m}_{iI} := \sum_{k_\mathcal{H}\in\mathbb{N}^d} k_i\cdot p_I(k_\mathcal{H}),\quad i\in\mathcal{H}.
 	\eeqnn
 	
 \end{enumerate}
 
 It is the impact of microbes on the host that has been widely considered in mathematical biology literature, e.g., releasing toxins and attacking the host cell.
 For instance, \textit{Candida albicans} in the gastrointestinal and genitourinary tract do not only release a  kind of toxins called \textit{Candidiasis} but also \textit{alkalinize phagosome} by physical rupture.
 We refer the impact of each microbe on the host as its \textit{characteristic}, which  usually is described as a non-negative function of its age and life-length.
 Denote by $\mathbb{T}$ the measurable space of all possible characteristic functions on $\mathbb{R}_+^2$. 
 
 \begin{enumerate}
 	\item[\namedlabel{P8}{(P8)}] {\it(Characteristic)} Each type-$i$ individual is endowed with a characteristic function randomly according to a probability law $\mathcal{P}_{\mathtt{T},i}(d\mathtt{T})$ on $\mathbb{T}$ with
 	\beqnn
 	\mathtt{T}_i(t):= \int_0^\infty \mathcal{P}_{\mathtt{L},i}(dy) \int_\mathbb{T} \mathtt{T}(t,y) \mathcal{P}_{\mathtt{T},i}(d\mathtt{T})<\infty,\quad t\geq 0.
 	\eeqnn
 	
 \end{enumerate}
 
 The branching particle system defined by these properties is a multi-type CMJI-process with initial state $\Xi_\mathcal{H}(0)$ and parameter $(p_{\mathcal{D}},\mathcal{P}_{\mathtt{L},\mathcal{H}}, \mathcal{P}_{\mathtt{B},\mathcal{H}}, \mathcal{P}_{\mathtt{T},\mathcal{H}})$. 
 In particular, when $p_{\mathcal{D}}(1)=1$ and $\mathcal{P}_{\mathtt{B},\mathcal{H}}(B(t,y)=\mathbf{1}_{\{y>t\}})=1$, it is often  known as a \textit{homogeneous, binary CMJI-process}.
 Denote by $\mathcal{I}_i$ the collection of all type-$i$ individuals in the population. 
 Associated with each individual $x\in\mathcal{I}_i$ is a quadruple $(\tau_x, \ell_x, \mathtt{B}_x, \mathtt{T}_x)$ that represents its birth/immigrating time, life-length, budding rate function and characteristic function.
 We are usually interested in the \textit{multi-type CMJI-process counted with random characteristic $\mathtt{T}$} ($\mathtt{T}$-CMJI-process), denoted by  $\{\mathbf{T}_\mathcal{H}(t):=(\mathbf{T}_i(t))_{i\in\mathcal{H}}:t\geq 0\}$ with 
 \beqnn
 \mathbf{T}_i(t):= \sum_{x\in\mathcal{I}_i} \mathtt{T}_x(t-\tau_x,\ell_x).
 \eeqnn
 Specially, if $\mathcal{P}_{\mathtt{T},i}(\mathtt{T}(t,y)=\mathbf{1}_{\{y>t\}})=1$, then the $\mathtt{T}$-CMJI-process reduces to the process of \textit{population size}, denoted as $\{\Xi_\mathcal{H}(t):t\geq 0\}$ with
 \beqnn
 \Xi_i(t):= \sum_{x\in\mathcal{I}_i} \mathbf{1}_{\{\ell_x>t-\tau_x\geq 0\}}, \quad t\geq 0,i \in \mathcal{H}.
 \eeqnn
 We end this section with several typical characteristic functions that are widely considered in biology and mathematics.  
 
 \begin{example}
 	For each $i\in\mathcal{H}$, if $\mathcal{P}_{\mathtt{T},i}$ is a probability measure on $\mathcal{M}(\mathbb{R}_+)$ and $\mathtt{T}(t,y):=\mathtt{T}(t\wedge y)$ is the mass of $\mathtt{T}$ on $[0,t\wedge y]$, then $\mathbf{T}_i(t)=\sum_{x\in\mathcal{I}_i} \mathtt{T}_x((t-\tau_x)\wedge \ell_x) $  is a {\rm multi-type CMJI-process counted with random measure}.
 	In particular,
 	\begin{enumerate}
 		\item[(1)] If $ \mathcal{P}_{\mathtt{T},i}(\mathtt{T}(t,y)=\mathbf{1}_{\{t\geq 0\}}) =1$, then
 		$ \mathbf{T}_i(t)=\sum_{x\in\mathcal{I}_i}\mathbf{1}_{\{t\geq \tau_x\}} $ is known as the {\rm total  progeny of type-$i$} up to time $t$; \vspace{3pt}
 		
 		\item[(2)]  If $\mathcal{P}_{\mathtt{T},i}(\mathtt{T}(t,y)=t^+\wedge y)=1$, then $\mathbf{T}_i$ is the {\rm integral of type-$i$ population}, i.e.,
 		\beqnn
 		\mathbf{T}_i(t)=\sum_{x\in\mathcal{I}_i} (t-\tau_x)^+ \wedge \ell_x =\sum_{x\in\mathcal{I}_i} \int_0^t  \mathbf{1}_{\{\ell_x>s-\tau_x\geq0\}}ds = \int_0^t  \Xi_i(s) ds, \quad t\geq 0.
 		\eeqnn
 		
 	\end{enumerate}
 \end{example}
 
 \begin{example}
 	For each $i\in\mathcal{H}$, let $\mathcal{P}_{\mathtt{T},i}$ be a probability measure on $L^1_{\rm TV}(\mathbb{R}_+)$.
 	Then $\mathbf{T}_i(t)=\sum_{x\in\mathcal{I}_i} \mathtt{T}_x(t-\tau_x)\mathbf{1}_{\{t-\tau_x\leq  \ell_x\}}$ is a {\rm multi-type CMJI-process counted with random integrable function}.
 	In particular, for some constant $\eta>0$,
 	
 	\begin{enumerate}
 		\item[(1)] If $ \mathcal{P}_{\mathtt{T},i}(\mathtt{T}(t,y)=\mathbf{1}_{\{0\leq t<\eta \wedge y \}})=1 $, then
 		$
 		\mathbf{T}_i(t)=\sum_{x\in\mathcal{I}_i} \mathbf{1}_{\{t-\tau_x \in[0,\eta\wedge \ell_x )\}}
 		$
 		is the total type-$i$ population alive at time $t$ which is younger than $\eta$;  \vspace{5pt}

 		\item[(2)] If $ \mathcal{P}_{\mathtt{T},i}(\mathtt{T}(t,y)=\mathbf{1}_{\{\eta\leq t<  y \}})=1 $, then
 		$
 		\mathbf{T}_i(t)=\sum_{x\in\mathcal{I}_i} \mathbf{1}_{\{\eta \leq t-\tau_x <\ell_x \}}
 		$
 		is the  total type-$i$ population alive at time $t$ which is older than $\eta$;  \vspace{5pt}
 		
 		\item[(3)] If $\mathcal{P}_{\mathtt{T},i}( \mathtt{T}(t,y)= \mathbf{1}_{\{0< y-t\leq \eta \}})=1$, then
 		$
 		\mathbf{T}_i(t)=\sum_{x\in\mathcal{I}_i} \mathbf{1}_{\{ 0<\ell_x-(t-\tau_x) \leq \eta \}}
 		$
 		is the  total type-$i$ population alive at time $t$ with residual life less than $\eta$.
 		
 	\end{enumerate}
 \end{example}

 \subsection{Hawkes representation}
 In this section, we link the foregoing multi-type CMJI-process to a self-excited dynamical system driven by multivariate MHPI-measures.
 Different to the early literature in which the population size is often studied first, we start by considering the \textit{total budding rate process} $\{ \mathbf{B}_\mathcal{H}(t):=(\mathbf{B}_i(t))_{i\in\mathcal{H}}:t\geq 0 \}$, where $\mathbf{B}_i(t)$ is the total budding rate of all type-$i$ individuals alive at time $t$, i.e.,
 \beqnn
 \mathbf{B}_i(t)=\sum_{x\in\mathcal{I}_i} \mathtt{B}_x(t-\tau_x, \ell_x).
 \eeqnn 
 Obviously, the process $ \mathbf{B}_\mathcal{H}$ is a $\mathtt{B}$-CMJI-process.
 However, this representation fails to clarify the population evolution dynamics and is not helpful to explore the long-term behavior of the population.
 We now establish a new representation based on a finer classification for individuals of various types. 
 Denote by $\{\tau_{I,k}\}_{k\geq 1}$ the immigrating times.
 For $i\in\mathcal{H}$, let $\{\tau_{i,k} \}_{k\geq 1}$ be the successive ages of all type-$i$ individuals. 
 From  property \ref{P4} and the mutual independence among individuals, we have $\tau_{i,k} < \tau_{i,k+1}$ and $\tau_{i,k}\neq \tau_{j,l}$ a.s. for any $(i,k),(j,l) \in \mathcal{D}\times \mathbb{Z}_+$ with $(i,k)\neq (j,l)$. 
 According to the origin of each type-$i$ individual, we can split $\mathcal{I}_i$ into three kinds of disjoint sets: for each $ j\in\mathcal{H}$ and $k\geq 1$,
 \begin{enumerate}
 	\item[$\bullet$] $\mathcal{I}_{i,0}$: Ancestors of type-$i$ at time $0$; \vspace{5pt}
 	
 	\item[$\bullet$] $\mathcal{I}_{i,I,k}$: Immigrants of type-$i$ entering into the population at the immigrating time $\tau_{I,k}$; \vspace{5pt}
 	
 	\item[$\bullet$] $\mathcal{I}_{i,j,k}$:  Offspring of type-$i$ produced by a type-$j$ mother individual at the successive age $\tau_{j,k}$.
 \end{enumerate}
 Notice that for each individual $x$ in $\mathcal{I}_{i,0}$ or $\mathcal{I}_{i,j,k}$ with $i\in\mathcal{H}$, $j\in\mathcal{D}$ and $k\geq 1$ , we have $\tau_x=0$ or $\tau_{j,k}$ respectively.
 Thus we can write the total budding rate of all type-$i$ individuals alive at time $t$ as
 \beqlb\label{eqn.HPB}
 \mathbf{B}_i(t)\ar=\ar \sum_{x\in\mathcal{I}_{i,0}} \mathtt{B}_{x}(t-\tau_x,\ell_{x}) + \sum_{\tau_{I,k}\leq t}\,\sum_{x\in \mathcal{I}_{i,I,k}} \mathtt{B}_x(t-\tau_{I,k}, \ell_x) \cr
 \ar\ar + \sum_{j\in\mathcal{H}}\, \sum_{\tau_{j,k}\leq t}\,\sum_{x\in \mathcal{I}_{i,j,k}} \mathtt{B}_x(t-\tau_{j,k}, \ell_x).
 \eeqlb
 Here the first sum on the right side of this equation is the total budding rate of all type-$i$ ancestors.
 The inner-sum in the second term is the total budding rate of all type-$i$ immigrants entering into the population in the $k$-th immigration.
 Similarly, the second inner-sum in the last term is the total budding rate of all type-$i$ offspring born at time $\tau_{i,k}$.
 Repeating the previous progress, we also can give representation analogous to (\ref{eqn.HPB}) for $\mathbf{T}_i(t)$ by replacing the budding rate function $\mathtt{B}_x$ with the characteristic function $\mathtt{T}_x$.
 
 To get a Hawkes representation for the CMJI-process, it remains to construct two random point measures to describe the arrivals and characteristics of immigration and reproduction respectively.
 Let $\mathbb{U}:= (\mathbb{N}\times \mathbb{R}_+^{\mathbb{N}}\times \mathbb{B}^\mathbb{N}\times \mathbb{T}^\mathbb{N})^d $.
 For each $j\in\mathcal{D}$ and $k\geq 1$, we introduce a notation $$\boldsymbol{u}_{j,k}:=(k_{\mathcal{H}},y_{\mathcal{H}},\mathrm{B}_{\mathcal{H}},\mathrm{T}_{\mathcal{H}})
 :=((k_{i})_{i\in\mathcal{H}},(y_{i})_{i\in\mathcal{H}},(\mathrm{B}_{i})_{i\in\mathcal{H}} ,(\mathrm{T}_{i})_{i\in\mathcal{H}} ) \in\mathbb{U}$$
 to describe the  new individuals getting into the population at time $\tau_{j,k}$, where
 \begin{enumerate}
 	\item[$\bullet$] $k_{i} \in\mathbb{N}$: Number of type-$i$ offspring/immigrants; \vspace{5pt}
 	
 	\item[$\bullet$] $y_{i}:=(\mathtt{y}_{i,l})_{l=1,\cdots, k_i}\in \mathbb{R}_+^{k_i}$: Life-lengths of type-$i$ offspring/immigrants; \vspace{5pt}
 	
 	\item[$\bullet$] $\mathrm{B}_{i}:=(\mathtt{B}_{i,l})_{l=1,\cdots, k_i}\in \mathbb{B}^{k_i}$: Budding rate functions of type-$i$ offspring/immigrants;\vspace{5pt}
 	
 	\item[$\bullet$] $\mathrm{T}_{i}:=(\mathtt{T}_{i,l})_{l=1,\cdots, k_i} \in \mathbb{T}^{k_i}$: Characteristic functions of type-$i$ offspring/immigrants.
 \end{enumerate} 
 At time $t$, the total budding rate and total characteristic of these new born/immigrating individuals of type-$i$ can be written as
 \beqlb
 \sum_{x\in \mathcal{I}_{i,j,k}} \mathtt{B}_x(t-\tau_{j,k}, \ell_x)\ar=\ar \sum_{l=1}^{k_{i}} \mathtt{B}_{i,l}(t-\tau_{j,k},\mathtt{y}_{i,l})=:\phi_i(t-\tau_{j,k}, \boldsymbol{u}_{j,k}), \label{Kernel01} \\
 \sum_{x\in \mathcal{I}_{i,j,k}} \mathtt{T}_x(t-\tau_{j,k}, \ell_x)\ar=\ar \sum_{l=1}^{k_{i}} \mathtt{T}_{i,l}(t-\tau_{j,k}, \mathtt{y}_{i,l})=:  \zeta_i(t-\tau_{j,k}, \boldsymbol{u}_{j,k}). \label{Kernel02}
 \eeqlb
 For each $j\in\mathcal{D}$, associated to the sequence $\{ (\tau_{j,k},\boldsymbol{u}_{j,k})  \}_{k\geq 1}$ we define an $(\mathscr{F}_t)$-random point measure on $(0,\infty)\times \mathbb{U}$
 \beqnn
 N_j(ds,d\boldsymbol{u})= \sum_{k=1}^\infty \mathbf{1}_{\{ \tau_{j,k}\in ds, \boldsymbol{u}_{j,k} \in d\boldsymbol{u} \}}.
 \eeqnn
 From properties \ref{P3}-\ref{P4} and \ref{P6}-\ref{P7}, we see that $N_j(ds,d\boldsymbol{u})$ has intensity  $ds \cdot \nu_I(d\boldsymbol{u})$ when $j=I$ or $\mathbf{B}_j(s-)\cdot ds\cdot \nu_j(d\boldsymbol{u})$ when $j\in\mathcal{H}$, where $\nu_j(d\boldsymbol{u})$ is a probability law on $\mathbb{U}$ defined by
 \beqnn
 \nu_j(d\boldsymbol{u}):= \sum_{n_\mathcal{H}\in\mathbb{N}^d}  p_j(n_\mathcal{H}) \cdot \delta_{n_\mathcal{H}} (dk_\mathcal{H})
 \prod_{i\in\mathcal{H}}\prod_{l=1}^{n_i} \mathcal{P}_{\mathtt{L},i}(d\mathtt{y}_{i,l}) \mathcal{P}_{\mathtt{B},i} ( d\mathtt{B}_{i,l}) \mathcal{P}_{\mathtt{T},i} (d\mathtt{T}_{i,l}),\quad j\in \mathcal{D}.
 \eeqnn
 
 We now give a more detailed description for ancestors. 
 For each $i\in\mathcal{H}$ and ancestor $x\in \mathcal{I}_{i,0}$, its life-length equals to the sum of its age $\mathtt{A}_{i,x}$ and residual life $\mathtt{R}_{i,x}$ at time $0$. 
 Thus we can write the total budding rate and total characteristic of all type-$i$ ancestors at time $t\geq 0$ as
 \beqlb 
 \mu_{i}(t)\ar:=\ar \sum_{x\in \mathcal{I}_{i,0}}\mathtt{B}_{x}(t+\mathtt{A}_{i,x}, \mathtt{R}_{i,x}+\mathtt{A}_{i,x} ),\label{InitialStateCMJ01} \\
 \psi_i(t)\ar:=\ar \sum_{x\in \mathcal{I}_{i,0}}\mathtt{T}_{x}(t+\mathtt{A}_{i,x}, \mathtt{R}_{i,x}+\mathtt{A}_{i,x} ). \label{InitialStateCMJ}
 \eeqlb
 Based on all preparations above, an argument similar to that in Section~\ref{BrachingRep} gives the following Hawkes representations for the total budding rate process $\mathbf{B}_\mathcal{H}$ and the $\mathtt{T}$-CMJI-process $ \mathbf{T}_\mathcal{H} $. 
 \begin{proposition}
 	$ N_\mathcal{H}(ds,d\boldsymbol{u})$ is a multivariate MHPI-random measure on $\mathbb{R}_+\times \mathbb{U}$ with mark distribution $\nu_\mathcal{H}(d\boldsymbol{u})$ and intensity process $\mathbf{B}_\mathcal{H}$ admitting the form of
 	\beqnn
 	\mathbf{B}_i(t)\ar=\ar \mu_{i}(t)+
 	\sum_{j\in\mathcal{D}} \int_0^t \int_\mathbb{U}  \phi_i(t-s,\boldsymbol{u}) N_j(ds,d\boldsymbol{u}),\quad t\geq 0,  i\in\mathcal{H}.\label{TotalBirthRate}
 	\eeqnn
 	Moreover, an analogous representation for $\mathbf{T}_\mathcal{H}$ can be obtained by replacing $(\mu_\mathcal{H},\phi_\mathcal{H})$ with  $(\psi_\mathcal{H},\zeta_\mathcal{H})$.
 \end{proposition}

 \subsection{Scaling limit theorems}
 
 In practice, the microbial population usually is very large and birth/death events occur at a high-frequency. 
 These make the low-frequency biological models (e.g. CMJ-processes and GW-processes) inefficient and the high-frequency models popular in the study of microbial population. 
 We introduce a parameter $n\in\mathbb{Z}_+$ to scale the population size   and assume that individuals are weighted by $1/n$. 
 Under some mild scaling assumptions, we now establish several limit theorems for  multi-type CMJI-processes by using the convergence results in Section~\ref{SLMHP}. 
 In the $n$-th model, the multi-type CMJI-process starts from $\Xi_{\mathcal{H}}^{(n)}(0)$ ancestors and has parameter
 $( p_{\mathcal{D}}^{(n)}, \mathcal{P}^{(n)}_{\mathtt{L},\mathcal{H}}, \mathcal{P}^{(n)}_{\mathtt{B},\mathcal{H}},\mathcal{P}^{(n)}_{\mathtt{T},\mathcal{H}})$.
 Quantities in the last two sections are defined similarly with superscript $(n)$. 
 Recall the constant $\alpha\in(1,2)$ in the hypothesis \ref{H1}.
 
 \subsubsection{Scaling limit for total budding rate processes}
 We first give some sufficient conditions on the initial state and parameters such that the rescaled CMJI-process converges to a non-degenerate limit.
 For any $i\in\mathcal{H}$ and $\boldsymbol{u}\in\mathbb{U}$, by (\ref{Kernel01}) we have 
 \beqnn
 \big\| \phi_i(\boldsymbol{u})\big\|_{L^1} =\sum_{l=1}^{k_{i}} \big\|\mathtt{B}_{i,l}( \mathtt{y}_{i,l})\big\|_{L^1},\quad
 \big\| \phi_i(\boldsymbol{u})\big\|_{\rm TV}\leq  \sum_{l=1}^{k_{i}} \big\|\mathtt{B}_{i,l}( \mathtt{y}_{i,l})\big\|_{\rm TV} 
 \eeqnn
 and
 \beqnn
 \int_0^\infty t\cdot \phi_i(t,\boldsymbol{u})dt =\sum_{l=1}^{k_{i}} \int_0^\infty t \cdot \mathtt{B}_{i,l}(t, \mathtt{y}_{i,l})dt.
 \eeqnn
 By H\"older's inequality, the hypothesis \ref{H1} is satisfied under the following condition.
 
 \begin{condition}\label{CMJMomentCondition}
 	For $i\in\mathcal{H}$ and $j\in\mathcal{D}$, there exists a constant $C>0$ such that for any $n\geq 1$,
 	\beqnn
 	\mathrm{v}_{\mathtt{L},i}^{(n)}+ \sum_{k_{\mathcal{H}}\in\mathbb{N}^d} \big|k_{\mathcal{H}}\big|^{2\alpha} p_j^{(n)}(k_{\mathcal{H}}) \leq C
 	\eeqnn
 	and 
 	\beqnn
 	\int_0^\infty \mathcal{P}^{(n)}_{\mathtt{L},i}(d\mathtt{y}) \int_\mathbb{B} \bigg(  \int_0^\infty t \cdot \mathtt{B}(t,\mathtt{y})dt + \big\|\mathtt{B}(\mathtt{y})\big\|_{\rm TV}   \bigg)^{2\alpha} \mathcal{P}_{\mathtt{B},i}^{(n)}(d\mathtt{B}) \leq C.
 	\eeqnn
 \end{condition}
 This condition means that both the branching and immigration mechanisms satisfy  the light-tailed condition. 
 More precisely, the number of new individuals in each immigration or born at each successive age  is light-tailed distributed. 
 Meanwhile, each individual is likely to give birth to its offspring in youth.  
 Since individuals give birth to their offspring independently, we have for each $i\in\mathcal{H}$ and $j\in\mathcal{D}$,
 \beqnn
 \phi_{ij}^{(n)}(t)=    \mathtt{B}_i^{(n)}(t)\cdot\mathrm{m}_{ij}^{(n)}, \quad  
 \big\|\phi_{ij}^{(n)}\big\|_{L^1}= \mathrm{m}_{\mathtt{B},i}^{(n)} \cdot  \mathrm{m}_{ij}^{(n)}, \quad 
 \int_0^\infty t \cdot \phi_{ij}^{(n)}(t)dt
 =\mathrm{d}_{\mathtt{B},i}^{(n)}\cdot \mathrm{m}_{ij}^{(n)} 
 \eeqnn
 and
 \beqnn  
 \int_\mathbb{U} \big\|\phi_i(\boldsymbol{u})\big\|_{L^1}^2 \nu_i^{(n)}(d\boldsymbol{u}) =  \mathrm{v}^{(n)}_{\mathtt{B},i} \cdot \mathrm{m}_{ii}^{(n)} + \big( \mathrm{v}_{ii}^{(n)}- \mathrm{m}_{ii}^{(n)}\big)\cdot \big|\mathrm{m}_{\mathtt{B},i}^{(n)}\big|^2. \qquad \qquad
 \eeqnn
 We now provide some asymptotic assumptions on the branching mechanism and immigration mechanism.

 \begin{condition}\label{ConvergenceParameterCMJ}
 Assume that hypothesis \ref{H2} holds and as $n\to\infty$,
 \begin{enumerate}
  \item[(1)] for each $i\in\mathcal{H}$, there exist constants $\mathrm{m}_{\mathtt{L},i}^*,\,\mathrm{v}_{ii}^*,\,\mathrm{d}_{\mathtt{B},i}^*,\,\mathrm{v}_{\mathtt{B},i}^*,\,\mathrm{m}_{ii}^*>0$ and $\mathrm{m}_{iI}^*\geq 0$ such that
  \beqnn
  \mathrm{m}_{\mathtt{L},i}^{(n)}\to \mathrm{m}_{\mathtt{L},i}^*,\quad \mathrm{v}_{ii}^{(n)}\to   \mathrm{v}_{ii}^*, \quad \mathrm{d}_{\mathtt{B},i}^{(n)} \to \mathrm{d}_{\mathtt{B},i}^*  
  \eeqnn
  and
  \beqnn
  \mathrm{v}^{(n)}_{\mathtt{B},i}\to  \mathrm{v}_{\mathtt{B},i}^*, \quad 
  \mathrm{m}_{ii}^{(n)} \to \mathrm{m}_{ii}^*, \quad
  \mathrm{m}_{iI}^{(n)} \to \mathrm{m}_{iI}^*;
  \eeqnn
 		
  \item[(2)] there exists a matrix $b_{\mathcal{H}^2}^*:=(b_{ij}^*)_{i,j\in \mathcal{H}}$ such that \beqnn
  n \Big(\mathrm{m}^{(n)}_{\mathcal{H}^2}-{\rm diag}(1/\mathrm{m}_{\mathtt{B},\mathcal{H}}^{(n)})\Big)\to b_{\mathcal{H}^2}^*.
  \eeqnn
 	\end{enumerate}
 \end{condition}
 
 The essence of this condition is that the rescaled branching and immigration mechanisms converge to a non-degenerate limit, i.e., immigrants enter into the population at the rate $\mathrm{m}_{\mathcal{H}I}^*$ and the net growth rate of the population is $b_{\mathcal{H}^2}^*$.
 We now provide some sufficient conditions on ancestors.
 For each $i\in\mathcal{H}$ and ancestor $x\in\mathcal{I}_{i,0}^{(n)}$,  denote by $\mathtt{A}^{(n)}_{i,x}$ and $\mathtt{R}^{(n)}_{i,x}$ its age and residual life at time $0$ respectively.
 Enlightened by the inspection paradox relating to the fact that observing a renewal interval at time $t$ gives an interval with average value larger than that of an average renewal interval; see Chapter~7.7 in \cite[p.460]{Ross2010}, we may assume that the residual life $\mathtt{R}^{(n)}_{i,x}$ is distributed as the \textit{excess life-length distribution}  of $\mathcal{P}^{(n)}_{\mathtt{L},i} $, also called \textit{forward recurrence time},  which is defined by
 \beqnn
 \breve{\mathcal{P}}^{(n)}_{\mathtt{L},i}(dy):=  \mathcal{P}^{(n)}_{\mathtt{L},i}[y,\infty)
 \cdot \frac{dy}{\mathrm{m}_{\mathtt{L},i}^{(n)}}.
 \eeqnn
 For an individual getting into the population at time $-t<0$, it stays alive at time $0$ with probability $\mathcal{P}^{(n)}_{\mathtt{L},i}[t,\infty)$. 
 Since the ancestor $x$ may get into the population at any time prior to time $0$, we may assume that its age $\mathtt{A}^{(n)}_{i,x}$ is distributed as
 $$
 \mathbf{P}( \mathtt{A}^{(n)}_{i,x}\in dt) = \mathcal{P}^{(n)}_{\mathtt{L},i}[t,\infty)
 \cdot \frac{dt}{\mathrm{m}_{\mathtt{L},i}^{(n)}} = \breve{\mathcal{P}}^{(n)}_{\mathtt{L},i}(dt).
 $$ 
 Taking these together, we assume that the next condition holds for ancestors.
 
 \begin{condition}\label{ConditionInitialStateCMJ}
 	For each $n\geq 1$ and $i\in\mathcal{H}$, assume that  $\mathtt{A}^{(n)}_{i,x}$ and $\mathtt{R}^{(n)}_{i,x}$ have  joint distribution
 	 \beqnn
 	 \breve{\mathcal{P}}_{\mathtt{AR},i}^{(n)}(dt,dy) := \mathbf{P}\big(\mathtt{A}^{(n)}_{i,x}\in dt, \mathtt{R}^{(n)}_{i,x}\in dy\big)= \frac{dt}{\mathrm{m}^{(n)}_{\mathtt{L},i}} \cdot \mathcal{P}^{(n)}_{\mathtt{L},i}(t+dy).
 	 \eeqnn
 \end{condition}
 Actually, Condition~\ref{ConditionInitialStateCMJ} is consistent with our previous assumptions on the age and residual-life distributions of ancestors.
 Indeed, it is easy to identify that the marginal distribution is
 \beqnn 
 \mathbf{P}(\mathtt{A}^{(n)}_{i,x}\in dt)=\breve{\mathcal{P}}_{\mathtt{AR},i}^{(n)}(dt,\mathbb{R}_+) = \mathcal{P}^{(n)}_{\mathtt{L},i}[t,\infty)\cdot \frac{dt}{\mathrm{m}^{(n)}_{\mathtt{L},i}}= \breve{\mathcal{P}}_{\mathtt{L},i}(dt).
 \eeqnn
 Moreover, for any $y\geq 0$ we also have
 \beqlb\label{marginaldistribution}
 \mathbf{P}( \mathtt{R}^{(n)}_{i,x}\geq y)=\breve{\mathcal{P}}_{\mathtt{AR},i}^{(n)}(\mathbb{R}_+,[y,\infty))
 \ar=\ar \int_0^\infty\mathcal{P}^{(n)}_{\mathtt{L},i}[t+y,\infty)\frac{dt}{\mathrm{m}^{(n)}_{\mathtt{L},i}}\cr \ar=\ar\int_y^\infty\mathcal{P}^{(n)}_{\mathtt{L},i}[t,\infty)\frac{dt}{\mathrm{m}^{(n)}_{\mathtt{L},i}}
 \eeqlb
 and hence $ \mathbf{P}( \mathtt{R}^{(n)}_{i,x}\in dy)= \breve{\mathcal{P}}^{(n)}_{\mathtt{L},i}(dy)$.
 We now give a scaling limit theorem for the total budding rate process.

 \begin{theorem}\label{ConvergenceBirthRateProcesses}
 	Under Condition~\ref{CMJMomentCondition}, \ref{ConvergenceParameterCMJ} and \ref{ConditionInitialStateCMJ}, if $\sup_{n\geq 1}\mathbf{E}[|\Xi^{(n)}_\mathcal{H}(0)/n|^{2\alpha}]<\infty$ and $ \Xi^{(n)}_\mathcal{H}(0)/n\overset{\rm d}\to \Xi^*_\mathcal{H}(0) \in \mathbb{R}_+^d$, then the rescaled process $\{\mathbf{B}^{(n)}_{\mathcal{H}}(nt)/n:t\geq 0\}$  converges weakly to  $\{\mathbf{B}_{\mathcal{H}}^*(t):t\geq 0\}$ in $\mathbf{D}([0,\infty),\mathbb{R}_+^d)$ as $n\to\infty$, where $\mathbf{B}_{\mathcal{H}}^*$ is the unique strong solution to (\ref{LimitDensity}) with
 	\beqnn
 	Z_i(0)= \frac{\Xi^*_i(0)}{\mathrm{m}_{\mathtt{L},i}^* \mathrm{m}_{ii}^*},\quad
 	a_i=\frac{\mathrm{m}^*_{iI}}{ \mathrm{m}^*_{ii}},\quad
 	b_{ij}= \frac{ b^*_{ij}}{ \mathrm{m}^*_{ii}}
 	\eeqnn
 	and 
 	\beqnn
 	c_i^2=\mathrm{v}^*_{\mathtt{B},i} \cdot \mathrm{m}^*_{ii}+ \frac{ \mathrm{v}^*_{ii}- \mathrm{m}^*_{ii}}{|\mathrm{m}_{ii}^*|^{2}}, \quad
 	\sigma_i= \mathrm{d}^*_{\mathtt{B},i}\cdot \mathrm{m}^*_{ii} ,\quad i,j\in\mathcal{H}.
 	\eeqnn
 \end{theorem}

 \subsubsection{Scaling limits for $\mathtt{T}$-CMJI-processes}
 In this section we study the convergence of rescaled multi-type CMJI-processes  counted with random characteristic by using the limit results in Section~\ref{Section.ShotNoise}. 
 For each $i\in\mathcal{H}$ and $j\in\mathcal{D}$, the mutual independence among individuals induces that in the $n$-th model, the mean total impact
 of type-$i$ offspring produced by a type-$j$ individual at each successive age is
 \beqnn
 \zeta_{ij}^{(n)}(t)= \int_{\mathbb{U}}\zeta_{i}(t,\boldsymbol{u})\nu^{(n)}_j(d\boldsymbol{u})=\mathtt{T}_i^{(n)}(t)\cdot \mathrm{m}_{ij}^{(n)},\quad t\geq 0.
 \eeqnn
 Condition~\ref{ConvergenceParameterCMJ} tells that the mean arrival rate of  type-$i$ immigrants is of the order of $1$.  
 From Theorem~\ref{ConvergenceBirthRateProcesses}, The mean rate of type-$j$ individuals giving birth to type-$i$ offspring is of the order of $n$ if $j=i$ and $1$ otherwise.
 Thus the main contribution to the $\mathtt{T}$-CMJI-process is made by individuals whose types are same to those of their parents.
 
 We first establish a convergence result for the rescaled process $\{ \mathbf{T}_{\mathtt{I},\mathcal{H}}^{(n)}(t) :t\geq 0\}$ with
 \beqnn
 \mathbf{T}_{\mathtt{I},\mathcal{H}}^{(n)}(t):= \frac{1}{n}\cdot \mathbf{T}_\mathcal{H}^{(n)}(nt),
 \eeqnn 
  in which the characteristic function of each individual represents the instantaneous rate at which it effects the host, e.g., toxin release rate and population size.
 From (\ref{Kernel02}), we have
 \beqlb
 \|\zeta_i^{(n)}(\boldsymbol{u})\|_{L^1} \ar=\ar \sum_{l=1}^{k_i} \|\mathtt{T}_{i,l}(\mathtt{y}_{i,l})\|_{L^1}, \label{eqn.3011} \\
 \|\zeta_i^{(n)}(\boldsymbol{u})\|_{\rm TV} \ar\leq\ar \sum_{l=1}^{k_i} \|\mathtt{T}_{i,l}(\mathtt{y}_{i,l})\|_{\rm TV} ,\label{eqn.301}
 \eeqlb
 for any $ \boldsymbol{u}\in\mathbb{U}$ and $ i\in\mathcal{H}$.
 It is easy to identify Condition~\ref{ConInstantaneousEffects02} and \ref{ConInstantaneousEffects} by the next condition.
 
 \begin{condition}\label{CMJInstantaneousConvergence}
 	Assume that 
 	\beqnn
 	\sup_{n\geq 1}\int_0^\infty \mathcal{P}^{(n)}_{\mathtt{L},i}(d\mathtt{y}) \int_\mathbb{T} \Big(  \|\mathtt{T}(\mathtt{y})\|_{L^1}^{\alpha}
 	+ \|\mathtt{T}(\mathtt{y})\|_{\rm TV}^{2\alpha}  \Big) \mathcal{P}_{\mathtt{T},i}^{(n)}(d\mathtt{T}) <\infty, \quad i\in\mathcal{H}.
 	\eeqnn
 	Moreover, there exist a constant $a^*_{\mathtt{I},\mathcal{H}}\in\mathbb{R}_+^d$ and a non-negative function $\bar{\mathtt{T}} \in L^1(\mathbb{R}_+) $ such that 
 	\beqnn
 	\lim_{n\to\infty}\big\|\mathtt{T}_\mathcal{H}^{(n)}\big\|_{L^1} =  a^*_{\mathtt{I},\mathcal{H}}
 	\quad\mbox{and}\quad
 	\sup_{n\geq 1}\big|\mathtt{T}_\mathcal{H}^{(n)}(t)\big|\leq \bar{\mathtt{T}}(t),\quad t\geq 0.
 	\eeqnn
 	 
 \end{condition}
 
 \begin{theorem}\label{ConvergenceCMJInstantaneousEffect}
 	Under Condition~\ref{CMJMomentCondition}, \ref{ConvergenceParameterCMJ}, \ref{ConditionInitialStateCMJ} and \ref{CMJInstantaneousConvergence}, we have
 	$\mathbf{T}_{\mathtt{I},\mathcal{H}}^{(n)}\overset{\rm d}\to \mathbf{T}_{\mathtt{I},\mathcal{H}}^*$ in $\mathbf{D}([0,\infty),\mathbb{R}_+^d)$ as $n\to\infty$ with $\mathbf{T}_{\mathtt{I},i}^*:= a^*_{\mathtt{I},i} \cdot\mathrm{m}_{ii}^*\cdot \mathbf{B}_i^*$ for each $i\in\mathcal{H}$.
 \end{theorem}
 
 As a corollary, we next give a scaling limit theorem for the population size process $ \Xi^{(n)}_\mathcal{H}$, which is a $\mathtt{T}$-CMJI-process with $\mathcal{P}^{(n)}_{\mathtt{T},\mathcal{H}}(\mathtt{T}(t,y)=\mathbf{1}_{\{y>t\}})=1$.
 In this case, we have for any $t\geq 0$ and $y>0$,
 \beqnn
 \|\mathtt{T}(y)\|_{\rm TV}=1, \quad \|\mathtt{T}(y)\|_{L^1}=y,\quad \mathtt{T}_\mathcal{H}(t)= \mathcal{P}_{\mathtt{L},\mathcal{H}}^{(n)}(t,\infty),\quad \|\mathtt{T}_\mathcal{H}^{(n)}\|_{L^1} = \mathrm{m}_{\mathtt{L},\mathcal{H}}^{(n)}.
 \eeqnn
 \begin{corollary}\label{Corollary.4.10} 
 	Under Condition~\ref{CMJMomentCondition}, \ref{ConvergenceParameterCMJ} and \ref{ConditionInitialStateCMJ}, we have  
 	\begin{enumerate}
 		\item[(1)]  The rescaled process $\{\Xi^{(n)}_{\mathcal{H}}(nt)/n:t\geq 0\}$ converges weakly to $\{\Xi^*_{\mathcal{H}}(t):= (\mathrm{m}_{\mathtt{L},i}^*\mathrm{m}^*_{ii} \cdot \mathbf{B}_i^*(t))_{i\in\mathcal{H}}:t\geq 0\}$ in $\mathbf{D}([0,\infty),\mathbb{R}_+^d)$  as $n\to\infty$;
 		
 		\item[(2)] If $\mathcal{P}^{(n)}_{\mathtt{L},\mathcal{H}} \overset{\rm d}\to\mathcal{P}^*_{\mathtt{L},\mathcal{H}}$ as $n\to\infty$, then for any constant $\eta >0$, the rescaled processes of total population which is younger than $\eta$, is older than $\eta$ or has residual life less than $\eta$ converge weakly  to 
 		\beqnn
 		\bigg( \frac{\int_0^\infty (y\wedge \eta) \mathcal{P}_{\mathtt{L},i}^*(dy) \cdot \Xi^*_i}{ \mathrm{m}_{\mathtt{L},i}^* }  \bigg)_{i\in\mathcal{H}} ,\quad 
 		\bigg( \frac{\int_\eta^\infty (y-\eta) \mathcal{P}_{\mathtt{L},i}^*(dy)\cdot \Xi^*_i}{ \mathrm{m}_{\mathtt{L},i}^* }  \bigg)_{i\in\mathcal{H}}
 		\eeqnn
 		or 
 		\beqnn
 		\bigg( \frac{\int_0^\infty (y\wedge \eta) \mathcal{P}_{\mathtt{L},i}^*(dy)\cdot \Xi^*_i}{ \mathrm{m}_{\mathtt{L},i}^* }  \bigg)_{i\in\mathcal{H}}  
 		\eeqnn
 		respectively in $\mathbf{D}([0,\infty), \mathbb{R}_+^d )$ as $n\to\infty$.  
 	\end{enumerate} 
 \end{corollary}
 
 We now consider the behavior at a large  time scale of the cumulative impact of microbes on the host, e.g. cumulative toxin release and total progeny. 
 By Corollary~\ref{Corollary.4.10}(1), the cumulative impact of individuals with same type as their parents on the host is of the order of $n^2$. However, the assumption  $ \Xi^{(n)}_\mathcal{H}(0)\sim \Xi^*_\mathcal{H}(0)\cdot n$ induces that the cumulative impact of ancestors is of the order of $n$ and can be asymptotically ignored. 
 Consequently, we have 
 \beqnn
 \mathbf{T}^{(n)}_i(nt)\sim \int_0^{nt}\int_\mathbb{U}\zeta_i(nt-s,\boldsymbol{u})N_i^{(n)}(ds,d\boldsymbol{u}),\quad t\geq 0,i\in\mathcal{H}
 \eeqnn
 as $n\to\infty$ and hence $\mathbf{E}[\mathbf{T}^{(n)}_i(nt)]$ is of the order of $n^2$.
 Thus it is natural to consider the weak convergence of the rescaled process $\{\mathbf{T}^{(n)}_{\mathtt{C},\mathcal{H}}(t) : t\geq 0 \}$ by using Theorem~\ref{ConvergenceCumulativeEffect}, where 
 \beqnn
 \mathbf{T}^{(n)}_{\mathtt{C},\mathcal{H}}(t):=\frac{1}{n^2}\cdot \mathbf{T}^{(n)}_\mathcal{H}(nt).
 \eeqnn
 From (\ref{Kernel02}), we have
 \beqnn
 \zeta_i(\infty,\boldsymbol{u})= \sum_{l=1}^{k_i} \mathtt{T}_{i,l}(y_{i,l},y_{i,l}) ,\quad \boldsymbol{u} \in\mathbb{U}.
 \eeqnn
 It is easy to see that Condition~\ref{MomentConCumulativeEffect} and \ref{ConditionCumulativeEffect} are satisfied under the following condition.
 \begin{condition}\label{CMJCumulative}
 	Assume that  $\mathtt{T}_{\mathcal{H}}^{(n)}(\infty)  \to a^*_{\mathtt{C},\mathcal{H}} \in\mathbb{R}_+^d$ as $n\to\infty$ and $\sup_{n\geq 1}|\mathtt{T}_{\mathcal{H}}^{(n)}(\infty)-\mathtt{T}_{\mathcal{H}}^{(n)}(t)|\to0$ as $t\to\infty$. Moreover, assume that
 	\beqnn
 	\sup_{n\geq1} \int_\mathbb{T}\mathcal{P}^{(n)}_{\mathtt{T},i}(d\mathtt{T})\int_0^\infty|\mathtt{T}(y,y)|^\alpha \mathcal{P}^{(n)}_{\mathtt{L},i}(dy)<\infty, \quad i\in\mathcal{H}.
 	\eeqnn 
 	
 \end{condition} 
 
 \begin{theorem}\label{ScalingLimitCumulativeCMJ}
 	Under Condition~\ref{CMJMomentCondition}, \ref{ConvergenceParameterCMJ}, \ref{ConditionInitialStateCMJ} and \ref{CMJCumulative},
 	we have $\mathbf{T}^{(n)}_{\mathtt{C},\mathcal{H}} \overset{\rm d}\to \mathbf{T}^*_{\mathtt{C},\mathcal{H}}$  in $\mathbf{D}([0,\infty),\mathbb{R}_+^d)$ as $n\to\infty$ with
 	\beqnn
 	\mathbf{T}^*_{\mathtt{C},i}(t):= a_{\mathtt{C},i}^* \cdot\mathrm{m}_{ii}^*\cdot \int_0^t \mathbf{B}_i^*(s)ds,\quad t\geq 0, i\in\mathcal{H}.
 	\eeqnn 
 \end{theorem}
 
 \begin{corollary}
 	Under Condition~\ref{CMJMomentCondition}, \ref{ConvergenceParameterCMJ}  and  \ref{ConditionInitialStateCMJ},
 	the two rescaled processes of total progeny and integral of population converge weakly in $\mathbf{D}([0,\infty),\mathbb{R}_+^d)$ to
 	\beqnn
 	\Big(  \frac{1}{\mathrm{m}_{\mathtt{L},i}^*} \cdot \int_0^t  \Xi^*_i(s)  ds\Big)_{i\in\mathcal{H}} 
 	\quad \mbox{and}\quad 
 	\Big(  \int_0^t  \Xi^*_i(s)  ds\Big)_{i\in\mathcal{H}}. 
 	\eeqnn
 	
 \end{corollary}

 \subsubsection{Scaling limits for population structures}
 In this section, we give some asymptotic results for the population structure of nearly critical multi-type CMJI-processes under the following condition.
 
 \begin{condition}\label{ConvergenceLifeDis}
 	Assume that $\mathcal{P}^{(n)}_{\mathtt{L},i} \overset{\rm d}\to \mathcal{P}^*_{\mathtt{L},\mathcal{H}}$ as $n\to\infty$ and 
 	\beqnn
 	\sup_{n\geq 1}\int_0^\infty y^{2\alpha}\mathcal{P}^{(n)}_{\mathtt{L},\mathcal{H}}(dy)<\infty. 
 	\eeqnn 
 \end{condition}
 In mathematical biology, the population structure is usually described by the {\it age-distribution} and {\it residual-life distribution} of all alive individuals in the population.
 In precise,  for  $i\in\mathcal{H}$, denote by $\mathcal{AR}^{(n)}_{i,t}(ds,dz) $ the joint distribution of age and residual life of all type-$i$ individuals alive at time $t$ in the $n$-th model, i.e.
 \beqnn
 \mathcal{AR}^{(n)}_{i,t}(ds,dz)\ar:=\ar \sum_{x\in\mathcal{I}^{(n)}_i} \mathbf{1}_{\{0\leq t- \tau_x< \ell_x \}} \cdot \delta_{(t- \tau_x,\ell_x-(t- \tau_x))}(ds,dz)
 \eeqnn
 is a measure on $\mathbb{R}_+^2$ with unit mass at the age and residual-life of each type-$i$ individual alive at time $t$.
 The two marginal measures
 \beqnn
 \mathcal{A}^{(n)}_{i,t}(ds):=\mathcal{AR}^{(n)}_{i,t}(ds,\mathbb{R}_+)
 \quad\mbox{and}\quad
 \mathcal{R}^{(n)}_{i,t}(dz):=\mathcal{AR}^{(n)}_{i,t}(\mathbb{R}_+,dz)
 \eeqnn
 are the corresponding age distribution and residual-life distribution at time $t$.
 Similarly, the life-length distribution $\mathcal{L}^{(n)}_{i,t}$ of all type-$i$ individuals alive at time $t$ is given by
 \beqnn
 \mathcal{L}^{(n)}_{i,t}(dy)\ar:=\ar \sum_{x\in\mathcal{I}^{(n)}_i} \mathbf{1}_{\{0\leq t- \tau_x< \ell_x \}} \delta_{\ell_x}(dy).
 \eeqnn
 We establish a scaling limit  for the population structure in the next theorem in collaboration with the following three probability laws 
 \beqnn
 \breve{\mathcal{P}}_{\mathtt{AR},i}^*(ds,dz)
 \ar:=\ar \frac{ds}{|\mathrm{m}^*_{\mathtt{L},i}|}  \cdot \mathcal{P}^*_{\mathtt{L},i}(s+dz),\cr
 \breve{\mathcal{P}}_{\mathtt{L},i}^*(dy) 
 \ar:=\ar \mathcal{P}^*_{\mathtt{L},i}[y,\infty) \cdot \frac{dy}{|\mathrm{m}^*_{\mathtt{L},i}|} , \cr
 \mathring{\mathcal{P}}_{\mathtt{L},i}^*(dy) 
 \ar:=\ar \frac{y}{|\mathrm{m}^*_{\mathtt{L},i}|}\cdot \mathcal{P}^*_{\mathtt{L},i}(dy).
 \eeqnn
 Wherein, $\mathring{\mathcal{P}}_{\mathtt{L},i}^*$ is usually known as the \textit{size-biased distribution} of $ \mathcal{P}_{\mathtt{L},i}^*$.

 \begin{theorem}\label{ConvergencePopulationStructure}
 	Under Condition~\ref{CMJMomentCondition},~\ref{ConvergenceParameterCMJ}, \ref{ConditionInitialStateCMJ} and \ref{ConvergenceLifeDis}, we have  as $n\to\infty$,
 	\begin{enumerate}
 		\item[(1)] $  \{\mathcal{AR}^{(n)}_{\mathcal{H},nt}/n:t\geq 0\} \overset{\rm d}\to \{( \Xi^*_i(t) \cdot \breve{\mathcal{P}}_{\mathtt{AR},i}^*)_{i\in\mathcal{H}} :t\geq 0\} $ in $\mathbf{D}([0,\infty),  \mathcal{M}(\mathbb{R}_+^2)^d)$;
 		
 		\item[(2)] both the two rescaled processes $  \{\mathcal{A}^{(n)}_{\mathcal{H},nt}/n :t\geq 0\}  $ and $  \{\mathcal{R}^{(n)}_{\mathcal{H},nt}/n  :t\geq 0\}  $ converge weakly to   $ \{( \Xi^*_i(t)\cdot \breve{\mathcal{P}}_{\mathtt{L},i}^* )_{i\in\mathcal{H}}:t\geq 0\} $ in $\mathbf{D}([0,\infty),  \mathcal{M}(\mathbb{R}_+)^d)$;
 		
 		\item[(3)] $\{\mathcal{L}^{(n)}_{\mathcal{H},nt}/n : t\geq 0\}\overset{\rm d}\to\{( \Xi^*_i(t)\cdot \mathring{\mathcal{P}}_{\mathtt{L},i}^* )_{i\in\mathcal{H}}:t\geq 0\}$ in $\mathbf{D}([0,\infty),  \mathcal{M}(\mathbb{R}_+)^d)$.
 	\end{enumerate} 
 \end{theorem}
 
 \begin{remark}
 	The essence of Theorem~\ref{ConvergencePopulationStructure} is that as the rescaled measure-valued process $\mathcal{AR}^{(n)}_{\mathcal{H},n\cdot}/n$ approaches to the limit, it can be recovered from the diffusion scaled population process $\Xi^{(n)}_{\mathcal{H}}(n\cdot)/n$ by the lifting map $ \breve{\mathcal{P}}_{\mathtt{AR},\mathcal{H}}^* $ from $\mathbb{R}_+^d$ to $\mathcal{M}(\mathbb{R}_+^2)^d$.
 	In other words, in a complex biological system enjoying short-memory property, the evolution of population can be fully described by the process of population size together with the life-length distribution, with more detailed information about the population not being necessary. 
 	This type of asymptotic behavior is well known as {\rm state space collapse}.
 	It was first systematically investigated in \cite{Bramson1998,Williams1998} in the study of multi-class queueing systems and since then has been widely observed in heavy traffic limits of various queuing systems; see \cite{Gromoll2004,Stolyar2004,VerloopAyestaNunez-Queija2011}.
 	
 \end{remark}
 
 \begin{remark}
  Compared to the complex structure of CMJI-processes,  the approximating models in the preceding theorems can be useful for several reasons. Firstly, they have simpler structures and are easier to be understood than CMJI-models counted with random characteristic. Each coefficient in the limit models has an intuitive and understandable interpretation. 
  Secondly, their properties are usually consistent with those of CMJI-processes, e.g. criticality, extinction and stationarity; see \cite{Jagers1975,KyprianouPalau2018,Xu2014}. 
  Thirdly, compared to the non-parametric estimation of CMJI-models, the approximating models are computationally more tractable and only few parameters are needed to be estimated.  
 \end{remark}

  \section{Proofs}\label{Proofs}
  
  In this section, we give the detailed proofs for the main results in the previous sections including Theorem~\ref{MainThm01}, \ref{ConvergenceCumulativeEffect}, \ref{ConvergenceInstantanuousEffect},  \ref{ConvergenceInstantanuousEffectAncestor}, \ref{ConvergenceBirthRateProcesses}, \ref{ConvergenceCMJInstantaneousEffect}, \ref{ScalingLimitCumulativeCMJ} and \ref{ConvergencePopulationStructure}.

  \subsection{Proof for Theorem~\ref{MainThm01}} \label{Section6.1}
  By the argument at the beginning of Section~\ref{Sec.SL}, it suffices to prove the weak convergence of the sequence $\{Z^{(n)}_{\beta,\mathcal{H}}\}_{n\geq 1}$ to $Z_{\beta,\mathcal{H}}$.
  In order to simplify the following statements and notation, we prove this result with $\lambda_b<0$ (equivalently, $b_{ii}<0$ for all $i\in\mathcal{H}$) and $ \beta=0$.
  The general case can be proved similarly.
  The asymptotic analysis in Section~\ref{Sec.AsymptoticAna} has shown that the time-scaled functions $R_{ii}^{(n)}(n\cdot)$, $R_{ij}^{(n)}(n\cdot)$, $R_{iI}^{(n)}(n\cdot)$ and $R_{i}^{(n)}(n\cdot,u)$ can be approximated respectively  by the corresponding exponential functions. 
  The errors are denoted as: for $i,j\in\mathcal{H}$ with $i\neq j$ and $(t,u)\in\mathbb{R}_+\times \mathbb{U}$,
  \beqnn
  \varepsilon_{R_{ii}}^{(n)}(t)\ar:=\ar R_{ii}^{(n)}(nt)- \frac{1}{\sigma_i}e^{\frac{b_{ii}}{\sigma_i}  t},  \quad \ \ \ \
  \varepsilon_{R_{ij}}^{(n)}(t):= n        R_{ij}^{(n)}(nt)- \frac{b_{ij}}{\sigma_i}e^{\frac{b_{ii}}{\sigma_i}t},\cr
  \varepsilon_{R_{iI}}^{(n)}(t)\ar:=\ar R_{iI}^{(n)}(nt)- \frac{a_{i}}{\sigma_i}e^{\frac{b_{ii}}{\sigma_i} t}, \quad \quad
  \varepsilon_{R_{i}}^{(n)}(t,u):= R_{i}^{(n)}(nt,u) - \frac{ \|\phi_i(u)\|_{L^1} }{\sigma_i}e^{\frac{b_{ii}}{\sigma_i}t}.
  \eeqnn
  The sum of the first two terms on the right side of (\ref{ScaledDensity}) can be approximated by $Z_i^{(n)}(0) e^{b_{ii}/\sigma_i\cdot t}$ and other terms can be approximated respectively by the corresponding (stochastic) integrals with the integrand replaced by the limit exponential function.
  Meanwhile, the error processes have the following representations respectively: for $i,j\in\mathcal{H}$ and $t\geq 0$,
  \beqnn
  \varepsilon^{(n)}_{\mu_{i}} (t)\ar:=\ar \frac{ \mu^{(n)}_{i}(nt)}{n} +    R^{(n)}_{ii}*\frac{\mu_{i}^{(n)} }{n}(nt) - Z_i^{(n)}(0) e^{\frac{b_{ii}}{\sigma_i} t  },\cr
  \tilde\varepsilon^{(n)}_{iI}(t)\ar:=\ar \int_0^t \varepsilon_{R_{iI}}^{(n)}(s)ds ,  \quad\quad  \  \varepsilon^{(n)}_{iI}(t) :=  \int_0^t \int_{\mathbb{U}} \frac{\varepsilon_{R_{i}}^{(n)}(t-s,u)}{n} \tilde{N}_I^{(n)}(n\cdot ds,du),\cr
  \tilde\varepsilon^{(n)}_{ij}(t)\ar:=\ar \varepsilon_{R_{ij}}^{(n)}*Z^{(n)}_j(t), \quad       \quad  \varepsilon^{(n)}_{ij}(t):= \int_0^t \int_{\mathbb{U}} \frac{\varepsilon_{R_{i}}^{(n)}(t-s,u)}{n} \tilde{N}_{j}^{(n)}(n\cdot ds,du).
  \eeqnn 
  Let $E^{(n)}_i :=\varepsilon^{(n)}_{\mu_{i}}  + \sum_{j\in\mathcal{D}_i}\tilde\varepsilon^{(n)}_{ij} + \sum_{j\in\mathcal{D}}\varepsilon^{(n)}_{ij} $.
  Based on these notation, we can write (\ref{ScaledDensity}) under the form
  \beqnn
  Z_{i}^{(n)}(t)
  \ar=\ar Z_i^{(n)}(0)  e^{ \frac{b_{ii}}{\sigma_i}t  }  + E^{(n)}_i(t)+ \int_0^t \frac{a_i}{\sigma_i} e^{ \frac{b_{ii}}{\sigma_i}(t-s)}  ds  + \sum_{j\in\mathcal{H}_i}\int_0^t \frac{b_{ij}}{\sigma_i} e^{\frac{b_{ii}}{\sigma_i} (t-s)} Z_{j}^{(n)}(s)  ds  \cr
  \ar\ar +   \sum_{j\in\mathcal{D}}\int_0^t \int_{\mathbb{U}}  e^{\frac{b_{ii}}{\sigma_i} (t-s)} \cdot \frac{\|\phi_i(u)\|_{L^1}}{n}\frac{1}{\sigma_i}  \tilde{N}_{j}^{(n)}(n\cdot ds,du),\quad t\geq 0, i\in\mathcal{H}.
  \eeqnn
  Using the fact that $e^{ \frac{b_{ii}}{\sigma_i}(t-s)}= 1+ \frac{b_{ii}}{\sigma_i} \int_s^t e^{\frac{b_{ii}}{\sigma_i} (r-s)}dr$ and Fubini's theorem, we also can write the foregoing equation into the following convenient form:
  \beqlb \label{SemiMarSVR}\qquad
  Z_{i}^{(n)}(t)
  \ar=\ar Z_i^{(n)}(0)+ E^{(n)}_i(t)    + \sum_{j\in\mathcal{D}}M_{ij}^{(n)} (t)\cr
  \ar\ar + \int_0^t  \Big(\frac{b_{ii}}{\sigma_i}E^{(n)}_i(s)+\frac{a_i}{\sigma_i}  + \sum_{j\in\mathcal{H}} \frac{b_{ij}}{\sigma_i} Z_{j}^{(n)}(s) \Big) ds  ,
  \quad i\in\mathcal{H},
  \eeqlb
  where $M_{ij}^{(n)} $ is an $(\mathscr{F}_{nt})$-local martingale defined as (\ref{MartMiI}) with $\beta=0$. By (\ref{Ni2N0i}) and Proposition~\ref{SVR}, for each $j\in\mathcal{H}$ we also can write $M_{ij}^{(n)}$  under the form
  \beqnn
  M_{ij}^{(n)}(t)= \int_0^t \int_\mathbb{U} \int_0^{Z^{(n)}_j(s-)} \frac{\|\phi_i(u)\|_{L^1}}{n}\frac{1}{\sigma_i}  \tilde{N}_{0,j}^{(n)}(n\cdot ds,du,n\cdot dz),\quad t\geq 0,
  \eeqnn
  where $\tilde{N}_{0,j}^{(n)}(n\cdot ds,du,n\cdot dz),\ j\in\mathcal{H}$ are $d$ mutually orthogonal compensated Poisson random measures on $(0,\infty)\times \mathbb{U}\times \mathbb{R}_+$ with intensity $ n^2\cdot ds\nu_j^{(n)}(du)dz$ respectively and also independent of  $N_I^{(n)}(ds,du)$.

  We now start to prove Theorem~\ref{MainThm01} by using the convergence results for infinite-dimensional stochastic differential equations established by Kurtz and Protter \cite[Theorem~7.5]{KurtzProtter1996}.
  The existence and uniqueness of solutions to (\ref{LimitDensity}) follow from Theorem 1 in \cite{YamadaWatanabe1971}.
  We now write (\ref{SemiMarSVR}) into the form of a stochastic integral and differential equation driven by an infinite-dimensional semimartingale; see Appendix~\ref{AppendixB}.
  Let $\mathbb{H}:= \mathbb{R} \times (L^{2}(\mathbb{R}_+))^d $ be a separable Banach space endowed with a norm $\|\cdot\|_{\mathbb{H}}$ defined by $ \|x\|_\mathbb{H}= |x_0|+ \sum_{i=1}^d\|x_i\|_{L^{2}}$ for  $x:=(x_0,x_1,\cdots, x_d) \in \mathbb{H}$.
  For each $n\geq 1$, we define a process $ U_\mathcal{H}^{(n)}$ by
  \beqnn
  U_i^{(n)}(t)=Z_i^{(n)}(0)+ E^{(n)}_i(t) +\int_0^t  \frac{b_{ii}}{\sigma_i}E^{(n)}_i(s)ds  +\sum_{j\in\mathcal{D}_i}M_{ij}^{(n)} (t)+ \frac{a_i}{\sigma_{i}}\cdot t,\quad t\geq 0, i\in\mathcal{H}
  \eeqnn 
  and a standard $\mathbb{H}^\#$-semimartingle $\boldsymbol{Y}^{(n)}:=(\boldsymbol{Y}_0^{(n)},W^{(n)}_1,\cdots,W^{(n)}_d)$ by $\boldsymbol{Y}_0^{(n)}(t):= t$ and 
  \beqnn 
  W_i^{(n)}(t):= \int_0^t \int_{\mathbb{U}}   \frac{\|\phi_i(u)\|_{L^1}}{n}\frac{1}{c_i}  \tilde{N}_{0,j}^{(n)}(n\cdot ds,du,n\cdot dz),\quad t\geq 0, i\in\mathcal{H}.
  \eeqnn
  We can rewrite the rescaled intensity process (\ref{SemiMarSVR}) as
  \beqnn
  Z^{(n)}_\mathcal{H}(t)= U_\mathcal{H}^{(n)}(t) + \boldsymbol{F}_\mathcal{H}(Z^{(n)}_\mathcal{H}(-))\cdot \boldsymbol{Y}^{(n)}(t),
  \eeqnn
  where $ \boldsymbol{F}_\mathcal{H}:=( \boldsymbol{F}_i)_{i\in\mathcal{H}} : \mathbb{R}_+^d \mapsto \mathbb{H}^d $ with the function $\boldsymbol{F}_i$ defined by
  \beqnn
  \boldsymbol{F}_i(x_\mathcal{H}):= \Big( \sum_{j\in\mathcal{H}} \frac{b_{ij}}{\sigma_i} x_{j}, 0,\cdots,0 ,\underbrace{\frac{c_i}{\sigma_i}\cdot\mathbf{1}_{\{z<x_i\}}}_{\tiny i\mbox{-th element}},0,\cdots,0\Big) \in \mathbb{H},\quad  x_\mathcal{H}\in\mathbb{R}_+^d.
  \eeqnn
  By \cite[Example~5.3]{KurtzProtter1991}, the function $ \boldsymbol{F}_\mathcal{H}$ satisfies conditions in \cite[Theorem~7.5]{KurtzProtter1996}. 
  The desired weak convergence in Theorem~\ref{MainThm01} follows immediately if the sequence of $\mathbb{H}^{\#}$-semimartingales $\{ \boldsymbol{Y}^{(n)} \}_{n\geq 1}$ is uniformly tight; see Definition~\ref{Definition.A1}, and $( U_\mathcal{H}^{(n)},\boldsymbol{Y}^{(n)})\Rightarrow(U_\mathcal{H},\boldsymbol{Y})$ as $n\to\infty$, where 
  \beqnn
  U_\mathcal{H}(t):= \Big(Z_i(0) + \frac{a_i}{\sigma_i}\cdot t\Big)_{i\in\mathcal{H}}
  \quad \mbox{and}\quad
  \boldsymbol{Y}^{(n)}(t):=\big(t,W_1(t,\cdot), \cdots,W_d(t,\cdot)\big).
  \eeqnn
  Actually, they follow directly from the next three claims:
  
  \begin{enumerate}
  	\item[$\bullet$] The two processes $E^{(n)}_\mathcal{H}$ and $\int_0^\cdot E^{(n)}_\mathcal{H}(s)ds $ converge weakly to $0$ in $\mathbf{D}([0,\infty),\mathbb{R}^d)$; see Section~\ref{Sec.ErrorP}. \vspace{7pt}
  	
  	\item[$\bullet$] For each $i\in\mathcal{H}$ and $j\in\mathcal{D}_i$, the local martingale $M_{ij}^{(n)}$ converges weakly to $0$ in $\mathbf{D}([0,\infty),\mathbb{R})$; see Section~\ref{ConvergenceMart}.   \vspace{7pt}
  	
  	\item[$\bullet$] The sequence of $((L^{2}(\mathbb{R}_+))^d)^{\#}$-local martingales $\big\{W_\mathcal{H}^{(n)}\big\}_{n\geq 1}$ is uniformly tight and  $W_\mathcal{H}^{(n)} \Rightarrow W_\mathcal{H}$ ; see Section~\ref{ConvergenceW}.
  \end{enumerate}

  \subsubsection{Negligible error functions $ \{\varepsilon_{R_{i}}^{(n)},\varepsilon_{R_{ij}}^{(n)}\}_{i\in\mathcal{H},j\in\mathcal{D}}$\,}
  
  In this section,  we prove the convergence of error functions $ \big\{\varepsilon_{R_{ij}}^{(n)}\big\}_{i\in\mathcal{H},j\in\mathcal{D}}$ and $\big\{ \varepsilon_{R_{i}}^{(n)}\big\}_{i\in\mathcal{H}}$ to $0$ in $L^1$ or $L^2$.  
  From the hypothesis~\ref{H2}, we have
  \beqlb\label{UpperBound.Rii}
  \sup_{n\geq 1} \big\| R^{(n)}_{ii} \big\|_{L^\infty} <\infty.
  \eeqlb
  Moreover, by extending the proofs of Lemma~4.2 and 4.4 in \cite{JaissonRosenbaum2015}, we can get the following helpful estimates for the Fourier transform $\hat\phi_{ii}^{(n)}$ with $i\in\mathcal{H}$; readers may refer to Appendix~\ref{AppendixA} for the detailed proof.

 \begin{proposition}\label{Thm.UpperBound} \it
  There exist constants $C_1,C_2,n_0>0$  such that for any $i\in\mathcal{H}$, $n\geq n_0$ and $\lambda\in\mathbb{R}$,
  \beqlb\label{LowerBound}
  \big|\hat\phi^{(n)}_{ii}(\lambda) \big| \leq C_1\big(|\lambda|^{-1}\wedge  1\big)
  \quad\mbox{and}\quad
  \big|1-\hat{\phi}_{ii}^{(n)}(\lambda)\big| \geq C_2(|\lambda|\wedge 1).
  \eeqlb
 \end{proposition}

 \begin{lemma}\label{L2ConvergenceRii}
 For each $i\in\mathcal{H}$, we have $	\big\|\varepsilon_{R_{ii}}^{(n)}\big\|_{L^2}\to0$ as $n\to\infty$.
 \end{lemma}
  \proof We first provide an upper bound for the Fourier transform of $R_{ii}^{(n)}(n\cdot)$.  By Condition~\ref{Con.PhiH},
  \beqnn
  \Big|\int_0^\infty e^{\mathrm{i}\lambda t} R_{ii}^{(n)}(nt)dt\Big| \leq \int_0^\infty R_{ii}^{(n)}(nt)dt 
  =\frac{\|\phi^{(n)}_{ii}\|_{L^1}}{n(1-\|\phi^{(n)}_{ii}\|_{L^1})},
  \eeqnn
  which converges to $-1/ b_{ii}>0$ as $n\to\infty$ and hence there  exist  constants $C,n_0>0$ such that
  \beqnn
  \sup_{n\geq n_0}\Big|\int_0^\infty e^{\mathrm{i}\lambda t} R_{ii}^{(n)}(nt)dt\Big| \leq C.
  \eeqnn
  On the other hand, by Proposition~\ref{Thm.UpperBound} we also have for large $n$,
  \beqnn
  \Big|\int_0^\infty e^{\mathrm{i}\lambda t} R_{ii}^{(n)}(nt)dt\Big|
  \leq \frac{|\hat\phi^{(n)}_{ii}(\lambda/n)|}{n|1-\hat\phi^{(n)}_{ii}(\lambda/n)|}
  \leq C \frac{\frac{n}{|\lambda|}\wedge 1}{|\lambda|\wedge n}= \frac{C}{|\lambda|}.
  \eeqnn
  Putting these two estimates together, there exist two constants $C,n_0>0$ such that for any $\lambda\in\mathbb{R}$,
  \beqlb\label{UpperBoundRii}
  \sup_{n\geq n_0} \Big| \int_0^\infty e^{\mathrm{i}\lambda t} R^{(n)}_{ii}(nt)dt\Big|\ar \leq\ar C(|\lambda|^{-1}\wedge 1).
  \eeqlb
  Thus the Fourier transforms of $R_{ii}^{(n)}(n\cdot)$ and $\varepsilon_{R_{ii}}^{(n)}$  are square integrable.
  By the Fourier isometry,
  \beqnn
  \|\varepsilon_{R_{ii}}^{(n)}\|_{L^2}^2
  =\int_0^\infty \big|\varepsilon_{R_{ii}}^{(n)}(t)\big|^2dt
  = \int_{\mathbb{R}}  \Big|\int_0^\infty e^{\mathrm{i}\lambda t} \varepsilon_{R_{ii}}^{(n)}(t) dt \Big|^2d\lambda.
  \eeqnn
  By the dominated convergence theorem, (\ref{UpperBoundRii}) and (\ref{ConvergenceFourierRii}), we can get the desired result immediately. 
  \qed
  
  \begin{lemma}\label{L1ConvergenceRij}
  	For any $T>0$, $i\in\mathcal{H}$ and $ j\in\mathcal{D}_i$, we have $ \|\varepsilon_{R_{ij}}^{(n)}\|_{L^1_T}\to0$ as $n\to\infty$.
  \end{lemma}
  \proof Here we just prove this lemma with $i,j\in\mathcal{H}$ and $i\neq j$.
  For the case of $j=I$, it can be proved in the same way.
  By H\"older's inequality,
  \beqnn
  \|\varepsilon_{R_{ij}}^{(n)}\|_{L^1_T}
  \ar\leq \ar \int_0^T n\phi^{(n)}_{ij}(nt)dt + \int_0^T  \Big| n\int_0^{nt}R^{(n)}_{ii}(nt-s)\phi^{(n)}_{ij}(s)ds - \frac{b_{ij}}{\sigma_i} e^{\frac{b_{ii}}{\sigma_i}t} \Big|dt\cr
  \ar\leq\ar \|\phi^{(n)}_{ij}\|_{L^1} + \sqrt{T}\Big(  \int_0^\infty  \Big| n\int_0^{nt}R^{(n)}_{ii}(nt-s)\phi^{(n)}_{ij}(s)ds - \frac{b_{ij}}{\sigma_i} e^{\frac{b_{ii}}{\sigma_i}t} \Big|^2dt\Big)^{1/2}.
  \eeqnn
  Here the first term on the right side of the last inequality vanishes as $n\to\infty$; see Condition~\ref{Con.PhiH}.
  By the convolution theorem and Condition~\ref{Con.PhiH}, we have as $n\to\infty$,
  \beqnn
  \int_0^\infty e^{\mathrm{i}\lambda t}dt \cdot n\int_0^{nt}R^{(n)}_{ii}(nt-s)\phi^{(n)}_{ij}(s)ds 
  \ar=\ar \frac{\hat{\phi}^{(n)}_{ii}(\frac{\lambda}{n})\cdot n\hat{\phi}^{(n)}_{ij}(\frac{\lambda}{n})}{n(1- \hat{\phi}^{(n)}_{ii}(\frac{\lambda}{n}))}\cr
  \ar\to\ar
  \frac{-b_{ij}}{ b_{ii}+\mathrm{i}\sigma_i \lambda} \cr
  \ar=\ar \int_0^\infty e^{\mathrm{i}\lambda t} \frac{b_{ij}}{\sigma_i} e^{\frac{b_{ii}}{\sigma_i}t}dt.
  \eeqnn
  Moreover, notice that $\sup_{n\geq 1}|n\hat{\phi}^{(n)}_{ij}(\lambda/n)| \leq \sup_{n\geq 1} n\|\phi^{(n)}_{ij}\|_{L^1}<\infty$, by Proposition~\ref{Thm.UpperBound} we have
  \beqnn
  \sup_{n\geq 1} \Big| \int_0^\infty e^{\mathrm{i}\lambda t}dt \cdot n\int_0^{nt}R^{(n)}_{ii}(nt-s)\phi^{(n)}_{ij}(s)ds  \Big| \leq C (|\lambda|^{-1}\wedge 1).
  \eeqnn
  By the Fourier isometry,
  \beqnn
  \lefteqn{\int_0^\infty  \Big| n\int_0^{nt}R^{(n)}_{ii}(nt-s)\phi^{(n)}_{ij}(s)ds - \frac{b_{ij}}{\sigma_i} e^{-\frac{b_{ii}}{\sigma_i}t} \Big|^2dt}\ar\ar\cr
  \ar=\ar \int_0^\infty  \Big| \frac{\hat{\phi}^{(n)}_{ii}(\lambda/n) \cdot n\hat{\phi}^{(n)}_{ij}(\lambda/n) }{n(1- \hat{\phi}^{(n)}_{ii}(\lambda/n))} + \frac{b_{ij}}{b_{ii}+\mathrm{i} \sigma_i \lambda}   \Big|^2d\lambda.
  \eeqnn
  By the dominated convergence theorem, it vanishes as $n\to\infty$ and the proof is completed.
  \qed
  
  \begin{lemma}\label{L2ConvergenceRi}
  	There exists a sequence $\{\epsilon_n\}_{n\geq 1}$ vanishing as $n\to\infty$ such that  for any $n\geq 1$, $u\in\mathbb{U}$ and $i\in\mathcal{H}$
  	\beqnn
  	\big\|\varepsilon_{R_{i}}^{(n)}(u)\big\|_{L^2} \leq \epsilon_n \cdot \Phi(u).
  	\eeqnn
  \end{lemma}
  \proof Notice that $|\varepsilon_{R_{i}}^{(n)}(t,u)|\leq \phi_i(nt,u)+ |A^{(n)}_1(t,u)|+ |A^{(n)}_2(t,u)|$, where
  \beqnn
  A^{(n)}_1(t,u)\ar:=\ar \int_0^t  n \phi_i(n(t-s),u) \varepsilon_{R_{ii}}^{(n)}(s)  ds, \\
  A^{(n)}_2(t,u)\ar:=\ar   \int_0^t  n \phi_i(n(t-s),u)  \frac{1}{\sigma_i}e^{\frac{b_{ii}}{\sigma_i} s}  ds - \frac{\|\phi_i(u)\|_{L^1}}{\sigma_i}  e^{\frac{b_{ii}}{\sigma_i}t}.
  \eeqnn
  By the hypothesis \ref{H1} we first have
  \beqnn
  \int_0^\infty |\phi_i(nt,u)|^2dt  \leq  \frac{ \|\phi_i(u)\|_{\rm TV} \cdot \|\phi_i(u)\|_{L^1} }{n}
  \leq \frac{  |\Phi(u)|^2}{n} .
  \eeqnn
  Moreover, by Young's convolution inequality and the hypothesis \ref{H1},  
  \beqnn
  \|A^{(n)}_1(u)\|_{L^2} \leq  \|\phi_i(u)\|_{L^1} \cdot \|\varepsilon_{R_{ii}}^{(n)}\|_{L^2} \leq \|\varepsilon_{R_{ii}}^{(n)}\|_{L^2} \cdot \Phi(u).
  \eeqnn 
  Taking Fourier transform of $A^{(n)}_2(t,u)$, we have
  \beqnn
  \int_0^\infty e^{\mathrm{i}\lambda t} A^{(n)}_2(t,u) dt \ar=\ar \frac{\|\phi_i(u)\|_{L^1}-\hat{\phi}_i(\lambda/n,u)}{b_{ii}+\mathrm{i}\sigma_i\lambda }
  \eeqnn
  and by the Fourier isometry,
  \beqnn
  \|A^{(n)}_2(u)\|_{L^2}^2 \ar=\ar \int_{\mathbb{R}} \Big| \frac{\|\phi_i(u)\|_{L^1}-\hat{\phi}_i(\lambda/n,u)}{b_{ii}+\mathrm{i}\sigma_i\lambda} \Big|^2d\lambda .
  \eeqnn
  From the facts that $|e^{\mathrm{i}\frac{\lambda}{n} t} - 1 | \leq (|\lambda|t/n)\wedge 2$ and $|b_{ii}+\mathrm{i}\sigma_i\lambda|^{-1}\leq C/(1+|\lambda|)$, we  also have
  \beqnn
  \|A^{(n)}_2(u)\|_{L^2}^2
  \ar\leq\ar C\int_{\mathbb{R}}  \Big| \int_0^\infty \Big(\frac{|\lambda|t}{n} \wedge 1\Big) \phi_i(t,u)dt  \Big|^2   \frac{d\lambda}{(1+|\lambda|)^{2}}\cr
  \ar\leq\ar C \int_{\mathbb{R}} \Big(\Big|\frac{\lambda}{n} \int_0^\infty t\phi_i(t,u)dt\Big| \wedge  \|\phi_i(u)\|_{L^1}\Big)^2  \frac{d\lambda}{(1+|\lambda|)^{2}} \cr
  \ar\leq\ar C \int_{\mathbb{R}} \Big(\frac{|\lambda|}{n}\wedge 1\Big)^2\frac{ d\lambda}{(1+|\lambda|)^2} \cdot  |\Phi
  (u)|^2
  \leq \frac{C}{n}\cdot  |\Phi(u)|^2.
  \eeqnn
  Putting these estimates together, we can immediately get the desired result with $\epsilon_n:= C\cdot ( \|\varepsilon_{R_{ii}}^{(n)}\|_{L^2} \vee n^{-1/2})$ for some large constant $C>0$.
  \qed
  
  \subsubsection{Weak convergence of error processes} \label{Sec.ErrorP} 
  If $E^{(n)}_{\mathcal{H}}\overset{\rm d}\to 0$ in $\mathbf{D}([0,\infty),\mathbb{R}^d)$, by Proposition~1.17(b) in \cite[p.328]{JS03} we have 
  \beqnn
  E^{(n)}_{\mathcal{H}}\overset{\rm u.c.p.}\longrightarrow 0
  \quad\mbox{and hence}
  \int_0^\cdot E^{(n)}_{\mathcal{H}}(s)ds \overset{\rm u.c.p.}\longrightarrow 0 . 
  \eeqnn 
  For the weak convergence of $E^{(n)}_{\mathcal{H}}$ to $0$,  by Corollary~3.33 in \cite[p.353]{JS03} it suffices to prove separately that  $ \varepsilon^{(n)}_{\mu_{i}} $, $\tilde\varepsilon^{(n)}_{ij}$ and $\varepsilon^{(n)}_{ij}$ converge weakly to $0$ in $\mathbf{D}([0,\infty),\mathbb{R})$ as $n\to\infty$ for each $i\in \mathcal{H}$ and $j\in \mathcal{D} $.

  \begin{lemma}\label{CovergenceInitialState}
  	For each $i\in\mathcal{H}$, we have $\varepsilon^{(n)}_{\mu_{i}}\overset{\rm u.c.p}\longrightarrow 0$  as $n\to\infty$.
  \end{lemma}
  \proof Let $\hat\varepsilon^{(n)}_{\mu_i}(t):= \mu^{(n)}_i(t)/n-\hat\mu^{(n)}_i(t)$ for $t\geq 0$. We can write $\varepsilon^{(n)}_{\mu_{i}}$ under the form
  \beqnn
  \varepsilon^{(n)}_{\mu_{i}}(t)= \hat\mu^{(n)}_i(nt) + R^{(n)}_{ii}* \hat\mu^{(n)}_i(nt) -Z^{(n)}(0)e^{\frac{b_{ii}}{\sigma_i}t} + \hat\varepsilon^{(n)}_{\mu_i}(nt)+   R^{(n)}_{ii}* \hat\varepsilon^{(n)}_{\mu_i}(nt).
  \eeqnn
  By (\ref{UpperBound.Rii}) and Condition~\ref{MomentConditionInitialState}, we have $ \| \hat\varepsilon^{(n)}_{\mu_i}\|_{L^{1,\infty}}\overset{\rm d}\to 0$ and $\| R^{(n)}_{ii}* \hat\varepsilon^{(n)}_{\mu_i}\|_{L^\infty} \leq C \|\hat\varepsilon^{(n)}_{\mu_i}\|_{L^1}\overset{\rm d}\to 0 $ as $n\to\infty$.
  By (\ref{bar_mu}) with $\beta =0$, 
  \beqnn
  \hat\mu^{(n)}_i(nt) + R^{(n)}_{ii}* \hat\mu^{(n)}_i(nt) = Z^{(n)}(0)\cdot n(1-\|\phi^{(n)}_{ii}\|_{L^1})    \int_t^\infty R^{(n)}_{ii}(ns)ds.
  \eeqnn
  By (\ref{ConvergenceFourierRii}) with $\beta=0$ and Condition~\ref{Con.PhiH}, we have 
  \beqnn
  \int_t^\infty R^{(n)}_{ii}(ns)ds \overset{\rm u.c.}\to \int_t^\infty \sigma_i^{-1}e^{b_{ii}/\sigma_i\cdot s}ds =\frac{1}{b_{ii}}\cdot e^{b_{ii}/\sigma_i\cdot t}
  \eeqnn and hence $ \hat\mu^{(n)}_i(nt) + R^{(n)}_{ii}* \hat\mu^{(n)}_i(nt) -Z^{(n)}(0)e^{\frac{b_{ii}}{\sigma_i}t}\overset{\rm u.c.p}\longrightarrow 0$.
  The desired result follows by putting these estimates together.
  \qed

  \begin{lemma}\label{MomentEstimate}
  	There exist constants $C,\vartheta>0$ such that for any $t\geq 0$ and $n\geq 1$,
  	\beqnn
  	\mathbf{E}\big[ \big|Z_{\mathcal{H}}^{(n)}(t) \big|^{2\alpha} \big]
  	\leq C e^{\vartheta t}. 
  	\eeqnn
  \end{lemma}
  \proof Obviously, it suffices to prove $\mathbf{E}[|Z_{\vartheta,\mathcal{H}}^{(n)}(t)|^{2\alpha} ] \leq C  $ for any $t\geq 0$ and $n\geq 1$. 
  Taking expectations on the both sides of (\ref{AdjustScaledDensity01}) with $\beta=\vartheta$,  we have
  \beqnn
  \mathbf{E}\big[ Z_{\vartheta,i}^{(n)}(t) \big]
  \ar=\ar \mathbf{E}\Big[\frac{ \mu^{(n)}_{\vartheta,i}(nt)}{n}\Big] +   \int_0^{nt}  R^{(n)}_{\vartheta,ii}(nt-s) \cdot \mathbf{E}\Big[\frac{\mu_{\vartheta,i}^{(n)}(s)}{n}\Big] ds \cr
  \ar\ar + \int_0^t R_{\vartheta,iI}^{(n)}(ns) ds + \sum_{j\in \mathcal{H}_i}\int_0^t nR_{\vartheta,ij}^{(n)}(n(t-s))  \mathbf{E}\big[ Z_{\vartheta,j}^{(n)}(s)\big]ds.
  \eeqnn
  From Condition~\ref{MomentConditionInitialState} and (\ref{UpperBound.Rii}), the first two terms on the right side of this equality are uniformly bounded.
  Moreover, by Young's convolution inequality,
  \beqnn
  \sup_{t\geq 0}\int_0^t nR_{\vartheta,ij}^{(n)}\big(n(t-s)\big) \cdot \mathbf{E}\big[ Z_{\vartheta,j}^{(n)}(s)\big]ds
  \ar\leq\ar \big\|R^{(n)}_{\vartheta,ij}\big\|_{L^1}\cdot \sup_{s\geq 0} \mathbf{E}\big[ Z_{\vartheta,j}^{(n)}(s)\big].
  \eeqnn
  From Condition~\ref{Con.PhiH} and (\ref{L1NormRii}), we have as $n\to\infty$,
  \beqlb\label{ConvergenceRtheta}
  \big\|R^{(n)}_{\vartheta,ij}\big\|_{L^1} \ar=\ar \frac{n\|\phi^{(n)}_{\vartheta,ij}\|_{L^1}}{n\big(1-\|\phi^{(n)}_{\vartheta,ii}\|_{L^1}\big)} \to \frac{b_{ij}}{\sigma_i\vartheta - b_{ii}}<\infty.
  \eeqlb
  Choosing $\vartheta$ large enough such that $\|R^{(n)}_{\vartheta,ij}\|_{L^1}\leq \frac{1}{2d}$ for any $n\geq 1$, we have
  \beqnn
  \sup_{t\geq 0}  \mathbf{E}\big[ Z_{\vartheta,i}^{(n)}(t)\big]\leq C + \frac{1}{2d} \sum_{j\in\mathcal{H}_i}\sup_{s\geq 0} \mathbf{E}\big[ Z_{\vartheta,j}^{(n)}(s)\big]
  \eeqnn
  and hence
  \beqnn
  \sum_{i\in\mathcal{H}}\sup_{t\geq 0}  \mathbf{E}\big[ Z_{\vartheta,i}^{(n)}(t)\big]\leq C + \frac{1}{2} \sum_{j\in\mathcal{H}}\sup_{s\geq 0} \mathbf{E}\big[ Z_{\vartheta,j}^{(n)}(s)\big],
  \eeqnn
  which immediately induces that $\sup_{t\geq 0}  \mathbf{E}[ |Z_{\vartheta,\mathcal{H}}^{(n)}(t)|] \leq C $. 
  Here the constant $C$ is independent of $n$ and $t$.
  We now start to give an upper bound for the second moment.
  Squaring both sides of (\ref{AdjustScaledDensity01}), using the Cauchy-Schwarz inequality and then taking expectations, we have
 \beqnn
 \mathbf{E}\big[ |Z_{\vartheta,i}^{(n)}(t)|^2 \big]
 \ar\leq\ar 2^{2d+2} \mathbf{E}\Big[\Big|\frac{ \mu^{(n)}_{\vartheta,i}(nt)}{n}\Big|^2 \Big] +2^{2d+2} \Big|\int_0^t R_{\vartheta,iI}^{(n)}(ns) ds\Big|^2 \cr
 \ar\ar
 +  2^{2d+2} \mathbf{E}\Big[\Big| \int_0^{nt}  R^{(n)}_{\vartheta,ii}(nt-s) \cdot \frac{\mu_{\vartheta,i}^{(n)}(s)}{n}ds \Big|^2 \Big]\cr
 \ar\ar  +\sum_{j\in \mathcal{H}_i} 2^{2d+2} \mathbf{E}\Big[\Big|\int_0^t  n  R_{\vartheta,ij}^{(n)}(n(t-s)) Z_{\vartheta,j}^{(n)}(s) ds \Big|^2 \Big] \cr
  \ar\ar +\sum_{j=1}^d 2^{2d+2} \mathbf{E}\Big[\Big| \int_0^t \int_{\mathbb{U}} R_{\vartheta,i}^{(n)}(n(t-s),u)\frac{e^{-\vartheta s}}{n} \tilde{N}_j^{(n)}(n\cdot ds,du )\Big|^2 \Big].
  \eeqnn
  Like the previous argument, we can prove that the first and third terms on the right side of this inequality are uniformly bounded.
  Notice that the integrand in the stochastic integral satisfies that for any $u\in\mathbb{U}$,
  \beqlb\label{UpperBoundRi}
  \sup_{n\geq 1} \big\|R_{\vartheta,i}^{(n)}(u)\big\|_{L^\infty} \leq \big\|\phi_i(u)\big\|_{\rm TV} + C \cdot \big\|\phi_i(u)\big\|_{L^1} \leq C \cdot \Phi(u).
  \eeqlb
  By the Burkholder-Davis-Gundy inequality and the uniform bound for the first moment of $\{Z_{\vartheta,\mathcal{H}}^{(n)}\}_{n\geq 1}$,
  \beqnn
  \lefteqn{\mathbf{E}\Big[\Big| \int_0^t \int_{\mathbb{U}} R_{\vartheta,i}^{(n)}\big(n(t-s),u\big)\frac{e^{-\vartheta s}}{n} \tilde{N}_j^{(n)}(n\cdot ds,du )\Big|^2 \Big]}\ar\ar\cr
  \ar\leq\ar C\int_0^t e^{-\vartheta s} \mathbf{E}\big[Z_{\vartheta,j}^{(n)}(s)\big] ds  \int_{\mathbb{U}}\big|\Phi(u)\big|^2 \nu_j^{(n)}(du) \leq C.
  \eeqnn
  Here the constant $C>0$ is independent of $n$ and $t$. 
  Moreover, by H\"older's inequality and then Young's convolution inequality,
  \beqnn
  \lefteqn{\sup_{t\geq 0}\mathbf{E}\Big[\Big|\int_0^t  n  R_{\vartheta,ij}^{(n)}\big(n(t-s)\big) Z_{\vartheta,j}^{(n)}(s) ds \Big|^2 \Big]}\ar\ar\cr
  \ar\leq\ar \big\|R_{\vartheta,ij}^{(n)}\big\|_{L^1} \cdot\sup_{t\geq 0} \int_0^t  n  R_{\vartheta,ij}^{(n)}\big(n(t-s)\big) \mathbf{E}\big[\big|Z_{\vartheta,j}^{(n)}(s)\big|^2\big] ds\cr
  \ar\leq\ar \big\|R_{\vartheta,ij}^{(n)}\big\|_{L^1}^2\cdot \sup_{t\geq 0}\mathbf{E}\big[\big|Z_{\vartheta,j}^{(n)}(t)\big|^2\big].
  \eeqnn
  Putting all estimates above together, we also have
  \beqnn
  \sup_{t\geq 0}\mathbf{E}\big[ \big|Z_{\vartheta,i}^{(n)}(t)\big|^2 \big]
  \ar\leq\ar C + \sum_{j\in \mathcal{H}_i} 2^{2d+2}  \big\|R_{\vartheta,ij}^{(n)}\big\|_{L^1}^2 \cdot \sup_{t\geq 0}\mathbf{E}\big[\big|Z_{\vartheta,j}^{(n)}(t)\big|^2\big].
  \eeqnn
  From (\ref{ConvergenceRtheta}), we  choose $\vartheta>0$ large enough such that $2^{2d+2}\big\|R_{\vartheta,ij}^{(n)}\big\|_{L^1}^2 \leq \frac{1}{2d}$ and then
  \beqnn
  \sum_{i=1}^d\sup_{t\geq 0}\mathbf{E}\big[ \big|Z_{\vartheta,i}^{(n)}(t)\big|^2 \big]
  \ar\leq\ar C + \frac{1}{2}\sum_{j=1}^d    \sup_{t\geq 0}\mathbf{E}\big[\big|Z_{\vartheta,j}^{(n)}(t)\big|^2\big],
  \eeqnn
  These induce that $\sup_{t\geq 0}  \mathbf{E}[ |Z_{\vartheta,\mathcal{H}}^{(n)}(t)|^2] \leq C $ with the constant $C$ independent of $n$ and $t$.
  Similarly, we also can prove that  for some  $\vartheta>0$,
  \beqnn
   \sup_{n\geq 1}\sup_{t\geq 0}  \mathbf{E}\big[ \big|Z_{\vartheta,\mathcal{H}}^{(n)}(t)\big|^{2\alpha}\big] \leq C . 
  \eeqnn 
  \qed

  \begin{lemma}\label{Error0ij}
  	For $i,j\in\mathcal{H}$ with $i\neq j$, we have $\tilde\varepsilon_{iI}^{(n)} \overset{\rm u.c.}\to 0$ and $\tilde\varepsilon_{ij}^{(n)}\overset{\rm u.c.p.} \longrightarrow 0$  as $n\to\infty$ .
  \end{lemma}
  \proof The first convergence follows directly from Lemma~\ref{L1ConvergenceRij}. For the second one, by H\"older's inequality we have for any $T>0$,
  \beqlb\label{ProofE0ij}
  \sup_{t\in[0,T]}\big|\tilde\varepsilon_{ij}^{(n)}(t)\big|^{\alpha}
  \leq \big\|\varepsilon_{R_{ij}}^{(n)}\big\|_{L^1_1}^{\alpha-1}\cdot \sup_{t\in[0,T]} \big|\varepsilon_{R_{ij}}^{(n)}\big|* \big|Z_j^{(n)}\big|^\alpha (t).
  \eeqlb
  By Young's inequality and Lemma~\ref{MomentEstimate},
  \beqnn
  \mathbf{E}\Big[\Big|\sup_{t\in[0,T]} \big|\varepsilon_{R_{ij}}^{(n)}\big|* \big|Z_j^{(n)}\big|^\alpha (t)\Big|^2 \Big] \ar\leq\ar  \big\|\varepsilon_{R_{ij}}^{(n)}\big\|_{L_{2T}^2}^2 \cdot \int_0^{2T} \mathbf{E}\big[\big|Z_j^{(n)}(s)\big|^{2\alpha} \big]ds \leq C \big\|\varepsilon_{R_{ij}}^{(n)}\big\|_{L_{2T}^2}^2 .
  \eeqnn
  From (\ref{ResovlentEqn01}), (\ref{UpperBound.Rii}) and Condition~\ref{Con.PhiH},  there exists a constant $C>0$ independent of $n$ such that
  \beqnn
  \int_0^{2T} \big|nR^{(n)}_{ij}(ns)\big|^2ds
  \ar\leq\ar  C n\int_0^\infty \big|\phi^{(n)}_{ij}(s)\big|^2ds  + C  \int_0^{2T}\big| n R^{(n)}_{ii}* \phi^{(n)}_{ij}(nt)\big|^2 dt \cr
  \ar\ar\cr
  \ar\leq\ar C n \big\|\phi^{(n)}_{ij}\big\|_{L^1}  \cdot \big\|\phi^{(n)}_{ij}\big\|_{\rm TV}  +  Cn^2 \big\|\phi^{(n)}_{ij}\big\|_{L^1}^2 <C
  \eeqnn
  and hence $\sup_{n\geq 1}\|\varepsilon_{R_{ij}}^{(n)}\|_{L_{2T}^2} <\infty $. Taking this back into (\ref{ProofE0ij}), from Lemma~\ref{L1ConvergenceRij} we have as $n\to\infty$,
  \beqnn
  \mathbf{E}\Big[\sup_{t\in[0,T]}\big|\tilde\varepsilon_{ij}^{(n)}(t)\big|^{2\alpha} \Big] 
  \leq C \big\|\varepsilon_{R_{ij}}^{(n)}\big\|_{L^1_T}^{2\alpha-2} \to 0
  \eeqnn
  and this proof is completed.
  \qed
  
  For each $i\in\mathcal{H}$ and $j\in\mathcal{D}$, we now start to prove the error $\varepsilon^{(n)}_{ij}$ vanishes as $n\to\infty$. 
  Like the proof of Lemma 5.7 in \cite{HX2019a}, we can prove the following moment estimate for the stochastic integral driven by MHPI-measures  by using the Burkholder-Davis-Gundy inequality.
  
  \begin{proposition} \label{PropBDG}
  	For any $T>0$, there exists a constant $C>0$ such that for any $\kappa\in[1,\alpha]$, $i\in\mathcal{D}$, $r,h\in[0,T]$ and measurable function $f(t,s,u)$  defined on $\mathbb{R}_+^2\times\mathbb{U}$,
  	\beqlb\label{BDG}
  	\lefteqn{\mathbf{E}  \Big[
  	\Big|	\int_r^{r+h} \int_{\mathbb{U} }   f(t,s,u)  \tilde{N}_{i}^{(n)}(n\cdot ds, du) \Big|^{2\kappa}
  	\Big]}\ar\ar\cr
  	\ar\leq\ar C  n^2\int_{\mathbb{U} }  \nu_{i}^{(n)}(du)
  	\int_r^{r+h} |f(t,s,u) |^{2\kappa}ds \cr
  	\ar\ar + C  \Big| n^2 \int_{\mathbb{U} }  \nu_{i}^{(n)}(du) 	\int_r^{r+h}   |f(t,s,u) |^2 ds\Big|^{\kappa}.
  	\eeqlb
  \end{proposition}
  \begin{corollary}\label{ConverFDD}
  	For each $i\in\mathcal{H}$ and $j\in\mathcal{D}$, we have $\varepsilon^{(n)}_{ij}\overset{\rm f.d.d.}\longrightarrow 0$ as $n\to\infty$.
  \end{corollary}
  \proof
  Applying (\ref{BDG}) together with   Lemma~\ref{L2ConvergenceRi} and \ref{H1} to $\mathbf{E} \big[|\varepsilon^{(n)}_{ij}(t)|^2\big]$  ,  we have for any $T>0$,
  \beqnn
  \sup_{t\in[0,T]}\mathbf{E} \big[\big|\varepsilon^{(n)}_{ij}(t)\big|^2\big]
  \leq C \int_{\mathbb{U}}  \big\|\varepsilon_{R_i}^{(n)}(u)\big\|_{L^2}^2  \nu_j^{(n)}(du)
  \leq C \big|\epsilon_n\big|^2 \cdot \int_{\mathbb{U}}  \big|\Phi(u)\big|^2  \nu_{ii}^{(n)}(du ),
  \eeqnn
  which goes to $0$ as $n\to\infty$ and the desired result follows.
  \qed

  We now start to prove the tightness of the sequence $\{\varepsilon^{(n)}_{ii} \}_{n\geq 1}$.
  The tightness of other sequences can be proved similarly.
  By  Corollary~3.33 in \cite[p.353]{JS03} and the definition of $\varepsilon^{(n)}_{ii}$, it suffices to prove that  $\{I_{ii}^{(n)}\}_{n\geq 1}$ is tight and $\{ J_{ii}^{(n)}\}_{n\geq 1}$ is $C$-tight, where
  \beqlb
  I_{ii}^{(n)}(t)\ar:=\ar  \int_0^t \int_{\mathbb{U}} \frac{R_i^{(n)}(n(t-s),u)}{n}  \tilde{N}_{i}^{(n)}(n\cdot ds,du),\label{ErrorIii}\\
  J_{ii}^{(n)}(t)\ar:=\ar \int_0^t \int_{\mathbb{U}} \frac{\|\phi_i(u)\|_{L^1}}{\sigma_i} \frac{ e^{\frac{b_{ii}}{\sigma_i}(t-s)}}{n}\tilde{N}_{i}^{(n)}(n\cdot ds,du), \label{ErrorJii}
  \eeqlb
  From the fact that $\exp\{\frac{b_{ii}}{\sigma_i}(t-s)\}= 1+ \frac{b_{ii}}{\sigma_i}\int_s^t \exp\{\frac{b_{ii}}{\sigma_i}(r-s)\}dr$, we can write $J_{ii}^{(n)} $ as
  \beqnn
  J_{ii}^{(n)} (t)=\frac{b_{ii}}{\sigma_i}\int_0^{t}J_{ii}^{(n)} (s)ds - \int_0^{t} \int_{\mathbb{U}}  \frac{ \|\phi_i(u)\|_{L^1} }{\sigma_i n}\tilde{N}_{i}^{(n)}(n\cdot ds,du), \quad t\geq 0.
  \eeqnn
  Obviously, $J_{ii}^{(n)} $ is an $(\mathscr{F}_{nt})$-semimartingale.

  \begin{proposition}\label{TightJii} \it
  	The sequence $\{ J_{ii}^{(n)}\}_{n\geq 1} $ is $C$-tight.
  \end{proposition}
  \proof
  As a preparation, we firstly give some moment estimates for $J_{ii}^{(n)}$. The exists a constant $C>0$ such that for any $n\geq 1$ and $T>0$,
  \beqnn
  \mathbf{E}\Big[  \sup_{t\in[0,T]} \big|  J_{ii}^{(n)}(t)\big|^{2}\Big]
  \ar\leq\ar C T\int_0^T  \mathbf{E}\Big[  \sup_{s\in[0,t]} \big|  J_{ii}^{(n)}(s)\big|^{2}\Big]dt \cr
  \ar\ar + C \mathbf{E}\Big[ \int_0^T \int_{\mathbb{U}}  \frac{ \|\phi_i(u)\|_{L^1}^2 }{\sigma_i^2 n^2}N_{i}^{(n)}(n\cdot ds,du) \Big]\cr
  \ar\leq\ar CT +C T\int_0^T  \mathbf{E}\Big[  \sup_{s\in[0,t]} \big|  J_{ii}^{(n)}(s)\big|^{2}\Big]dt.
  \eeqnn
  Here the first inequality follows from H\"older's inequality and (\ref{BDG}). 
  The second one follows from Lemma~\ref{MomentEstimate} and \ref{H1}.
  By Gronwall's inequality,
  \beqlb\label{MomentJ01}
  \sup_{n\geq 1}\mathbf{E}\Big[ \sup_{t\in[0,T]} \big|  J_{ii}^{(n)}(t) \big|^2\Big]<\infty.
  \eeqlb
  We now prove the tightness of $\big\{ J_{ii}^{(n)}\big\}_{n\geq 1} $. For any bounded stopping time $\tau\leq T$  and $h\in(0,1)$, we have
  \beqlb\label{MomentDifference}
  \mathbf{E}\big[\big|\Delta_h J_{ii}^{(n)} (\tau)\big|^2\big]\ar\leq\ar C \mathbf{E}\Big[\Big|\int_{\tau}^{\tau+h}J_{ii}^{(n)} (s)ds\Big|^2\Big] \cr
  \ar\ar + C\mathbf{E}\Big[\Big|\int_{\tau}^{\tau+h}\int_{\mathbb{U}}  \frac{ \|\phi_i(u)\|_{L^1} }{\sigma_i n}\tilde{N}_{i}^{(n)}(n\cdot ds,du)\Big|^2\Big].
  \eeqlb
  From H\"older's inequality and (\ref{MomentJ01}), we have
  \beqnn
  \mathbf{E}\Big[\Big|\int_{\tau}^{\tau+h}J_{ii}^{(n)} (s)ds\Big|^2\Big] 
  \ar\leq\ar h \mathbf{E}\Big[\int_{\tau}^{\tau+h}\big|J_{ii}^{(n)} (s)\big|^2ds \Big] \cr 
  \ar\leq\ar h \int_0^{T+1}\mathbf{E}\big[ \big|J_{ii}^{(n)} (s)\big|^2 \big]ds\leq Ch.
  \eeqnn
  Moreover, applying (\ref{BDG}) again to the last expectation in (\ref{MomentDifference}), it can be bounded by
  \beqnn
  \lefteqn{C \mathbf{E}\Big[ \int_{\tau}^{\tau+h}\int_{\mathbb{U}}  \frac{ \|\phi_i(u)\|^2_{L^1} }{ n^2}N_{i}^{(n)}(n\cdot ds,du) \Big]}\ar\ar\cr
  \ar\leq\ar C \mathbf{E}\Big[ \int_{\tau}^{\tau+h} Z_i^{(n)}(s)ds\Big] \cdot \int_{\mathbb{U}}  \big\|\phi_i(u)\big\|^2_{L^1}  \nu_{i}^{(n)}(du).
  \eeqnn
  Applying H\"older's inequality, Jensen's inequality and then using Lemma~\ref{MomentEstimate}, we also have
  \beqnn
  \mathbf{E}\Big[ \int_{\tau}^{\tau+h} Z_i^{(n)}(s)ds\Big] \ar\leq\ar \sqrt{h} \mathbf{E}\Big[ \Big(\int_{\tau}^{\tau+h} \big|Z_i^{(n)}(s)\big|^2ds\Big)^{1/2}\Big]\cr
  \ar\leq\ar  \sqrt{h} \Big(\int_0^{T+1} \mathbf{E}\big[\big|Z_i^{(n)}(s)\big|^2\big]ds\Big)^{1/2}
  \leq C\sqrt{h}.
  \eeqnn
  Putting all estimates above together, we have $  \mathbf{E}[|\Delta_h J_{ii}^{(n)} (\tau)|^2] \leq C \sqrt{h}$
  with the constant $C$ independent of $n$ and $\tau$.
  The criterion of Aldous; see \cite{Aldous1978}, yields the tightness of $\{J_{ii}^{(n)} \}_{n\geq 1}$ directly.
  
  It remains to prove the continuity of cluster points.
  For any cluster point $J^*_{ii}$, it suffices to prove that $ \mathbf{E}[\sum_{t\in[0,T]}|\Delta_-J^*_{ii}(t)|^{2\alpha}  ]=0$ for any $T>0$.
  There exists a subsequence of $\{J_{ii}^{(n)}\}$, still denoted by itself, such that $J_{ii}^{(n)} \overset{\rm d}\to J^*_{ii}$ in $\mathbf{D}([0,\infty),\mathbb{R})$. 
  For $\epsilon>0$, let $g_\epsilon $ be a  continuous function on $\mathbb{R}$ vanishing in a neighborhood of $0$ and satisfying that $g_\epsilon(x) $ increases to $ |x|^{2\alpha}$ as $\epsilon\to 0$ for any $x\in\mathbb{R}$. By the monotone convergence theorem and then Proposition~3.16 in \cite[p.349]{JS03},
  	\beqnn
  	\mathbf{E}\Big[\sum_{t\in[0,T]} \big|\Delta_- J^*_{ii}(t)\big|^{2\alpha}  \Big] 
  	\ar= \ar \lim_{\epsilon \to 0+}\mathbf{E}\Big[\sum_{t\in[0,T]} g_\epsilon\big(\Delta_- J^*_{ii}(t)\big)\Big] \cr
  	\ar=\ar \lim_{\epsilon \to 0+} \lim_{n\to\infty}\mathbf{E}\Big[\sum_{t\in[0,T]} g_\epsilon\big(\Delta_- J^{(n)}_{ii}(t)\big)\Big] \cr
  	\ar\leq\ar \lim_{n\to\infty}\mathbf{E}\Big[\sum_{t\in[0,T]} \big|\Delta_- J^{(n)}_{ii}(t)\big|^{2\alpha}\Big] .
  	\eeqnn 
  By the hypothesis \ref{H1} and the properties of stochastic integrals with respect to a random point measure,
  \beqnn
  \mathbf{E}\Big[\sum_{t\in[0,T]}\big|\Delta_- J^{(n)}_{ii}(t) \big|^{2\alpha}  \Big]\ar=\ar \mathbf{E}\Big[ \int_0^T \int_{\mathbb{U}}  \frac{ \|\phi_i(u)\|_{L^1}^{2\alpha} }{|\sigma_i n|^{2\alpha}} N_{i}^{(n)}(n\cdot ds,du)  \Big]\leq \frac{C}{n^{2\alpha-2}} ,
  \eeqnn
  which goes to $0$ as $n\to\infty$ and hence the sequence $\{J^{(n)}_{ii}\}_{n\geq 1}$ is $C$-tight.
  \qed

  We now start to prove the tightness of $\{I_{ii}^{(n)} \}_{n\geq 1} $.
  For some  $\theta>2$, let $I_{ii,\theta}^{(n)}$ be a linear interpolation of  $I_{ii}^{(n)}$ defined as follows
  \beqlb\label{LinearInterSmallError}
  I_{ii,\theta}^{(n)}(t)\ar:=\ar I_{ii}^{(n)}\Big(\frac{[n^{\theta} t]}{n^{\theta}}\Big)
  +\Big(n^{\theta}t-
  [n^{\theta}t]\Big)
  \Big[I_{ii}^{(n)}\Big(\frac{[n^{\theta} t]+1}{n^{\theta}}\Big)
  -I_{ii}^{(n)}\Big(\frac{[n^{\theta} t]}{n^{\theta}}\Big)\Big],
  \eeqlb
  for $ t\geq 0$. 
  We now start to prove that the sequence $\{I_{ii,\theta}\}_{n\geq 1}$ is tight and a good approximation for $\{I^{(n)}_{ii}\}_{n\geq 1}$. This will induce the tightness of $\{I^{(n)}_{ii}\}_{n\geq 1}$ immediately. 
  As a preparation, the next proposition gives some upper bound estimates for the shifted resolvent.

  \begin{proposition}\label{ResolventDifferenceL2}
  	For any $\kappa\geq 1$,	there exists a constant $C>0$ such that for any $u\in\mathbb{U}$ and $h\in[0,1]$,
  	\beqlb\label{BoundDiff}
  	\sup_{n\geq 1}\int_0^\infty \big| \Delta_{nh}  R_i^{(n)}(nt,u)\big|^{2\kappa} dt \leq C \big|\Phi(u)\big|^{2\kappa}\cdot h.
  	\eeqlb
  \end{proposition}
  \proof  We first prove this result with $\kappa=1$. By the Fourier isometry,
  \beqnn
  \int_0^\infty \big| \Delta_{n h}  R_i^{(n)}(nt,u)\big|^2 dt= \int_{\mathbb{R}}
  \Big|\big(e^{\mathtt{i}\lambda  h }-1\big)\int_{\mathbb{R}}e^{\mathtt{i}\lambda t}  R_i^{(n)}(nt,u)dt\Big|^2 d\lambda.
  \eeqnn
  Similarly as in the proof of Lemma~\ref{L2ConvergenceRii}, we can prove that for any $u\in\mathbb{U}$ and $\lambda\in\mathbb{R}$,
  \beqnn
  \Big|\int_0^\infty e^{\mathrm{i}\lambda t} \phi^{(n)}_{i}(t,u) dt \Big| +
  \Big| \int_0^\infty e^{\mathrm{i}\lambda t} R^{(n)}_{i}(nt,u)dt\Big| \leq C\cdot \Phi(u) \Big(\frac{1}{|\lambda|}\wedge 1 \Big). 
  \eeqnn
  From this and the fact that $|e^{\mathtt{i}\lambda h }-1|\leq |\lambda h|\wedge 2$,
  \beqnn
  \int_0^\infty \big| \Delta_{nh}  R_i^{(n)}(nt,u)\big|^2 dt\leq C  \big|\Phi(u)\big|^2 \cdot \int_{\mathbb{R}}  \big(|\lambda h|^2\wedge 1 \big) \cdot  \Big(\frac{1}{|\lambda|^2}\wedge 1 \Big)d\lambda .
  \eeqnn
  A simple calculation deduces that the last integral can be bounded by $6h$ and hence the inequality (\ref{BoundDiff}) holds for $\kappa=1$.
  When $\kappa>1$, we have
  \beqnn
  \int_0^\infty \big| \Delta_{nh}  R_i^{(n)}(nt,u)\big|^{2\kappa} dt  \ar\leq\ar 2^{2\kappa-2} \big\|R_i^{(n)}(u)\big\|_{L^\infty}^{2\kappa-2}\int_0^\infty \big| \Delta_{nh}  R_i^{(n)}(nt,u)\big|^{2} dt.
  \eeqnn
  Then (\ref{BoundDiff}) with $\kappa>1$ follows directly from (\ref{UpperBoundRi}) and the previous result.
  \qed
  
  \begin{proposition}\label{WellApproximationIii}
  	We have $|I_{ii}^{(n)} - I_{ii,\theta}^{(n)}| \overset{\rm u.c.p.}\longrightarrow 0$ as $n\to\infty$.
  \end{proposition}
  \proof Here we just prove $| I_{ii}^{(n)} - I_{ii,\theta}^{(n)}| \overset{\rm p}\to 0$ uniformly on $[0,1]$.
  From the definition of $ I_{ii,\theta}$ and the triangle inequality, we have for any $k\geq 0$ and  $t\in [k/n^{\theta},(k+1)/n^{\theta}]$,
  \beqlb\label{UpperBoundApproximation01}
  \big|I_{ii}^{(n)} (t)- I_{ii,\theta}^{(n)}(t)\big| 
  \ar\leq\ar  \big|I_{ii}^{(n)} (t)-  I_{ii}^{(n)}(kn^{-\theta})\big|+ \big|I_{ii}^{(n)} (t)-  I_{ii}^{(n)}(kn^{-\theta})\big| \cr
  \ar\ar +  \big|I_{ii}^{(n)}((k+1)n^{-\theta})-I_{ii}^{(n)} (kn^{-\theta})\big|,
  \eeqlb
  which can be bounded by $3\sup_{  h\leq n^{-\theta}} \big|\Delta_h I_{ii}^{(n)} (kn^{-\theta}) \big|$
  and hence
  \beqlb\label{UpperBoundApproximation02}
  \sup_{t\in[0,1]}\big|I_{ii}^{(n)} (t)- I_{ii,\theta}^{(n)}(t)\big|
  \ar\leq\ar 3 \sup_{k=0,\cdots,[n^\theta]; h\leq n^{-\theta}} \big|\Delta_h I_{ii}^{(n)} (kn^{-\theta}) \big|.
  \eeqlb
  For any $u\in\mathbb{U}$, since $\|\phi_i(u)\|_{\rm TV}<\infty$, we have $\phi_i(t,u) =\phi_i^+(t,u)-\phi_i^-(t,u) $, where   $\phi_i^+(t,u)$ and $\phi_i^-(t,u)$ are two  non-negative, non-increasing functions\footnote{For any non-negative function $f$ on $\mathbb{R}_+$ with $\|f\|_{\rm TV}<\infty$, by the Jordan decomposition there exists two  non-negative, non-decreasing functions $f^+_0$ and $f^-_0$ on $\mathbb{R}_+$ such that $f=f^+_0-f^-_0$, $f^+_0\geq f^-_0$ and $\|f\|_{\rm TV}= f^+_0(\infty)-f^+_0(0)+ f^-_0(\infty)-f^-_0(0)<\infty$. 
  	Thus $f=f^+-f^- $, where $f^+ :=\|f\|_{\rm TV}-f^-_0  $ and $f^- :=\|f\|_{\rm TV}-f^+_0  $ are non-negative, non-increasing on $\mathbb{R}_+$. } 
  on $\mathbb{R}$ with $\phi_i^+(t,u)=\phi_i^-(t,u)=\phi_i^+(0,u)$ for $t< 0$.
  From (\ref{ErrorIii}), we have   $ |\Delta_h I_{ii}^{(n)} (t) | \leq  A^{(n)}_{1}(t,h) +  A^{(n)}_{2}(t,h) + A^{(n)}_{3}(t,h)+ A^{(n)}_{4}(t,h)$ for any $t,h\in[0,1]$ with
  \beqnn
  A^{(n)}_{1}(t,h) \ar:=\ar  n \int_{\mathbb{U}}     \nu_{ii}^{(n)}(du)   \int_0^{t+h} Z_i^{(n)}(s)  |\Delta_{nh}R_{i}^{(n)}(n(t-s),u)|ds , \\
  A^{(n)}_{2}(t,h)  \ar:=\ar   \int_0^{t}  \int_{\mathbb{U}}    \frac{ |\Delta_{nh}\big(R_{ii}^{(n)}* \phi_i(n(t-s),u)\big)| }{n}N_i^{(n)}(n\cdot ds,du) ,\\
  A^{(n)}_{3}(t,h)  \ar:=\ar \int_0^{t+h}  \int_{\mathbb{U}}    \frac{ |\Delta_{nh}\phi_i^+(n(t-s),u)|}{n}N_i^{(n)}(n\cdot ds,du),\\
  A^{(n)}_{4}(t,h)  \ar:=\ar \int_0^{t+h}  \int_{\mathbb{U}}   \frac{ |\Delta_{nh}\phi_i^-(n (t-s),u)|}{n}N_i^{(n)}(n\cdot ds,du) .
  \eeqnn 
  Thus it suffices to prove that  for any $\eta>0$ and $j\in\{1,2,3 ,4 \}$,
  \beqlb\label{ErrorI0}
  \lim_{n\to\infty} \mathbf{P} \Big(  \sup_{k=0,\cdots,[n^\theta]; h\leq n^{-\theta}} A^{(n)}_j(kn^{-\theta}, h)  \geq \eta \Big) = 0.
  \eeqlb
  In  the sequel of this proof, the constant $C>0$ is independent of $(n,t,u,h)$ and may vary from line to line. 
  
  {\it Step 1.} We first prove  (\ref{ErrorI0}) with $j=1$.
  By Young's convolution inequality and  Proposition~\ref{ResolventDifferenceL2},
  \beqnn
  \lefteqn{\sup_{t\in[0,1]} \Big| \int_0^{t+h} Z_i^{(n)}(s) \cdot  \big|\Delta_{nh}R_i^{(n)}(n(t-s),u)\big| ds  \Big|^2}\ar\ar\cr
  \ar \leq \ar  \int_0^2 |Z_i^{(n)}(r)|^2 dr  \cdot  \int_0^\infty   \big|\Delta_{nh}R_i^{(n)}(ns,u)\big|^{2}  ds  \cr
  \ar\leq\ar C  \int_0^2 \big|Z_i^{(n)}(s)\big|^2ds   \cdot \big|\Phi(u)\big|^2\cdot h
  \eeqnn
  and hence
  \beqnn
  \sup_{h\leq n^{-\theta}; t\in[0,1]}  \big| A^{(n)}_{1}(t,h)\big|^2 \leq C   \int_0^2 \big|Z^{(n)}(s)\big|^2ds  \cdot   n^{2-\theta} .
  \eeqnn
  From Chebyshev's inequality and Lemma~\ref{MomentEstimate},
  \beqnn
   \mathbf{P} \Big( \sup_{k=0,\cdots,[n^\theta]; h<n^{-\theta}}   A^{(n)}_{1}(kn^{-\theta},h)  \geq \eta \Big) \leq \frac{1}{\eta^2}  \mathbf{E} \Big[  \sup_{h\leq n^{-\theta}; t\in[0,1]}  \big| A^{(n)}_{1}(t,h)\big|^2 \Big] \leq \frac{C}{\eta^2}n^{2-\theta } ,
  \eeqnn
  which vanishes as $n\to\infty$ since $\theta>2$. \smallskip
  
  {\it Step 2.} We now prove (\ref{ErrorI0}) with $j=2$.	
  By (\ref{UpperBound.Rii}), we have for any $t,h\in[0,1]$,
  \beqnn
  \big|\Delta_{nh}\big(R_{ii}^{(n)}* \phi_i(nt,u)\big)\big|
  \ar\leq \ar    \int_0^{nt} R_{ii}^{(n)}(s)\big|\Delta_{nh}\phi_i(nt-s,u)\big|ds \cr
  \ar\ar + \int_{nt}^{n(t+h)} R_{ii}^{(n)}(s)\phi_i\big(n(t+h)-s,u\big)ds \cr
  \ar\leq\ar C\int_{ 0}^{nh}   \phi_i(s,u)ds
  + C\int_0^{nt}    |\phi_i(nh+s,u)-\phi_i(s,u)|ds.
  \eeqnn
  The first term on the right side of the last inequality can be bounded by $C \|\phi_i(u)\|_{\rm TV}\cdot nh$. 
  By the preceding decomposition of $\phi_i$,  the second term can be bounded by
  \beqlb\label{DifferencePhi}
  \int_0^{nt}\big[\phi_i^+(s,u)-\phi_i^+(nh+s,u)\big]ds +
  \int_0^{nt}\big[\phi_i^-(s,u)-\phi_i^-(nh+s,u)\big]ds ,
  \eeqlb
  which can be bounded by $4 \|\phi_i(u)\|_{\rm TV}\cdot nh$.
  Putting these estimates together, we have
  $
  |\Delta_{nh}(R_{ii}^{(n)}* \phi_i)(nt,u)| \leq C  \|\phi_i(u)\|_{\rm TV}\cdot n h $
  and hence 
  \beqnn
  \sup_{t\in[0,1]; h\leq n^{-\theta}} A^{(n)}_{2}(t,h)  \leq  \frac{C}{n^{\theta}} \int_0^2 \int_{\mathbb{U}}   \big\|\phi_i(u)\big\|_{\rm TV} N_i^{(n)}(n\cdot ds,du).
  \eeqnn
  By Chebyshev's inequality and hypothesis \ref{H1},
  \beqnn
  \mathbf{P} \Big(  \sup_{k=0,\cdots,[n^\theta]; h\leq n^{-\theta}} A^{(n)}_2(kn^{-\theta}, h)  \geq \eta \Big)
  \leq \frac{1}{\eta} \mathbf{E} \Big[  \sup_{t\in[0,1]; h\leq n^{-\theta}} A^{(n)}_2(t, h) \Big] 
  \leq \frac{C}{\eta}\cdot n^{2-\theta},
  \eeqnn
  which goes to $0$ as $n\to\infty$ since $\theta>2$.\smallskip
  
  {\it Step 3.} We now prove (\ref{ErrorI0}) with $j=3$.
  For the case of $j=4$, it can be proved in the same way.
  Notice that
  $\sup_{h\leq n^{-\theta}} A_3^{(n)}(t,h) = A^{(n)}_{3,1}(t)+ A^{(n)}_{3,2}(t) $  with
  \beqnn
  A^{(n)}_{3,1}(t)\ar:=\ar n \int_{\mathbb{U}}  \nu^{(n)}_i(du)  \int_0^{t+n^{-\theta} } Z_i^{(n)}(s) \cdot \big| \Delta_{n^{1-\theta}} \phi_i^+(n(t-s),u) \big| ds, \\
  A^{(n)}_{3,2}(t) \ar:=\ar  \int_0^{t+n^{-\theta}} \int_{\mathbb{U}} \frac{| \Delta_{ n^{1-\theta}} \phi_i^+(n(t-s),u) | }{n}\tilde{N}_i^{(n)}(n\cdot ds,du ).
  \eeqnn
  By Young's convolution inequality, we can bound $\sup_{t\in[0,1]} A^{(n)}_{3,1}(t)$ by
  \beqlb\label{BoundI321}
  n\int_{\mathbb{U}} \nu^{(n)}_i(du)  \cdot   \Big(\int_0^2\big|\Delta_{ n^{1-\theta}} \phi_i^+(ns,u)\big|^{2} ds\Big)^{ 1/2} \cdot \Big(\int_0^{2}  \big|Z_i^{(n)}(s)\big|^2 ds\Big)^{1/2} .
  \eeqlb
  Since  $\phi_i^+$ is non-increasing, we have
  \beqnn
  \int_0^2 \big|\Delta_{ n^{1-\theta}} \phi_i^+(ns,u)\big|^{2} ds
  \ar\leq\ar
  2\big\| \phi_i(u)\big\|_{\rm TV}\int_0^2\big[ \phi_i^+(ns,u)-  \phi_i^+\big(n(s+ n^{-\theta}),u\big)\big]ds \cr
  \ar=\ar 2\big\| \phi_i(u)\big\|_{\rm TV}\Big( \int_0^2 \phi_i^+(ns,u)ds - \int_{n^{-\theta}}^{2+ n^{-\theta}} \phi_i^+(ns,u) ds \Big) \cr
  \ar\leq\ar \frac{4}{ n^{\theta}}\cdot \big\| \phi_i(u)\big\|_{\rm TV}^2
  \eeqnn
  and hence by hypothesis \ref{H1},
  \beqnn
  \sup_{t\in[0,1]} \big|A^{(n)}_{3,1}(t)\big|^2 \ar\leq \ar \frac{C}{n^{\theta-2}} \cdot \int_0^{2}  \big|Z_i^{(n)}(s)\big|^2 ds
  .
  \eeqnn
  Applying Chebyshev's inequality again, we have as $n\to\infty$,
  \beqnn
  \mathbf{P}\Big(  \sup_{t\in[0,1]} A^{(n)}_{3,1}(t)\geq \eta\Big)
  \ar\leq\ar  \frac{C}{\eta^2} \mathbf{E} \Big[ \sup_{t\in[0,1]} \big|A^{(n)}_{3,1}(t)\big|^2 \Big] \leq \frac{C}{\eta^2}  \frac{1}{n^{\theta-2}} \to 0.
  \eeqnn
  Applying  (\ref{BDG}) to $ A^{(n)}_{3,2}(t) $, similarly as in (\ref{BoundI321}) we also have
  \beqnn
  \sup_{t\in[0,1]}\mathbf{E}\big[ \big|A^{(n)}_{3,2}(t)\big|^{2\alpha}\big]
  \ar\leq\ar
  C\Big( \int_{\mathbb{U}}  \nu^{(n)}_i(du)  \int_0^2 \big|\Delta_{n^{1-\theta}} \phi_i^+(ns,u) \big|^2 ds \Big)^{\alpha}\cr
  \ar\ar + C n^{2-2\alpha} \int_{\mathbb{U}}  \nu^{(n)}_i(du)   \int_0^2 \big|\Delta_{n^{1-\theta}} \phi_i^+(ns,u) \big|^{2\alpha}ds \cr
  \ar\ar\cr
  \ar\leq\ar  C \cdot \big(n^{-\alpha\theta} +n^{2-2\alpha- \theta}\big) .
  \eeqnn
  From this and Chebyshev's inequality,
  \beqnn
  \mathbf{P}\Big( \sup_{k=0,\cdots,[n^\theta]} \big|A^{(n)}_{3,2}(kn^{-\theta})\big|\geq\eta\Big) 
  \ar\leq\ar  \frac{1}{\eta^{2\alpha}}  \sum_{k=0}^{[n^\theta]} \mathbf{E}\big[ \big|A^{(n)}_{3,2}(kn^{-\theta}) \big|^{2\alpha}\big] \cr \ar\leq\ar  \frac{C}{\eta^{2\alpha}}  \big(n^{\theta(1-\alpha)} +n^{2-2\alpha}\big),
  \eeqnn
  which vanishes as $n\to\infty$ since $\alpha\in(1,2 )$.
  \qed

  \begin{proposition}\label{MomentIncre}
  	For any $T>0$, there exists a constant $C>0$  such that for any $h\in[0,1]$ and $n\geq 1$,
  	\beqnn
  	\sup_{t\in[0,T]}	\mathbf{E}\big[\big| \Delta_h I_{ii}^{(n)}(t) \big|^{2\alpha}\big] 
  	\leq  C \cdot  \big(n^{2-2\alpha} \cdot h + h^\alpha \big).
  	\eeqnn 	
  \end{proposition}
  \proof 
  Applying the inequality (\ref{BDG}) to $ \Delta_h I_{ii}^{(n)}(t) $ and then using Proposition~\ref{ResolventDifferenceL2},
  \beqnn
  \mathbf{E}\big[\big| \Delta_h I_{ii}^{(n)}(t) \big|^{2\alpha}\big] \ar\leq\ar
  C \Big| \int_{\mathbb{U}} \nu^{(n)}_i(du) \int_0^\infty  \big|\Delta_{nh} R_i^{(n)}(ns,u)\big|^2 ds \Big|^{\alpha}\cr
  \ar\ar + Cn^{2-2\alpha} \int_{\mathbb{U}} \nu_i^{(n)}(du) \int_0^\infty  \big|\Delta_{nh} R_i^{(n)}(ns,u)\big|^{2\alpha}  ds  \cr
  \ar\leq\ar C h^{\alpha}  \Big| \int_{\mathbb{U}} \big|\Phi(u)\big|^2 \nu^{(n)}_i(du) \Big|^{\alpha}
  + \frac{Ch}{n^{2\alpha-2}} \int_\mathbb{U} \big|\Phi(u)\big|^{2\alpha} \nu^{(n)}_i(du)
  \eeqnn
  and the desired result follows directly from hypothesis \ref{H1}.
  \qed
  
  \begin{proposition}\label{TightnessIii}
  	The sequence $\big\{ I^{(n)}_{ii,\theta}\big\}_{n\geq1}$ is tight.
  \end{proposition}
  \proof From Proposition~10.3 in \cite[p.149]{EthierKurtz2005},
  it suffices to prove that there exist constants $C>0$ and $\epsilon\in(0,(2\alpha-2)/\theta)$ such that for any $t,h\in[0,1]$,
  \beqnn
  \sup_{n\geq 1}\mathbf{E}\Big[\big|\Delta_h I^{(n)}_{ii,\theta}(t)\big|^{2\alpha}\Big]\leq  C \cdot h^{1+\epsilon}.
  \eeqnn 
  If $jn^{-\theta}\leq t<t+h\leq (j+1)n^{-\theta} $ for some $j\geq 0$, from (\ref{LinearInterSmallError}) and Proposition~\ref{MomentIncre} we have $\Delta_h I^{(n)}_{ii,\theta}(t) = n^\theta h\cdot \Delta_{n^{-\theta}} I^{(n)}_{ii}(jn^{-\theta})$ and
  \beqnn
  \mathbf{E}\Big[\big|\Delta_h I^{(n)}_{ii,\theta}(t)\big|^{2\alpha}\Big]
  \ar=\ar n^{2\alpha\theta}h^{2\alpha} \mathbf{E}\Big[\big|\Delta_{n^{-\theta}} I^{(n)}_{ii}(jn^{-\theta})\big|^{2\alpha}\Big] \cr
  \ar\ar\cr
  \ar\leq\ar C   n^{2\alpha\theta-2\alpha-\theta+2}h^{2\alpha} \leq C h^{1+\epsilon}.
  \eeqnn
  Here the constant $C>0$ is independent of $n$ and $h$. 
  Similarly, if  $jn^{-\theta}\leq t \leq (j+1) n^{-\theta}\leq t+h\leq (j+2) n^{-\theta}$ for some $j\geq 0$, we also have
  \beqnn
  \mathbf{E}\Big[\big|\Delta_h I^{(n)}_{ii,\theta}(t)\big|^{2\alpha}\Big]
  \ar\leq\ar C\mathbf{E}\Big[|I^{(n)}_{ii,\theta}(t+h)-I^{(n)}_{ii,\theta}\big((i+1) n^{-\theta}\big)|^{2\alpha}\Big]\cr
  \ar\ar\cr
  \ar\ar + C\mathbf{E}\Big[\big|I^{(n)}_{ii,\theta}\big((i+1)  n^{-\theta}\big)-I^{(n)}_{ii,\theta}(t)\big|^{2\alpha}\Big]
  \leq C\cdot h^{1+\epsilon}.
  \eeqnn
  Finally, if $jn^{-\theta}\leq t \leq (j+1) n^{-\theta} $ and $l n^{-\theta}\leq t+h\leq (l+1) n^{-\theta}$ for some $j<l$, we have
  \beqnn
  \mathbf{E}\Big[\big|\Delta_hI^{(n)}_{ii,\theta}(t)\big|^{2\alpha}\Big]
  \ar\leq\ar
  C\mathbf{E}\Big[\big|I^{(n)}_{ii,\theta}(t+h)-I^{(n)}_{ii,\theta}(  l\cdot n^{-\theta})\big|^{2\alpha}\Big] \cr 
  \ar\ar +C\mathbf{E}\Big[\big|I^{(n)}_{ii,\theta}\big( (j+1) n^{-\theta}\big)-I^{(n)}_{ii,\theta}(t)\big|^{2\alpha}\Big]\cr
  \ar\ar +C\mathbf{E}\Big[\big|I^{(n)}_{ii,\theta}( l \cdot n^{-\theta})-I^{(n)}_{ii,\theta}\big( (j+1) n^{-\theta}\big)\big|^{2\alpha}\Big] .
  \eeqnn
  From the foregoing two results, the first two terms on the right side of this inequality can be bounded by $C\cdot h^{1+\epsilon}$.
  For the third term, notice that 
  $$I^{(n)}_{ii,\theta}( l n^{-\theta})-I^{(n)}_{ii,\theta}\big( (j+1) n^{-\theta}\big)=I^{(n)}_{ii }( l n^{-\theta})-I^{(n)}_{ii }\big( (j+1) n^{-\theta}\big),$$
  by Proposition~\ref{MomentIncre} it can be bounded by $C (h^{\alpha} +n^{2-2\alpha}h)\leq C h^{1+\epsilon}.$
  Here we have finished the proof.
  \qed
  
  We now summarize the results in Corollary~\ref{ConverFDD}, Proposition~\ref{TightJii}, \ref{WellApproximationIii} and \ref{TightnessIii} to get the weak convergence of $\{\varepsilon^{(n)}_{ij}\}_{n\geq 1}$ to $0$.

  \begin{lemma}
  	For each $i\in\mathcal{H}$ and $j\in\mathcal{D}$, we have $\varepsilon^{(n)}_{ij}\overset{\rm d}\to 0$ in $\mathbf{D}([0,\infty),\mathbb{R})$ as $n\to\infty$.
  \end{lemma}

  \subsubsection{Weak convergence of $(M^{(n)}_{ij})_{i\in\mathcal{H},j\in\mathcal{D}_i}$ to $0$\,}\label{ConvergenceMart}
  For each $i\in\mathcal{H}$ and $j\in\mathcal{D}_i$, we now prove the weak convergence of $\{ M^{(n)}_{ij}\}_{n\geq 1}$  to $0$ by using the Burkholder-Davis-Gundy inequality.
  
  \begin{lemma}
  	For any $i\in\mathcal{H}$ and $j\in\mathcal{D}_i$, we have $M^{(n)}_{ ij}\overset{\rm u.c.p.}\longrightarrow 0$  as $n\to\infty$.
  \end{lemma}
  \proof For any $T\geq 0$, by the Burkholder-Davis-Gundy inequality and the hypothesis \ref{H1}
  we  have
  \beqnn
  \mathbf{E}\Big[\sup_{t\in[0,T]}\big|M^{(n)}_{ iI}(t)\big|^2\Big] \ar\leq \ar \frac{C}{n} \int_{\mathbb{U}} \big\|\phi_i(u)\big\|^2_{L^1}  \nu_I^{(n)}(du),
  \eeqnn
  which goes to $0$ as $n\to\infty$.
  Similarly,  for $j\in\mathcal{H}_i$, by Lemma~\ref{MomentEstimate} we also have
  \beqnn
  \mathbf{E}\Big[\sup_{t\in[0,T]}\big|M^{(n)}_{ij}(t)\big|^2 \Big] \ar\leq\ar C \mathbf{E}\Big[\int_0^T  \int_{\mathbb{U}} \frac{\big\|\phi_i(u)\big\|^2_{L^1}}{n^2} N_j^{(n)}(n\cdot ds,du) \Big] \cr
  \ar\leq\ar  C \int_{\mathbb{U}} \big\|\phi_i(u)\big\|^2_{L^1}  \nu_j^{(n)}(du).
  \eeqnn
  For any $K>0$, let $U_K:= \{u\in\mathbb{U}: \|\phi_i(u)\|_{L^1}\leq K  \}$ and $U_K^{\rm c}$ be its complement. We have
  \beqnn
  \int_{\mathbb{U}} \big\|\phi_i(u)\big\|^2_{L^1}  \nu_j^{(n)}(du) \ar\leq\ar  K \int_{U_K} \big\|\phi_i(u)\big\|_{L^1}  \nu_j^{(n)}(du) \cr
  \ar\ar + \frac{1}{K^{2\alpha-2}} \int_{U^{\rm c}_K} \big\|\phi_i(u)\big\|^{2\alpha}_{L^1}  \nu_j^{(n)}(du).
  \eeqnn
  The first term on the right side of this inequality can be bounded by $K\|\phi_{ij}^{(n)}\|_{L^1}$, which goes to $0$ as $n\to\infty$; see Condition~\ref{Con.PhiH}.
  By the hypothesis \ref{H1}, the second term can be uniformly bounded by $C/K^{2\alpha-2}$, which can be ignored for large $K$ and hence $ \mathbf{E}[\sup_{t\in[0,T]}|M^{(n)}_{ij}(t)|^2 ] \to0$ as $n\to\infty$.
  \qed

  \subsubsection{Uniform tightness and weak convergence of $\{W^{(n)}_\mathcal{H}\}_{n\geq 1}$\,}\label{ConvergenceW}
  By the mutual independency among $W^{(n)}_i$, $i\in\mathcal{H}$, it suffices to prove the uniform tightness and weak convergence of $L^{2}(\mathbb{R}_+)^{\#}$-martingales $\{  W^{(n)}_i\}_{n\geq 1}$ separately for each $i\in\mathcal{H}$.

  \begin{lemma}
  	For each $i\in\mathcal{H}$, we have $W^{(n)}_i \Rightarrow W_i$ as $n\to \infty$.
  \end{lemma}
  \proof By the continuity of $W_i$ and Corollary 3.33 in \cite[p.353]{JS03}, it suffices to prove $W^{(n)}_i(f) \overset{\rm d}\to W_i(f)$ in $\mathbf{D}([0,\infty),\mathbb{R})$ for any $f\in L^{2}(\mathbb{R}_+)$.
  Similarly as in the proof of Proposition~\ref{TightJii}, we can prove that $\{W^{(n)}_i(f)\}_{n\geq 1}$ is $C$-tight and 
  \beqnn
  \sup_{n\geq 1} \mathbf{E}\bigg[\sup_{t\in[0,T]}\Big|W^{(n)}_i(f,t)\Big|^2\bigg]<\infty,
  \quad T\geq 0.
  \eeqnn  
  We now start to characterize the cluster points.
  Without loss of generality, we may assume $W^{(n)}_i(f)$ converges to a limit process $X_f$ weakly and hence uniformly on compacts in probability.
  By the Skorokhod representation theorem, we may assume $W^{(n)}_i(f)\overset{\rm u.c.}\longrightarrow  X_f$ a.s. and hence in $L^2([0,T])$, which induces that  
  \beqnn
  \big|W^{(n)}_i(f,t)\big|^2- \frac{c_i^n}{c_i^2}\cdot \big\|f\big\|_{L^2}^2 \cdot t \overset{\rm u.c.}\longrightarrow \big|X_f(t)\big|^2- \big\|f\big\|_{L^2}^2 \cdot t, 
  \eeqnn
  a.s. and hence in $L^1([0,T])$.
  Thus both  $X_f$ and $\{ |X_f(t)|^2- \|f\|_{L^2}^2 \cdot t:t\geq 0 \}$ are martingales.
  In conclusion, we have  $X_f$ is a continuous martingale with quadratic variation $\langle X_f\rangle_t= \|f\|_{L^2}^2 \cdot t$ for $t\geq 0$.
  By Theorem III-7 in \cite{ElKarouiMeleard1990}, there exists a Gaussian white noise $W_i(ds,dz)$ on $(0,\infty)^2$ with intensity $dsdz$ such that
  \beqnn
  X_f(t)= \int_0^t \int_0^\infty f(z)W_i(ds,dz)=W_i(f,t),\quad t\geq 0.
  \eeqnn
  \qed
  
  \begin{lemma}
  	For each $i\in\mathcal{H}$, the sequence of $L^2(\mathbb{R}_+)^{\#}$-martingales $\{W_i^{(n)}\}_{n\geq 1}$ is uniformly tight.
  \end{lemma}
  \proof Let $\mathcal{S}$ be the collection of all $(\mathscr{F}_t)$-predictable $L^2(\mathbb{R}_+)$-valued processes.
  By the definition of uniform tightness, it suffices to prove that for any $T>0$,
  \beqnn
  \bigcup_{n=1}^\infty \bigg\{ \sup_{t\in[0,T]}| W_i^{(n)}(X,t)|: X\in \mathcal{S} \mbox{ with } \sup_{t\in[0,T]}\|X(t)\|_{L^2}\leq 1 \bigg\}
  \eeqnn
  is stochastically bounded. Actually, using Chebyshev's inequality and then the Burkholder-Davis-Gundy inequality together with the  hypothesis \ref{H1}, we have for any $\eta>0$,
  \beqnn
  \mathbf{P}\Big( \sup_{t\in[0,T]}| W_i^{(n)}(X,t)|\geq \eta\Big) \ar\leq \ar \eta^{-2}\mathbf{E}\Big[\sup_{t\in[0,T]}\Big| \int_0^t \int_0^\infty X(s,z)W_i^{(n)}(ds,dz) \Big|^2\Big]\cr
  \ar\leq\ar \frac{C}{\eta^2} \int_0^T \mathbf{E}[\|X(s)\|_{L^2}^2]ds \leq \frac{C}{\eta^2}\cdot T.
  \eeqnn
  This upper bound holds uniformly in $n\geq 1$ and $X\in \mathcal{H}$ with $\sup_{t\in[0,T]}\|X(t)\|_{L^2}\leq 1$. Thus the sequence $\{W_i^{(n)}\}_{n\geq 1}$ is uniformly tight. 
  \qed
  
  \subsection{Proof for Theorem~\ref{ConvergenceCumulativeEffect}}\label{Section6.6}
  
  By Condition~\ref{ConditionCumulativeEffect},  the rescaled process $S_{\mathtt{C},\mathcal{H}}^{(n)}$ can be well approximated by $\hat S_{\mathtt{C},\mathcal{H}}^{(n)}$ with
  \beqnn
  \hat S_{\mathtt{C},i}^{(n)}(t):= \int_0^t \int_\mathbb{U}  \frac{\zeta_i(\infty, u)}{n^2}N_i^{(n)}(n\cdot ds,du),\quad t\geq 0, i\in\mathcal{H}.
  \eeqnn
  The error process is denoted as $\varepsilon^{(n)}_{\mathtt{C},\mathcal{H}}:=\hat{S}_{\mathtt{C},\mathcal{H}}^{(n)}-S_{\mathtt{C},\mathcal{H}}^{(n)}$.
  By  Corollary~3.33 in \cite[p.353]{JS03}, Theorem~\ref{ConvergenceCumulativeEffect} follows directly from the following weak convergence results for the two sequences $\{\hat S_{\mathtt{C},\mathcal{H}}^{(n)}\}_{n\geq 1}$ and $\{\varepsilon^{(n)}_{\mathtt{C},\mathcal{H}}\}_{n\geq 1}$.
  The first one has been widely studied in \cite[Chapter~IX]{JS03} under Condition~\ref{MomentConCumulativeEffect}.
  
  \begin{lemma}
  	Theorem~\ref{ConvergenceCumulativeEffect} holds with $ S^{(n)}_{\mathtt{C},\mathcal{H}}$ replaced by $\hat S^{(n)}_{\mathtt{C},\mathcal{H}} $.
  \end{lemma}
  
  \begin{lemma}
  	We have $\varepsilon^{(n)}_{\mathtt{C},\mathcal{H}}  \overset{\rm d}\to0$ in $\mathbf{D}([0,\infty),\mathbb{R}^d)$ as $n\to\infty$.
  \end{lemma} 
  \proof For  $t\geq 0$ and $u\in\mathbb{U}$, let  $\zeta_i^{\mathrm{c}}(t,u):=\zeta_i(\infty,u)-\zeta_i(t,u)$.
  For any $\epsilon\in(0,1)$, we split $ \varepsilon^{(n)}_{\mathtt{C},i}(t)$ into the following two parts:
  \beqnn
  \bar\varepsilon^{(n)}_{\mathtt{C},i}(t,\epsilon)\ar:=\ar\int_0^{(t-\epsilon)^+} \int_{\mathbb{U}}  \frac{ \zeta_i^{\mathrm{c}}(n(t-s),u)}{n^2} N^{(n)}_i(n\cdot ds,du),\quad \label{eqnA.12}\\
  \underline\varepsilon^{(n)}_{\mathtt{C},i}(t,\epsilon)\ar:=\ar  \int_{(t-\epsilon)^+}^{t} \int_{\mathbb{U}}  \frac{ \zeta_i^{\mathrm{c}}(n(t-s),u)}{n^2}N^{(n)}_i(n\cdot ds,du).\label{eqnA.13}
  \eeqnn
  The monotonicity of $\zeta_i^{\mathrm{c}}(\cdot,u)$ induces that
  \beqnn
  |\bar\varepsilon^{(n)}_{\mathtt{C},i}(t,\epsilon)|
  \ar\leq\ar   \int_0^{t} \int_{\mathbb{U}}     \frac{ |\zeta_i^{\mathrm{c}}( n\epsilon,u) |}{n^2} N^{(n)}_i(n\cdot ds,du).
  \eeqnn
  From Lemma~\ref{MomentEstimate}  and Condition~\ref{ConditionCumulativeEffect}, we have for any $T>0$,
  \beqnn
  \mathbf{E}\Big[ \sup_{t\in[0,T]}  |\bar\varepsilon^{(n)}_{\mathtt{C},i}(t,\epsilon)| \Big]
  \ar\leq\ar  \int_0^T \mathbf{E}\big[ Z_i^{(n)}(s)\big]ds\cdot   \int_{\mathbb{U}}|\zeta_i^{\mathrm{c}}( n\epsilon,u) |\nu^{(n)}_i(du)
  \leq C\cdot \zeta_{ii}^{(n)\mathrm{c}}( n\epsilon),
  \eeqnn
  which vanishes  as $n\to\infty$.
  Moreover, by the monotonicity of $\zeta_i^{\mathrm{c}}(\cdot,u)$ again, 
  \beqnn
  \sup_{t\in[0,T]}  |\underline\varepsilon^{(n)}_{\mathtt{C},i}(t,\epsilon)|\ar\leq \ar    \sup_{t\in[0,T]} \int_{(t-\epsilon)^+}^{t} \int_{\mathbb{U}}   \frac{\zeta_i(\infty,u)}{n^2}  N^{(n)}_i(n\cdot ds,du)\cr
  \ar\leq\ar \sup_{0\leq j\leq [T/\epsilon]} \int_{j\epsilon}^{(j+2)\epsilon}  \int_{\mathbb{U}} \frac{\zeta_i(\infty,u)}{n^2}  N^{(n)}_i(n\cdot ds,du).
  \eeqnn
  By Chebyshev's inequality, we have for any $\eta>0$,
  \beqlb\label{eqnB.07}
  \mathbf{P}\Big( \sup_{t\in[0,T]}  |\underline\varepsilon^{(n)}_{\mathtt{C},i}(t,\epsilon)|\geq \eta \Big)
  \ar\leq \ar  \sum_{j=0}^{ [T/\epsilon]}\mathbf{P}\Big(\int_{j\epsilon}^{(j+2)\epsilon}  \int_{\mathbb{U}}    \frac{\zeta_i(\infty,u)}{n^2}  N^{(n)}_i(n\cdot ds,du)\geq \eta \Big) \cr
  \ar\leq\ar
  \frac{1}{\eta^\alpha} \sum_{j=0}^{ [T/\epsilon]}\mathbf{E}\Big[ \Big| \int_{j\epsilon}^{(j+2)\epsilon}  \int_{\mathbb{U}} \frac{\zeta_i(\infty,u)}{n^2}  N^{(n)}_i(n\cdot ds,du)\Big|^{\alpha} \Big].
  \eeqlb
  Using the Cauchy-Schwarz inequality inequality to the last expectation, it can be bounded by
  \beqnn
  \lefteqn{C \cdot \mathbf{E}\Big[ \Big| \int_{j\epsilon}^{(j+2)\epsilon} Z^{(n)}_i(s)ds \int_{\mathbb{U}} \zeta_i(\infty,u)  \nu^{(n)}_i(du)\Big|^{\alpha} \Big]}  \ar\ar\cr
  \ar\ar  + C \cdot \mathbf{E}\Big[ \Big| \int_{j\epsilon}^{(j+2)\epsilon}  \int_{\mathbb{U}} \frac{ \zeta_i(\infty,u)}{n^2}  \tilde{N}^{(n)}_i(dns,du)\Big|^{\alpha} \Big] ,
  \eeqnn
  which also, by (\ref{BDG}) and then H\"older's inequality and Lemma~\ref{MomentEstimate}, can be bounded by 
  \beqnn
   C\epsilon^{\alpha-1}\int_{j\epsilon}^{(j+2)\epsilon} \mathbf{E}\big[ \big| Z^{(n)}_i(s) \big|^\alpha\big]ds + \frac{C}{n^{2\alpha-2}}\int_{j\epsilon}^{(j+2)\epsilon} \mathbf{E}\big[ Z^{(n)}_i(s) \big]ds 
  \ar\leq\ar C\Big( \epsilon^{\alpha}  + \frac{\epsilon}{n^{2\alpha-2} } \Big).
  \eeqnn 
  Taking this back into (\ref{eqnB.07}), 
  \beqnn
  \mathbf{P}\Big( \sup_{t\in[0,T]}  |\underline\varepsilon^{(n)}_{\mathtt{C},i}(t,\epsilon)|\geq \eta \Big) \leq C\epsilon^{\alpha-1}  + C n^{2-2\alpha},
  \eeqnn
  which vanishes as $n\to\infty$ and then $\epsilon\to 0+$. We have finished the proof.
  \qed

  \subsection{Proofs for Theorem~\ref{ConvergenceInstantanuousEffect} and \ref{ConvergenceInstantanuousEffectAncestor}}
  
  Based on our asymptotic analysis before Theorem~\ref{ConvergenceInstantanuousEffect}, the rescaled process $S_{\mathtt{I},\mathcal{H}}^{(n)}$ can be well approximated by $\hat{S}_{\mathtt{I},\mathcal{H}}^{(n)}$ and the error process is  denoted as $\varepsilon_{\mathtt{I},\mathcal{H}}^{(n)}$.
  Firstly, we prove the weak convergence of $\{\hat{S}_{\mathtt{I},\mathcal{H}}^{(n)} \}_{n\geq 1}$ and $\{\varepsilon_{\mathtt{I},\mathcal{H}}^{(n)}\}_{n\geq 1}$ in the next two lemmas.

  \begin{lemma}\label{ConvergenceHatSI}
  	Theorem~\ref{ConvergenceInstantanuousEffect}  holds with $S_{\mathtt{I},\mathcal{H}}^{(n)}$ replaced by $\hat{S}_{\mathtt{I},\mathcal{H}}^{(n)}$. 
  \end{lemma}
  \proof By the Skorokhod representation theorem, we may assume that $Z_\mathcal{H}^{(n)}$ converges to $ Z_\mathcal{H}$ a.s. in $\mathbf{D}([0,\infty),\mathbb{R}_+^d)$ and hence uniformly on compacts.
  Thus it suffices to prove that for each $i\in\mathcal{H}$,
  \beqnn
  \int_0^{nt} \zeta_{ii}^{(n)}( s) Z_i^{(n)}(t-s/n) ds - b_{\mathtt{I},i}\cdot Z_i(t)
  \eeqnn
  goes to $0$ a.s. in $\mathbf{D}([\delta,1],\mathbb{R}_+)$.
  Subtracting the integral $\int_0^{nt} \zeta_{ii}^{(n)}(s) ds \cdot Z_i(t)$ and then adding it back, we can write the preceding quantity into
  \beqlb \label{eqn.165}
  \int_0^{nt} \zeta_{ii}^{(n)}(s) [Z_i^{(n)}(t-s/n)-Z_i(t)]  ds - I_{\zeta,ii}^{(n)}(nt)	 \cdot   Z_i(t) + \big(  \| \zeta_{ii}^{(n)}\|_{L^1} -b_{\mathtt{I},i} \big)  Z_i(t).
  \eeqlb
  Condition~\ref{ConInstantaneousEffects} implies that the last two terms go to $0$ as $n\to\infty$ uniformly on  $[\delta,1]$ and $[0,1]$  respectively.
  For any $\epsilon\in(0,1)$, the first term can be bounded by
  \beqnn
  \ar\ar\big\|\zeta_{ii}^{(n)}\big\|_{L^1} \cdot \Big( \sup_{t\in[0,1]} \big| Z_i^{(n)}(t)-Z_i(t) \big| + \sup_{t\in[0,1]} \sup_{s\in[0,\epsilon]} \big|\Delta_s Z_i(t) \big|\Big) \cr
  \ar\ar \qquad + I_{\zeta,ii}^{(n)}(n\epsilon) \cdot \sup_{t\in[0,1]}\Big( Z_i^{(n)}(t) +   Z_i(t)\Big) ,
  \eeqnn
  which goes to $0$ a.s. as $n\to\infty$, since $Z_i^{(n)} \overset{\rm a.s.}\to Z_i$ uniformly on $[0,1]$ and  $Z_i$ is uniformly continuous on $[0,2]$. 
  Putting these estimates together, we can get the first desired result.
  For the second one, like the preceding argument it suffices to prove that the second term in (\ref{eqn.165}) converges to $0$ uniformly on $[0,1]$ as $n\to\infty$. 
  Indeed, by the fact that $\sup_{n\geq 1}\|\zeta_{ii}^{(n)}\|_{L^1} +\sup_{t\in[0,1]}Z(t)<\infty$ a.s., we have for any $\epsilon\in(0,1)$,
  \beqnn
  \sup_{t\in[0,1]}  I_{\zeta,ii}^{(n)}(nt)	\cdot   Z_i(t)
  \ar\leq\ar \big\|\zeta_{ii}^{(n)}\big\|_{L^1} \cdot \sup_{t\in[0,\epsilon]}  Z_i(t) +  I_{\zeta,ii}^{(n)}(n\epsilon)	\cdot \sup_{t\in[0,1]}  Z_i(t)\cr
  \ar\leq\ar C  \sup_{t\in[0,\epsilon]} Z_i(t) + C  \int_{n\epsilon}^\infty \bar\zeta(s)ds ,
  \eeqnn
  which goes to $0$ a.s. as $n\to\infty$ and then $\epsilon\to 0+$ because of the continuity of $Z_i$ and the integrability of $\bar\zeta$. Here we have finished the proof.
  \qed

  \begin{lemma}\label{ErrorHatSI}
  	We have  $ \varepsilon_{\mathtt{I},\mathcal{H}}^{(n)}\overset{\rm d}\to 0$ in   $\mathbf{D}([0,\infty),\mathbb{R}^d)$ as $n\to\infty$.
  \end{lemma}
  \proof  By the Burkholder-Davis-Gundy inequality, the inequality $(x+y)^{\alpha/2} \leq |x|^{\alpha/2}+|y|^{\alpha/2}$ and Lemma~\ref{MomentEstimate}, we have for any $t\geq 0$,
  \beqnn
  \mathbf{E}\big[\big|\varepsilon_{\mathtt{I},i}^{(n)}(t) \big|^\alpha\big]
  \ar\leq\ar C \cdot \mathbf{E}\Big[ \Big( \int_0^{nt}\int_{\mathbb{U}} \frac{|\zeta_i(nt-s,u)|^2}{n^2}N^{(n)}_i(ds,du) \Big)^{\alpha/2} \Big]\cr
  \ar\leq\ar C \cdot \mathbf{E}\Big[   \int_0^{nt}\int_{\mathbb{U}} \frac{|\zeta_i(nt-s,u)|^\alpha}{n^\alpha}N^{(n)}_i(ds,du)  \Big] \cr
  \ar\leq\ar \frac{C}{n^{\alpha-1}}\cdot \int_\mathbb{U} \big\|\zeta_i(u)\big\|_{L^\alpha}^\alpha\nu_i^{(n)}(du).
  \eeqnn
  Notice that 
  \beqnn
  \big\|\zeta_i(u)\big\|_{L^\alpha}^\alpha
  \leq  \big\|\zeta_i(u)\big\|_{\rm TV}^{\alpha-1} \cdot \big\|\zeta_i(u)\big\|_{L^1} 
  \leq  C\big(\big\|\zeta_i(u)\big\|_{\rm TV}^{\alpha} + \big\|\zeta_i(u)\big\|_{L^1}^\alpha \big).
  \eeqnn
  By Condition~\ref{ConInstantaneousEffects02}, we have $\mathbf{E}\big[\big|\varepsilon_{\mathtt{I},i}^{(n)}(t) \big|^\alpha\big]\to0$ as $n\to\infty$ and hence $\varepsilon_{\mathtt{I},i}^{(n)} \overset{\rm f.d.d.} \longrightarrow 0$.
  We now start to prove the tightness of $\big\{ \varepsilon^{(n)}_{\mathtt{I},i} \big\}_{n\geq 1}$ on $[0,1]$ and the general case can be proved similarly.
  Since $\big\|\zeta_i(u)\big\|_{\rm TV}<\infty$ for any $u\in\mathbb{U}$, similarly as in  the proof of Proposition~\ref{WellApproximationIii}, it suffices to prove the case in which $\zeta_i(t,u)$ decreases in $t$.
  \smallskip
  
  {\it Step 1.} We first show that $\varepsilon^{(n)}_{\mathtt{I},i}$ can be well approximated by its linear interpolation $ \varepsilon^{(n)}_{\mathtt{I},i,\theta}$ defined as (\ref{LinearInterSmallError}), i.e.	$\varepsilon^{(n)}_{\mathtt{I},i}- \varepsilon^{(n)}_{\mathtt{I},i,\theta}\overset{\rm u.c.p.}\longrightarrow 0$.
  Like  (\ref{UpperBoundApproximation01})-(\ref{UpperBoundApproximation02}), we also have
  \beqlb \label{eqn.158}
  \sup_{t\in[0,1]}  \big|\varepsilon^{(n)}_{\mathtt{I},i}(t) - \varepsilon^{(n)}_{\mathtt{I},i,\theta}(t) \big| 
  \ar\leq\ar 3 \sup_{k=0,\cdots,[n^\theta]; h\leq n^{-\theta}} \big[A_{\mathtt{I},1}^{(n)}(kn^{-\theta},h) +
  A_{\mathtt{I},2}^{(n)}(kn^{-\theta},h)\big],
  \eeqlb
  where
  \beqnn
  A_{\mathtt{I},1}^{(n)}(t,h)\ar:=\ar n \int_0^{t+h} Z_i^{(n)}(s)ds \int_\mathbb{U} \big|\Delta_{nh}\zeta_i(n(t-s),u)\big| \nu_i^{(n)}(du)  ,\cr
  A_{\mathtt{I},2}^{(n)}(t,h)\ar:=\ar \int_0^t \int_\mathbb{U} \frac{|\Delta_{nh}\zeta_i(n(t-s),u)| }{n}N_i^{(n)}(n\cdot ds,du) .
  \eeqnn
  For any $t,h\in[0,1]$, we first have
  \beqnn
  A_{\mathtt{I},1}^{(n)}(t,h)\ar \leq\ar \sup_{r\in[0,2]}Z_i^{(n)}(r)\cdot  \int_\mathbb{U} \nu_i^{(n)}(du)\cdot  \int_0^\infty \big|\Delta_{nh}\zeta_i(s,u) \big|  ds. 
  \eeqnn
  Similarly as in (\ref{DifferencePhi}), we have $\int_0^\infty |\Delta_{nh}\zeta_i(s,u) |  ds \leq 2\|\zeta_i(u)\|_{\rm TV}\cdot nh$. By Condition~\ref{ConInstantaneousEffects02}, we have uniformly in $h\in[0,1]$,
  \beqnn
  \sup_{t\in[0,1]} A_{\mathtt{I},1}^{(n)}(t,h) \ar\leq\ar  C \sup_{r\in[0,2]}Z_i^{(n)}(r)\cdot nh.
  \eeqnn
  Since $ Z_i$ is continuous, we have $\sup_{n\geq 1}\sup_{r\in[0,2]}Z_i^{(n)}(r)<\infty$ a.s. and hence as $n\to\infty$,
  \beqnn
  \sup_{k=0,\cdots,[n^\theta]; h\leq n^{-\theta}}
  A_{\mathtt{I},1}^{(n)}(kn^{-\theta},h)  \leq \sup_{r\in[0,2]}Z_i^{(n)}(r)\cdot   n^{1-\theta}\overset{\rm a.s.}\to 0.
  \eeqnn 
  Moreover, like Step~3 in the proof of Proposition~\ref{WellApproximationIii},  we can prove that for any $\eta>0$,
  \beqnn
  \lim_{n\to\infty}\mathbf{P}\Big(  \sup_{k=0,\cdots,[n^\theta]; h\leq n^{-\theta}}   A_{\mathtt{I},2}^{(n)}(kn^{-\theta},h)\geq \eta   \Big) = 0.
  \eeqnn
  Taking these two estimates back into (\ref{eqn.158}), we  have $\sup_{t\in[0,1]}|\varepsilon^{(n)}_{\mathtt{I},i}(t)-\varepsilon^{(n)}_{\mathtt{I},i,\theta}(t)|\overset{\rm p}\to 0$ as $n\to\infty$.
  \smallskip
  
  {\it Step 2.} We now prove the tightness of the sequence $\{  \varepsilon^{(n)}_{\mathtt{I},i,\theta} \}_{n\geq 1}$.
  By (\ref{BDG}), there exits a  constant $C>0$ such  that for any $n\geq 1$ and $t,h\in[0,1]$, 
  \beqnn
  \mathbf{E}\big[\big|\Delta_h\varepsilon^{(n)}_{\mathtt{I},i}(t) \big|^{2\alpha}\big]
  \ar\leq\ar C \Big(  \int_\mathbb{U} \nu^{(n)}_i(du) \int_0^{t+h} |\Delta_{nh}\zeta_i(n(t-s),u)|^2 ds  \Big)^\alpha\cr
  \ar\ar +  \frac{C}{n^{2\alpha-2}} \int_\mathbb{U} \nu^{(n)}_i(du) \int_0^{t+h} |\Delta_{nh}\zeta_i(n(t-s),u)|^{2\alpha} ds. 
  \eeqnn
  Notice that $\int_0^{t+h} |\Delta_{nh}\zeta_i(n(t-s),u)|^{2\alpha} ds  \leq 2\|\zeta_i(u)\|_{\rm TV}^{2\alpha}\cdot h$ uniformly in $t,h\in[0,1]$.
  By Condition~\ref{ConInstantaneousEffects02}, there exists a constant $C>0$ such that 
  \beqnn 
  \mathbf{E}\big[\big|\Delta_h\varepsilon^{(n)}_{\mathtt{I},i}(t) \big|^{2\alpha}\big]\leq C\cdot\big( n^{2-2\alpha} h  +h^\alpha\big),
  \eeqnn
  for any $t,h \in[0,1]$ and $n\geq 1$.
  Like the proof of Proposition~\ref{TightnessIii},  we can prove the tightness of the sequence $\big\{  \varepsilon^{(n)}_{\mathtt{I},i,\theta} \big\}_{n\geq 1}$ in the same way.  
  Consequently, the sequence $\big\{ \varepsilon^{(n)}_{\mathtt{I},i} \big\}_{n\geq 1}$ is tight in $\mathbf{D}([0,\infty),\mathbb{R})$ and the whole proof is end.
  \qed
  
  \textsc{Proofs for Theorem~\ref{ConvergenceInstantanuousEffect} and \ref{ConvergenceInstantanuousEffectAncestor}. }
  By  Corollary~3.33 in \cite[p.353]{JS03}, we can get Theorem~\ref{ConvergenceInstantanuousEffect} directly from Lemma~\ref{ConvergenceHatSI} and \ref{ErrorHatSI}.
  For Theorem~\ref{ConvergenceInstantanuousEffectAncestor},
  by Condition~\ref{Con.InitialState01} and Theorem~\ref{ConvergenceInstantanuousEffect} it suffices to prove that 
  \beqnn
  \hat\psi_{\mathtt{I},i}^{(n)}(nt) +  \hat{S}_{\mathtt{I},i}^{(n)}(t)\overset{\rm d}\to b_{\mathtt{I},i} \cdot Z_i(t),
  \eeqnn
  in $\mathbf{D}([0,1],\mathbb{R}_+)$ as $n\to\infty$.
  Like the proof of Lemma~\ref{ConvergenceHatSI}, it suffices to prove that
  \beqnn
  \lefteqn{Z_i^{(n)}(0) I_{\zeta,ii}^{(n)}(nt)+ \int_0^{nt} \zeta_{ii}^{(n)}(s) Z_i^{(n)}(t-s/n) ds - b_{\mathtt{I},i}\cdot Z_i(t)}\ar\ar\cr
  \ar=\ar  I_{\zeta,ii}^{(n)}(nt)\cdot   \big(Z_i^{(n)}(0) -Z_i(t) \big)+  \big(  \big\| \zeta_{ii}^{(n)}\big\|_{L^1} -b_{\mathtt{I},i} \big)  Z_i(t) \cr
  \ar\ar  + \int_0^{nt} \zeta_{ii}^{(n)}(s) \cdot \big[Z_i^{(n)}(t-s/n)-Z_i(t)\big]  ds
  \eeqnn
  goes to $0$ a.s. in $\mathbf{D}([0,1],\mathbb{R})$. From the proof of Lemma~\ref{ConvergenceHatSI}, the last two terms on the right side of this equality go to $0$ uniformly on compacts. 
  For any $\epsilon\in(0,1)$, the first term can be bounded by
  \beqnn
  I_{\zeta,ii}^{(n)}(n\epsilon)\cdot   \sup_{t\in[\epsilon,1]} \big|Z_i^{(n)}(0) -Z_i(t)\big| +\sup_{t\in[0,\epsilon]}I_{\zeta,ii}^{(n)}(nt)\cdot   \big|Z_i^{(n)}(0) -Z_i(t)\big|.
  \eeqnn
  Since $\sup_{t\in[\epsilon,1]} |Z_i^{(n)}(0) -Z_i(t)|<\infty$ a.s. and $I_{\zeta,ii}^{(n)}(n\epsilon)\to0$, the first term in this sum goes to $0$  as $n\to\infty$.
  Moreover, the second term can be bounded by
  \beqnn
  \lefteqn{\sup_{t\in[0,\epsilon]}I_{\zeta,ii}^{(n)}(nt)\cdot  \big( \big|Z_i^{(n)}(0) -Z_i(0)\big| +   \big|Z_i(0) -Z_i(t)\big|\big)}\ar\ar\cr
  \ar\leq\ar C\cdot \big|Z_i^{(n)}(0)-Z_i(0)\big|+C\sup_{t\in[0,\epsilon]}\big|Z_i(0) -Z_i(t)\big|,
  \eeqnn
  which goes to $0$ a.s. as $n\to\infty$ and then $\epsilon\to0+$. The proof is end.
  \qed

  \subsection{Proof for Theorem~\ref{ConvergenceBirthRateProcesses}}
  
  Before proving Theorem~\ref{ConvergenceBirthRateProcesses} by using Theorem~\ref{MainThm01}, it remains to identify the total budding rate function of ancestors satisfies Condition~\ref{MomentConditionInitialState}. Indeed, by Condition~\ref{ConditionInitialStateCMJ}, the mean budding rate function of each type-$i$ ancestor, denoted as $\breve{\mathtt{B}}^{(n)}_i$, is
  \beqlb\label{eqn114}
  \breve{\mathtt{B}}^{(n)}_i(t)=\int_\mathbb{B}  \mathcal{P} _{\mathtt{B},i}^{(n)}(d\mathtt{B})
  \int_{\mathbb{R}_+^2} \mathtt{B}(t+s,y+s)  \breve{\mathcal{P}}_{\mathtt{A}\mathtt{R},i}^{(n)}(ds,dy),\quad t\geq 0.
  \eeqlb
  By Fubini's theorem and the fact that $\mathtt{B}(t,y)=0$ for $t\geq y$,
  \beqlb\label{eqn115}
  \int_{\mathbb{R}_+^2} \mathtt{B}(t+s,y+s)  \breve{\mathcal{P}}_{\mathtt{AR},i}^{(n)}(ds,dy)
  \ar=\ar \frac{1}{\mathrm{m}_{\mathtt{L},i}^{(n)}} \int_0^\infty ds \int_s^\infty \mathtt{B}(t+s,y) \mathcal{P}^{(n)}_{\mathtt{L},i}(dy)\cr
  \ar=\ar \frac{1}{\mathrm{m}_{\mathtt{L},i}^{(n)}} \int_t^\infty ds \int_0^\infty \mathtt{B}(s,y)  \mathcal{P}^{(n)}_{\mathtt{L},i}(dy)
  \eeqlb
 and hence
 \beqlb\label{eqn116}
 \breve{\mathtt{B}}^{(n)}_i(t)
 = \frac{1}{\mathrm{m}_{\mathtt{L},i}^{(n)}} \int_t^\infty ds \int_{\mathbb{B}}  \mathcal{P}^{(n)}_{\mathtt{B},i}(d\mathtt{B}) \int_0^\infty \mathtt{B}(s,y)  \mathcal{P}^{(n)}_{\mathtt{L},i}(dy)  
 \ar=\ar \frac{1}{\mathrm{m}_{\mathtt{L},i}^{(n)}} \int_t^\infty \mathtt{B}^{(n)}_i(s)ds\cr 
 \ar=\ar \frac{I^{(n)}_{\phi,ii}(t)}{\mathrm{m}_{\mathtt{L},i}^{(n)}\mathrm{m}_{ii}^{(n)}}.
 \eeqlb
 By the law of large numbers, it is natural to believe that $\mu^{(n)}_{\mathcal{H}}/n$ can be well approximated by 
 \beqnn
 \hat\mu^{(n)}_{\mathcal{H}}:=\big(Z_i^{(n)}(0)\cdot I^{(n)}_{\phi,ii}\big)_{i\in\mathcal{H}}
 \quad\mbox{with}\quad
 Z_i^{(n)}(0):=  \frac{\Xi_i^{(n)}(0)/n}{\mathrm{m}_{\mathtt{L},i}^{(n)}\mathrm{m}_{ii}^{(n)}},
 \quad i\in\mathcal{H}.
 \eeqnn
  
  \begin{lemma}\label{ConvergenceInitialStateCMJ}
  	We have $ \big\|\mu^{(n)}_{\mathcal{H}}/n- \hat\mu^{(n)}_{\mathcal{H}} \big\|_{L^{1,\infty}}\overset{\rm d}\to0$ as $n\to\infty$.
  \end{lemma}
  \proof
  For each $i\in\mathcal{H}$, let $\tilde\varepsilon^{(n)}_{\mu_i}:= \mu^{(n)}_i/n- \hat\mu^{(n)}_i$.
  For any $\eta>0$ and $K>0$, we have
 \beqnn
 \mathbf{P} \big( \big\|\tilde\varepsilon^{(n)}_{\mu_i}\big\|_{L^{1,\infty}}>\eta\big)
 \leq \mathbf{P}\big( Z_i^{(n)}(0) >K  \big)
 +\mathbf{P}\big(\big\|\tilde\varepsilon^{(n)}_{\mu_i}\big\|_{L^{1,\infty}}>\eta,Z_i^{(n)}(0) \leq K \big).
 \eeqnn
 By Condition~\ref{ConvergenceParameterCMJ} and the assumption that $ \Xi_i^{(n)}(0)/n \overset{\rm d}\to \Xi_i^*(0)$ as $n\to\infty$, 
 \beqnn
 \lim_{K\to \infty}\sup_{n\geq 1}\mathbf{P}\big(\big|Z_i^{(n)}(0)\big| >K \big) =0 .
 \eeqnn
 Thus it suffices to prove this lemma with $\{\Xi_i^{(n)}(0)/n\}_{n\geq 1}$ being deterministic and uniformly bounded.
 The following proof follows closely that of Theorem~4.1 in \cite{Acosta1981}.
 In detail, let $\{Y_{i,k}^{(n)}:k=1,2,\cdots, \Xi_i^{(n)}(0)\}$ be a sequence of i.i.d. function-valued random variables with 
  \beqnn
  Y_{i,k}^{(n)}(t):=\mathtt{B}^{(n)}_{i,k} \big(t+\mathtt{A}^{(n)}_{i,k},\mathtt{R}_{i,k}^{(n)}+\mathtt{A}^{(n)}_{i,k} \big)  - \breve{\mathtt{B}}^{(n)}_i(t) ,
  \quad t\geq 0. 
  \eeqnn 
  From (\ref{InitialStateCMJ}), we have 
  \beqnn
  \tilde\varepsilon^{(n)}_{\mu_i}(t) =  \frac{\mu^{(n)}_{i}(t)}{n}- \frac{\Xi_i^{(n)}(0)}{n} \cdot \breve{\mathtt{B}}^{(n)}_i(t) =\frac{1}{n}\sum_{k=1}^{\Xi_i^{(n)}(0)}  Y_{i,k}^{(n)}(t),\quad t\geq 0, i\in\mathcal{H}
  \eeqnn 
  From (\ref{eqn114})-(\ref{eqn116}), we have  $\mathbf{E}\big[Y_{i,k}^{(n)}(t)\big]=0$ for any $t\geq 0$. 
  By (\ref{eqn116}), Fubini's theorem and Condition~\ref{CMJMomentCondition}, there exists a constant $C>0$ such that for any $n\geq 1$ and $i\in\mathcal{H}$,
  \beqnn
  \big\|\breve{\mathtt{B}}^{(n)}_i\big\|_{L^1}=    \int_0^\infty  \frac{s \cdot \mathtt{B}^{(n)}_i(s)}{\mathrm{m}_{\mathtt{L},i}^{(n)}}ds \leq C
  \quad \mbox{and}\quad 
  \big\|\breve{\mathtt{B}}^{(n)}_i\big\|_{L^\infty} = \breve{\mathtt{B}}^{(n)}_i(0) \leq C.
  \eeqnn
  Notice that 
  \beqnn
  \int_0^\infty\mathtt{B}^{(n)}_{i,k} (t +\mathtt{A}^{(n)}_{i,k},  \mathtt{R}_{i,k}^{(n)}+\mathtt{A}^{(n)}_{i,k})dt\leq \big\|\mathtt{B}^{(n)}_{i,k} (  \mathtt{R}_{i,k}^{(n)}+\mathtt{A}^{(n)}_{i,k})\big\|_{L^1}
  \eeqnn
  and 
  \beqnn
  \sup_{t\geq 0}\big|\mathtt{B}^{(n)}_{i,k} (t+\mathtt{A}^{(n)}_{i,k},\mathtt{R}_{i,k}^{(n)}+\mathtt{A}^{(n)}_{i,k})\big|
  \leq \big\|\mathtt{B}^{(n)}_{i,k} (\mathtt{R}_{i,k}^{(n)}+\mathtt{A}^{(n)}_{i,k})\big\|_{\rm TV}. 
  \eeqnn 
  By (\ref{eqn115}), Fubini's theorem, the inequality $2|xy| \leq x^2+y^2 $ and Condition~\ref{CMJMomentCondition},  
  \beqnn
  \mathbf{E}\big[\big\|\mathtt{B}^{(n)}_{i,k} (\mathtt{R}_{i,k}^{(n)}+\mathtt{A}^{(n)}_{i,k})\big\|_{L^1}^\alpha\big]
  \ar=\ar \int_{\mathbb{B}} \mathcal{P}^{(n)}_{\mathtt{B},i}(d\mathtt{B}) \int_{\mathbb{R}^2_+} \big\|\mathtt{B}  (s+y)\big\|_{L^1}^\alpha\breve{\mathcal{P}}_{\mathtt{AR},i}^{(n)}(ds,dy) \cr
  \ar=\ar  \int_{\mathbb{B}} \mathcal{P}^{(n)}_{\mathtt{B},i}(d\mathtt{B}) \int_0^\infty \frac{ds}{\mathrm{m}_{\mathtt{L},i}^{(n)}}\int_s^\infty \big\|\mathtt{B}  (y)\big\|_{L^1}^\alpha\breve{\mathcal{P}}_{\mathtt{L},i}^{(n)}( dy)\cr
  \ar=\ar  \int_{\mathbb{B}} \mathcal{P}^{(n)}_{\mathtt{B},i}(d\mathtt{B}) \int_0^\infty \frac{y\cdot \big\|\mathtt{B}  (y)\big\|_{L^1}^\alpha}{\mathrm{m}_{\mathtt{L},i}^{(n)}} \breve{\mathcal{P}}_{\mathtt{L},i}^{(n)}( dy)\cr
  \ar\leq\ar    \int_0^\infty \frac{y^2}{\mathrm{m}_{\mathtt{L},i}^{(n)}} \breve{\mathcal{P}}_{\mathtt{L},i}^{(n)}( dy) 
   +  \int_{\mathbb{B}} \mathcal{P}^{(n)}_{\mathtt{B},i}(d\mathtt{B}) \int_0^\infty \frac{ \big\|\mathtt{B}  (y)\big\|_{L^1}^{2\alpha}}{\mathrm{m}_{\mathtt{L},i}^{(n)}} \breve{\mathcal{P}}_{\mathtt{L},i}^{(n)}( dy),
  \eeqnn
  which is bounded uniformly in $n\geq 1$ and $i\in\mathcal{H}$. 
  Similarly, we also have 
  \beqnn
  \sup_{n\geq 1}\mathbf{E}\big[\big\|\mathtt{B}^{(n)}_{i,k} (\mathtt{R}_{i,k}^{(n)}+\mathtt{A}^{(n)}_{i,k})\big\|_{\rm TV}^\alpha\big]  <\infty,\quad i\in\mathcal{H}.
  \eeqnn 
  Putting all these estimates together and then using  Minkowski's inequality, we have 
  \beqnn
  \sup_{n\geq 1}\mathbf{E}\Big[\big\|Y_{i,k}^{(n)}\big\|_{L^{1,\infty}}^\alpha\Big]<\infty,\quad i\in\mathcal{H}.
  \eeqnn
  For $K>0$, by Chebyshev's inequality and Minkowski's inequality, 
  \beqnn
  \sup_{n\geq 1}\mathbf{P}\big(\big\| \tilde\varepsilon^{(n)}_{\mu_i}\big\|_{L^{1,\infty}}\geq K\big) 
  \ar\leq\ar \frac{1}{K}\cdot\sup_{n\geq 1} \mathbf{E}\big[ \big\| \tilde\varepsilon^{(n)}_{\mu_i}\big\|_{L^{1,\infty}} \big] \cr
  \ar\leq\ar \frac{1 }{K}\sup_{n\geq 1} \frac{\Xi_i^{(n)}(0)}{n} \cdot \sup_{n\geq 1}\mathbf{E}\big[\big\|Y_{i,k}^{(n)}\big\|_{L^{1,\infty}}\big]\leq \frac{C}{K},
  \eeqnn 
  which vanishes as $K\to\infty$. 
  Thus the sequence $\{\varepsilon^{(n)}_{\mu_i}\}_{n\geq 1}$ is tight and hence flatly concentrated; see Definition~2.1~in \cite{Acosta1970}.
  Notice that $L^1(\mathbb{R}_+)\cup L^\infty(\mathbb{R}_+)$ is the dual space of $L^{1,\infty}(\mathbb{R}_+)$.
  For any $f\in L^1(\mathbb{R}_+)\cup L^\infty(\mathbb{R}_+)$, we see that $\big\{f(Y_{i,k}^{(n)}):k=1,\cdots,\Xi_i^{(n)}(0)\big\}_{n\geq 1}$ is an array of row-wise independent random variables with 
  \beqnn
  \mathbf{E}\big[\big|f(Y_{i,k}^{(n)})\big|^\alpha\big] \leq C \mathbf{E}\big[\big\|Y_{i,k}^{(n)})\big\|_{L^{1,\infty}}^\alpha\big] \leq C,
  \eeqnn
  uniformly in $n\geq 1$ and hence this array is uniformly integrable.
  By using the main theorem in \cite{Gut1992}, we have as $n\to\infty$,
  \beqnn
  f( \tilde\varepsilon^{(n)}_{\mu_i})= \frac{1}{n}\sum_{k=1}^{\Xi_i^{(n)}(0)} f(Y_{i,k}^{(n)})\overset{\rm p} \to 0
  \eeqnn
  By Theorem~2.4 in \cite{Acosta1970}, it follows that $\big\| \tilde\varepsilon^{(n)}_{\mu_i} \big\|_{L^{1,\infty}} \overset{\rm d}\to 0$ as $n\to\infty$.
  \qed

  \subsection{Proof for Theorem~\ref{ConvergenceCMJInstantaneousEffect}}
  It is obvious that this theorem follows directly from Theorem~\ref{ConvergenceInstantanuousEffectAncestor} and the following two auxiliary results.
  
  {\bf (1)}  {\it The impact of immigrants and offspring born by parents of different type can be ignored}, i.e.,
  for each $i\in\mathcal{H}$ and $j\in\mathcal{D}_i$, the rescaled process $\varepsilon^{(n)}_{N,ij}$ defined as below converges weakly to $0$ in $\mathbf{D}([0,\infty),\mathbb{R})$:
  \beqnn
  \varepsilon^{(n)}_{N,ij}(t):=\frac{1}{n}\int_0^{nt} \int_\mathbb{U}  \zeta_i(nt-s,\boldsymbol{u}) N^{(n)}_j(ds,d\boldsymbol{u}), \quad t\geq 0.
  \eeqnn
  Here we prove it with $j\in\mathcal{H}_i$. For the case of $j=I$, it can be proved similarly. 
  By Lemma~\ref{MomentEstimate} and Condition~\ref{CMJInstantaneousConvergence}, there exist two constants $C,\vartheta>0$ such that $ \sup_{n\geq 1} \mathbf{E}[|\mathbf{B}_\mathcal{H}^{(n)}(t)|^{2\alpha}]\leq Ce^{\vartheta t}$ for any $t\geq 0$ and
  \beqnn
  \mathbf{E}[|\varepsilon^{(n)}_{N,ij}(t)|]
  \ar\leq\ar C\int_\mathbb{U}  \|\zeta_i(\boldsymbol{u})\|_{L^1} \nu^{(n)}_j(d\boldsymbol{u}) \cr
  \ar=\ar C \cdot \mathrm{m}^{(n)}_{ij}\int_0^\infty \mathcal{P}^{(n)}_{\mathtt{L},i}(d\mathtt{y}) \int_\mathbb{T}   \|\mathtt{T}(\mathtt{y})\|_{L^1} \mathcal{P}_{\mathtt{T},i}^{(n)}(d\mathtt{T})
  \leq  C \cdot \mathrm{m}^{(n)}_{ij} ,
  \eeqnn
  which goes to $0$; see Condition~\ref{ConvergenceParameterCMJ}.
  Hence $\varepsilon^{(n)}_{N,ij}\overset{\rm f.d.d.}\longrightarrow 0$.
  We now prove the tightness of $\{\varepsilon^{(n)}_{N,ij}\}_{n\geq 1}$.
  Because of $\|\mathtt{T}(y)\|_{\rm TV}<\infty$ for any $\mathtt{T}\in\mathbb{T}$ and $y\geq 0$, we have  $\|\zeta_i(\boldsymbol{u})\|_{\rm TV}<\infty$ for any $\boldsymbol{u}\in\mathbb{U}$.
  Similarly as in the proof of Lemma \ref{ErrorHatSI}, it suffices to consider the case with $\zeta_i(t,\boldsymbol{u})$ being non-increasing in $t$.
  Let $ \varepsilon^{(n)}_{N,ij,\theta}$ be the linear interpolation of $\varepsilon^{(n)}_{N,ij}$  defined as (\ref{LinearInterSmallError}).
  Like the argument in (\ref{UpperBoundApproximation01})-(\ref{UpperBoundApproximation02}), we have
  \beqnn
  \lefteqn{\sup_{t\in[0,1]}  \Big|\varepsilon^{(n)}_{N,ij}(t) - \varepsilon^{(n)}_{N,ij,\theta} (t) \Big| }\ar\ar\cr 
  \ar\leq\ar 3 \sup_{k=0,\cdots,[n^\theta]; h\leq n^{-\theta}}
  \int_0^{kn^{-\theta}} \int_\mathbb{U} \frac{|\Delta_{nh}\zeta_i(n(kn^{-\theta}-s),\boldsymbol{u})| }{n}N_j^{(n)}(n\cdot ds,d\boldsymbol{u}).
  \eeqnn
  Proceeding as in Step~3 of the proof for Lemma~\ref{WellApproximationIii}, we have the foregoing supremum goes to $0$ in probability as $n\to\infty$.
  Like the proof of Proposition~\ref{TightnessIii}, we now turn to prove the $C$-tightness of the sequence $\{ \varepsilon^{(n)}_{N,ij,\theta}\}_{n\geq 1}$. 
  For any $t,h\in[0,1]$,
  using the Cauchy-Schwarz inequality and then (\ref{BDG}), 
  \beqlb\label{eqn.300} 
  \mathbf{E}\big[\big|\Delta_h\varepsilon^{(n)}_{N,ij}(t)\big|^{2\alpha}\big]\ar \leq \ar  C  \Big(n\int_0^{t+h} ds \int_\mathbb{U} \big|\Delta_{nh}\zeta_i(n(t-s),\boldsymbol{u}) \big| \nu_j^{(n)}(d\boldsymbol{u}) \Big)^{2\alpha} \cr
  \ar\ar + C  \Big(\int_0^{t+h} ds \int_\mathbb{U} \big|\Delta_{nh}\zeta_i(n(t-s),\boldsymbol{u}) \big|^2 \nu_j^{(n)}(d\boldsymbol{u}) \Big)^{\alpha} \cr
  \ar\ar + C \int_0^{t+h} ds \int_\mathbb{U} \frac{|\Delta_{nh}\zeta_i(n(t-s),\boldsymbol{u}) |^{2\alpha}}{n^{2\alpha-2}} \nu_j^{(n)}(d\boldsymbol{u}). 
  \eeqlb
  The monotonicity of $\zeta_i(\cdot, \boldsymbol{u})$ induces that
  \beqnn
  \int_0^{t+h} \big|\Delta_{nh}\zeta_i\big(n(t-s),\boldsymbol{u}\big) \big| ds = \int_t^{t+h}  \zeta_i(ns,\boldsymbol{u}) ds 
  \leq \big\| \zeta_i(\boldsymbol{u})\big\|_{\rm TV}\cdot h. 
  \eeqnn
  By (\ref{eqn.301}), the first term on the ride side of (\ref{eqn.300}) can bounded by
  \beqnn
  \lefteqn{C  \Big(n \cdot  h \int_\mathbb{U}  \big\| \zeta_i(\boldsymbol{u})\big\|_{\rm TV} \nu_j^{(n)}(d\boldsymbol{u}) \Big)^{2\alpha} }\ar\ar\cr
  \ar\leq\ar C h^{2\alpha} \Big( n\cdot\mathrm{m}_{ij}^{(n)}\int_0^\infty \mathcal{P}^{(n)}_{\mathtt{L},i}(d\mathtt{y}) \int_\mathbb{T} \big\|\mathtt{T}(\mathtt{y})\big\|_{\rm TV} \mathcal{P}_{\mathtt{T},i}^{(n)}(d\mathtt{T})  \Big)^{2\alpha},
  \eeqnn
  which can be  uniformly bounded by $Ch^{2\alpha}$; see Condition~\ref{ConvergenceParameterCMJ} and \ref{CMJInstantaneousConvergence}.
  Similarly, the second term on the ride side of (\ref{eqn.300}) also can be bounded by
  \beqnn
  C  \Big( \int_\mathbb{U} \big\|\zeta_i(\boldsymbol{u}) \big\|_{\rm TV}^2 \nu_j^{(n)}(d\boldsymbol{u}) \Big)^{\alpha} \cdot h^\alpha \leq Ch^\alpha
  \eeqnn
  and the third term can be bounded by $Ch/n^{2\alpha-2}$.
  Putting all estimates above together, we have   $\mathbf{E}[|\Delta_h\varepsilon^{(n)}_{N,ij}(t)|^{2\alpha}] \leq C(h^\alpha+ h/n^{2\alpha-2})$.
  Like the proof of Proposition~\ref{TightnessIii}, we have  $\{ \varepsilon^{(n)}_{N,ij,\theta}\}_{n\geq 1}$ is $C$-tight and so is the sequence $\{\varepsilon^{(n)}_{N,ij}\}_{n\geq 1}$.
  
  {\bf (2)} {\it The impact of ancestors, $\psi_\mathcal{H}^{(n)}$, satisfies Condition~\ref{Con.InitialState01}.}
  Similarly as in (\ref{eqn114})-(\ref{eqn116}), the mean instantaneous characteristic of a type-$i$ ancestor at time $t$, denoted as $\breve{\mathtt{T}}^{(n)}_i(t)$, equals to
  \beqnn
  \breve{\mathtt{T}}^{(n)}_i(t)
  \ar:=\ar \int_\mathbb{T} \mathcal{P}_{\mathtt{T},i}^{(n)}(d\mathtt{T}) \int_{\mathbb{R}_+^2}  \mathtt{T}(t+s,\mathtt{y}+s) \breve{\mathcal{P}}_{\mathtt{A}\mathtt{R},i}^{(n)}(ds,d\mathtt{y}) = \frac{1}{\mathrm{m}^{(n)}_{\mathtt{L},i}}\int_t^\infty  \mathtt{T}_i^{(n)}(s) ds
  \eeqnn
  and hence $\psi_i^{(n)}/n$ can be well approximated by $\hat\psi_i^{(n)}$ with
  \beqnn
  \hat\psi_i^{(n)}(t)
  \ar:=\ar \frac{ \Xi_i^{(n)}(0)/n}{\mathrm{m}^{(n)}_{\mathtt{L},i}}\int_t^\infty  \mathtt{T}_i^{(n)}(s) ds\cr
  \ar=\ar \frac{\Xi_i^{(n)}(0)/n }{\mathrm{m}^{(n)}_{\mathtt{L},i}\mathrm{m}^{(n)}_{ii}}\int_t^\infty \zeta^{(n)}_{ii}(s)ds = Z^{(n)}_i(0) I^{(n)}_{\zeta,ii}(t),\quad t\geq 0.
  \eeqnn
  Similarly as in proof of Lemma~\ref{ConvergenceInitialStateCMJ}, we  can prove $\|\psi_\mathcal{H}^{(n)}/n-\hat\psi_\mathcal{H}^{(n)}\|_{L^{1,\infty}}\overset{\rm d}\to0$ as $n\to\infty$. 
  
  \subsection{Proof for Theorem~\ref{ScalingLimitCumulativeCMJ}}
  
  We first consider the cumulative characteristic of all ancestors.
  From Condition~\ref{ConditionInitialStateCMJ}, for any $i\in\mathcal{H}$ we have $\Xi^{(n)}_i(0)=O(n)$ as $n\to\infty$ and hence
  \beqnn
  \mathbf{E}\Big[ \sup_{t\geq 0}\psi^{(n)}_{i}(nt)/n^{2}\Big]
  \leq \mathbf{E}\big[ n^{-2}\Xi_i^{(n)}(0)\big] \cdot \mathtt{T}_{i}^{(n)}(\infty) \to 0.
  \eeqnn
  
  For $i\in\mathcal{H}$ and $j\in\mathcal{H}_i$, we now prove that the cumulative impact of all type-$i$ offspring born by type-$j$ parents can be ignored. Indeed, since $\zeta_i(t,\boldsymbol{u})$ is non-decreasing in $t$, by Lemma~\ref{MomentEstimate} we have for  $j\in\mathcal{H}$, 
  \beqnn
  \mathbf{E}\Big[  \sup_{t\in[0,1]}	\frac{1}{n^2}\int_0^{nt} \int_\mathbb{U}  \zeta_i(nt-s,\boldsymbol{u}) N^{(n)}_j(ds,d\boldsymbol{u}) \Big] 
  \ar\leq\ar
  C  \int_\mathbb{U}  \psi_i(\infty,\boldsymbol{u}) \nu^{(n)}_j(d\boldsymbol{u}) \cr
  \ar=\ar C \cdot \mathtt{T}^{(n)}_i(\infty) \cdot \mathrm{m}^{(n)}_{ij},
  \eeqnn
  which goes to $0$ as $n\to\infty$ since $\mathrm{m}_{ij}^{(n)}\to0$; see Condition~\ref{ConvergenceParameterCMJ}.
  Similarly, the cumulative characteristic of immigrants also can be ignored, i.e.,
  \beqnn
  \mathbf{E}\Big[  \sup_{t\in[0,1]}	\frac{1}{n^2}\int_0^{nt} \int_\mathbb{U}  \zeta_i(nt-s,\boldsymbol{u}) N^{(n)}_I(ds,d\boldsymbol{u}) \Big] 
  \ar\leq\ar
  \frac{C}{n} \int_\mathbb{U} \zeta_i(\infty,\boldsymbol{u}) \nu^{(n)}_I(d\boldsymbol{u})\cr
  \ar=\ar \frac{C}{n} \cdot \mathtt{T}_{i}^{(n)}(\infty)\cdot \mathrm{m}_{iI}^{(n)},
  \eeqnn
  which goes to $0$ as $n\to\infty$.
  In conclusion, the asymptotic behavior of $\mathbf{T}^{(n)}_i(nt)/n^2$ is fully determined by
  \beqnn
  \frac{1}{n^2}\int_0^{nt} \int_\mathbb{U}  \zeta_i(nt-s,\boldsymbol{u}) N^{(n)}_i(ds,d\boldsymbol{u}),
  \eeqnn
  whose weak convergence can be obtained by using Theorem~\ref{ConvergenceCumulativeEffect}.

  \subsection{Proof for Theorem~\ref{ConvergencePopulationStructure}}
  Let $C_{\rm Lip}(\mathbb{R}_+^2)$ be the space of Lipschitz continuous functions on $\mathbb{R}_+^2$ with compact support, which is dense in $C_0(\mathbb{R}_+^2)$.
  From Theorem~9.1 in \cite[p.142]{EthierKurtz2005},  it suffices to prove that the following two claim hold.
  \begin{enumerate}
  	\item[(a)]  The sequence $\{\mathcal{AR}^{(n)}_{i,nt}/n:t\geq 0\}_{n\geq 1}$ satisfies the compact containment condition, i.e.,
  	for any $\eta, T>0$ there exists a compact set $\Gamma_{\eta,T} \subset \mathcal{M}(\mathbb{R}_+^2)$ such that
  	\beqnn
  	\inf_{n\geq 1} \mathbf{P}\Big( \frac{1}{n}\cdot \mathcal{AR}^{(n)}_{i,nt}  \in \Gamma_{\eta,T}\mbox{ for }t\in[0,T] \Big) \geq 1-\eta;
  	\eeqnn
  	
  	\item[(b)] For any $f\in C_{\rm Lip}( \mathbb{R}_+^2)$, the sequence $\{\mathcal{AR}^{(n)}_{i,nt}(f)/n:t\geq 0\}_{n\geq 1}$ converges weakly to $\{ \Xi_i^*(t) \cdot \breve{\mathcal{P}}_{\mathtt{A}\mathtt{R},i}^*(f):t\geq 0\}$ in $\mathbf{D}([0,\infty),\mathbb{R})$ as $n\to\infty$.
  	
  \end{enumerate}
  
  Indeed, for any $f\in C_{\rm Lip}(\mathbb{R}_+^2)$ or $f \equiv 1$, we have $\mathcal{A}\mathcal{R}^{(n)}_{i,t}(f)= \mathbf{T}^{(n)}_i(t)$ with $\mathcal{P}^{(n)}_{\mathtt{T},i}(\mathtt{T}(t,y)= f (t, y-t)\cdot\mathbf{1}_{\{y-t>0\}})=1$.
  In this case, we have $\|\mathtt{T}(y)\|_{L^1} + \|\mathtt{T}(y)\|_{\rm TV} \leq  C (1+y)$ for any $y\geq 0$. By Condition~\ref{ConvergenceLifeDis}, we see that Condition~\ref{CMJInstantaneousConvergence} is satisfied with
  \beqnn
  \|\mathtt{T}_i^{(n)}\|_{L^1} 
  \ar=\ar\int_0^\infty dt \int_0^\infty f (t, y-t)\cdot\mathbf{1}_{\{y-t>0\}} \mathcal{P}^{(n)}_{\mathtt{L},i}(dy)\cr
  \ar=\ar \int_0^\infty dt \int_0^\infty  f(t,z)  \mathcal{P}^{(n)}_{\mathtt{L},i}(t+dz) \cr
  \ar\to \ar \int_0^\infty dt \int_0^\infty  f(t,z)  \mathcal{P}^*_{\mathtt{L},i}(t+dz)=\mathrm{m}_{\mathtt{L},i}^* \cdot \breve{\mathcal{P}}_{\mathtt{AR},i}^*(f),
  \eeqnn
 as $n\to\infty$. 
  By Theorem~\ref{ConvergenceCMJInstantaneousEffect} and Corollary~\ref{Corollary.4.10}, we have
  \beqnn
  \frac{1}{n}\cdot \mathcal{AR}^{(n)}_{i,nt}(f)  \overset{\rm d}\to  \mathrm{m}_{\mathtt{L},i}^*\mathrm{m}_{ii}^*\cdot \mathbf{B}_i^*(t) \cdot \breve{\mathcal{P}}_{\mathtt{AR},i}^*(f) =\Xi^*_i(t)  \cdot \breve{\mathcal{P}}_{\mathtt{AR},i}^*(f),
  \eeqnn
  in $\mathbf{D}([0,\infty),\mathbb{R})$ as $n\to\infty$. Here we have got claim (b). 
  By the Skorokhod representation theorem, we may assume 
  \beqnn
  \{\mathcal{AR}^{(n)}_{i,nt}(f)/n:t\geq 0\}_{n\geq 1} \overset{\rm a.s.}\to \Xi^*_i  \cdot \breve{\mathcal{P}}_{\mathtt{AR},i}^*(f),
  \eeqnn
  in $\mathbf{D}([0,\infty),\mathbb{R})$. 
  In particular, by Proposition~1.17 in \cite[p.328]{JS03} and the continuity of $\Xi^*_i$, 
  \beqnn
  \{\mathcal{AR}^{(n)}_{i,nt}(1)/n:t\geq 0\}_{n\geq 1} \overset{\rm a.s.}\to \Xi^*_i  \cdot \breve{\mathcal{P}}_{\mathtt{AR},i}^*(1),
  \eeqnn
  uniformly on compacts. For any $T>0$ and $\eta\in(0,1)$, there exists a constant $K>0$ such that
  \beqnn
  \mathbf{P}\bigg(\sup_{t\in[0,T]} \frac{1}{n}\cdot \mathcal{AR}^{(n)}_{i,nt}(1)  \leq K \bigg) \geq 1-\eta
  \eeqnn
  and hence claim (a) holds.
  Here we have got the first convergence result in Theorem~\ref{ConvergencePopulationStructure}.
  
  Similarly as in (\ref{marginaldistribution}), we have $\breve{\mathcal{P}}_{\mathtt{AR},i}^*(dy,\mathbb{R}_+)=\breve{\mathcal{P}}_{\mathtt{AR},i}^*(\mathbb{R}_+,dy)=\breve{\mathcal{P}}_{\mathtt{L},i}^*(dy)$. Hence the second desired convergence result follows from the fact that $\mathcal{A}^{(n)}_{i,t}$ and $\mathcal{R}^{(n)}_{i,t}$ are marginal measures of $\mathcal{AR}^{(n)}_{i,t}$.
  Finally, the third desired convergence result follows directly from the first one together with the facts that
  \beqnn
  \mathcal{L}^{(n)}_{i,t}(dy)= \int_{\mathbb{R}_+^2}\delta_{s+z}(dy) \mathcal{AR}^{(n)}_{i,t}(ds,dz)
  \quad \mbox{and}\quad 
  \int_{\mathbb{R}_+^2}\delta_{s+z}(dy)  \breve{\mathcal{P}}_{\mathtt{AR},i}^*(ds,dz)= \mathring{\mathcal{P}}_{\mathtt{L},i}^* (dy).
  \eeqnn

  \begin{appendix}
  	\section*{Stochastic differential equations driven by $\mathbb{H}^\#$-semimartingales} \label{AppendixB}
  	
  	In this section we give a brief introduction to stochastic differential equations driven by infinite-dimensional semimartingales; readers may refer to \cite{KurtzProtter1996} for more details. Let $\mathbb{H}$ be an arbitrary, separable Banach space endowed with norm $\|\cdot\|_{\mathbb{H}}$.
  	We now give the definition of  $\mathbb{H}^\#$-semimartingales.
  	
  	\begin{definition}
  		$\boldsymbol{Y}$ is an $(\mathscr{F}_t)$-adapted \textsl{$\mathbb{H}^\#$-semimartingale}, if it is a stochastic process indexed by $\mathbb{H} \times \mathbb{R}_+$ such that
  		\begin{enumerate}
  			\item[$\bullet$] for each $f\in \mathbb{H}$,  $\{ \boldsymbol{Y}(f,t):t\geq 0 \}$ is a   c\'adl\'ag $(\mathscr{F}_t)$-semimartingale with $\boldsymbol{Y}(f,0)\overset{\rm a.s.}=0$;
  			\smallskip
  			
  			\item[$\bullet$] for each $t\geq 0$, $\alpha_1,\cdots,\alpha_m\in\mathbb{R}$ and $f_1,\cdots,f_m\in\mathbb{H}$, 
  			\beqnn
  			\boldsymbol{Y}\bigg(\sum_{k=1}^m \alpha_k f_k,t \bigg)\overset{\rm a.s.}=\sum_{k=1}^m \alpha_k \boldsymbol{Y}\big(f_k,t\big).
  			\eeqnn
  		\end{enumerate}
  	\end{definition}
  	
 Let $\mathbb{H}_0$ be a dense subset of $\mathbb{H}$ and $\mathcal{S}_0$ the collection of $\mathbb{H}$-valued stochastic processes of  the form
 \beqnn
 \boldsymbol{X}(t):= \sum_{k=1}^m \xi_k(t)\varphi_k \quad \mbox{with}\quad \xi_k(t):=\sum_{i=0}^{\infty} \eta_i^k \cdot \mathbf{1}_{[\tau_i^k,\tau_{i+1}^k)}(t),
 \quad t\geq 0,
 \eeqnn
 where $m\geq 1$, $\varphi_1,\cdots,\varphi_m\in\mathbb{H}_0$, $\{\tau_i^k\}_{i\geq 0}$ is a sequence of non-decreasing $(\mathscr{F}_t)$-stopping times and $\eta_i^k \in\mathbb{R}^d$ is  $\mathscr{F}_{\tau_i^k}$-measurable. For any $\boldsymbol{X} \in\mathcal{S}_0$, we define
 \beqnn
 \boldsymbol{X}_-\cdot \boldsymbol{Y}(t) =\sum_{k=1}^m \int_0^t \xi_k(s-)d \boldsymbol{Y}(t,\varphi_k), \quad t\geq 0.
 \eeqnn
 The $\mathbb{H}^\#$-semimartingale $\boldsymbol{Y}$ is \textit{standard} if 
 \beqnn
 \mathcal{H}_t:=  \bigg\{ \sup_{s\leq t}|\boldsymbol{X}_-\cdot \boldsymbol{Y}(s)| : \boldsymbol{X} \in\mathcal{S}_0,\, \sup_{s\leq t}\|\boldsymbol{X}(s)\|_{\mathbb{H}}\leq 1  \bigg\}
 \eeqnn
 is stochastically bounded for each $t\geq 0$. In this case, for any $\mathbb{H}$-valued c\'adl\'ag process $\boldsymbol{X}$, we can find a sequence $\{\boldsymbol{X}^\epsilon\}_{\epsilon>0}\subset \mathcal{S}_0$ such that as $\epsilon \to0$,
 \beqnn
 \sup_{t\in[0,T]}\|\boldsymbol{X}^\epsilon(t)-\boldsymbol{X}(t)\|_\mathbb{H}\overset{\rm  a.s. }\to0
 \quad \mbox{and} \quad
 \boldsymbol{X}_-\cdot \boldsymbol{Y} \equiv \lim_{\epsilon\to 0+} \boldsymbol{X}^\epsilon_-\cdot \boldsymbol{Y}
 \eeqnn
 exists a.s. in the sense that $\sup_{t\in[0,T]} |\boldsymbol{X}_-\cdot \boldsymbol{Y}(t)- \boldsymbol{X}^\epsilon_-\cdot \boldsymbol{Y}(t) | \overset{\rm p}\to 0$.
 Moreover, the limit process $\boldsymbol{X}_-\cdot \boldsymbol{Y}$ is c\'adl\'ag,  independent of $\{\boldsymbol{X}^\epsilon\}_{\epsilon>0}$ and called the \textit{stochastic integral} of $\boldsymbol{X}$ with respect to $\boldsymbol{Y}$. For any $(\mathscr{F}_t)$-stopping time $\sigma$, we have 
 \beqnn
 \boldsymbol{X}_-\cdot \boldsymbol{Y}(t\wedge \sigma)= \boldsymbol{X}_-^\sigma\cdot \boldsymbol{Y}(t)
 \quad \mbox{and}\quad 
 \boldsymbol{X}_-^\sigma(t):= \boldsymbol{X}_-(t) \mathbf{1}_{[0,\sigma)}(t),
 \quad t\geq 0. 
 \eeqnn  
  	
 \begin{definition}\label{Definition.A1}
 We say  a sequence of $\mathbb{H}^\#$-semimartingales $\{\boldsymbol{Y}_n\}_{n\geq 1}$ is \textsl{uniformly tight} if $\{\mathcal{H}_{n,t}\}_{n\geq 1}$ is uniformly stochastically bounded for any $t\geq 0$ with $\mathcal{H}_{n,t}$ is defined as $\mathcal{H}_{t}$ with $\boldsymbol{Y}$ replaced by $\boldsymbol{Y}_n$.
  		
 Moreover, we say it \textsl{converges weakly} to $\boldsymbol{Y}$ and write $\boldsymbol{Y}_n \Rightarrow \boldsymbol{Y}$ if 
 \beqnn
 (\boldsymbol{Y}_n(f_1),\cdots, \boldsymbol{Y}_n(f_m))\overset{\rm d}\to(\boldsymbol{Y}(f_1),\cdots, \boldsymbol{Y}(f_m)),
 \eeqnn
 in $\mathbf{D}([0,\infty),\mathbb{R}^m)$ as $n\to\infty$, for any $m\geq 1$ and $f_1,\cdots,f_m\in\mathbb{H}$.
 \end{definition}

  \end{appendix}

   \begin{appendix}
  	\section*{Proof for Proposition~\ref{Thm.UpperBound}} \label{AppendixA}
   
   We first have $ \big|\hat\phi^{(n)}_{ii}(\lambda) \big|\leq \|\phi^{(n)}_{ii}\|_{L^1} $ for any $n\geq1$ and $\lambda\in\mathbb{R}$.
   Moreover, for any $\epsilon>0$ we can find a nonnegative smooth  function $g_\epsilon(t) $ on $\mathbb{R}_+ $ satisfying that  $ \|\phi^{(n)}_{ii}-g_\epsilon\|_{L^1}\leq \epsilon$ and $ \|g_\epsilon\|_{\rm TV}\leq \|\phi^{(n)}_{ii}\|_{\rm TV}$.
   Denote by $\hat{g}_\epsilon $ the Fourier transform of $g_\epsilon$.
   Then we have
   \beqnn
   \big| \hat\phi^{(n)}_{ii}(\lambda)\big|
   \ar\leq\ar
   \big|\hat\phi^{(n)}_{ii}(\lambda)-\hat{g}_\epsilon (\lambda)\big|+  \big|\hat{g}_\epsilon (\lambda)\big|
   \leq \|\phi^{(n)}_{ii}-g_\epsilon \|_{L^1}+\big|\hat{g}_\epsilon (\lambda)\big|.
   \eeqnn
   By the differentiation property of Fourier transform, we have
   \beqnn
   \hat{g}_\epsilon (\lambda)= \frac{1}{\lambda}\int_0^\infty  e^{\mathrm{i}\lambda t} \frac{\partial}{\partial t} g_\epsilon (t)  dt
   \quad\mbox{and hence}\quad
   \big|\hat{g}_\epsilon (\lambda)\big| \leq \frac{\|g_\epsilon \|_{\rm TV}}{|\lambda|} \leq \frac{\|\phi^{(n)}_{ij}\|_{\rm TV}}{|\lambda|} . 
   \eeqnn
   From these estimates and the arbitrariness of $\epsilon$, we have
   $| \hat\phi^{(n)}_{ii}(\lambda) | \leq  \|\phi^{(n)}_{ii}\|_{\rm TV}/|\lambda|$ and hence
   \beqnn
   \big| \hat\phi^{(n)}_{ii}(\lambda)\big| \leq \|\phi^{(n)}_{ii}\|_{L^1} \wedge  \frac{\|\phi^{(n)}_{ii}\|_{\rm TV}}{|\lambda|}.
   \eeqnn
   The first inequality in (\ref{LowerBound}) follows directly from (\ref{ShortMemory}).
   
   We now start to prove the second inequality in (\ref{LowerBound}). From the hypothesis \ref{H2} and Condition~\ref{Con.PhiH}, there exist constants $n_0\geq 1$ and $T_0>0$ such that for any $n\geq n_0$,
   \beqnn
   \int_{T_0}^\infty t \phi_{ii}^{(n)}(t)  dt\leq   \int_{T_0}^\infty t  \bar{\phi}_{i}(t) dt  \leq \frac{1}{8} \cdot\sigma_i
   \quad \mbox{and}\quad
   \int_0^\infty t  \phi_{ii}^{(n)}(t) dt \geq \frac{3}{4}\cdot\sigma_i .
   \eeqnn
   Since $\cos(x)\geq 1/2$ for any $|x|\leq 1$, we have for any $|\lambda|\leq 1/T_0$,
   \beqnn
   \frac{\partial}{\partial \lambda}\int_0^\infty \sin(\lambda t) \phi_{ii}^{(n)}(t)  dt\ar=\ar \int_0^\infty \cos(\lambda t)\cdot t\cdot \phi_{ii}^{(n)}(t) dt \cr
   \ar\geq\ar  \int_0^{T_0} \frac{1}{2}\cdot t \phi_{ii}^{(n)}(t)  dt-\int_{T_0}^\infty t  \phi_{ii}^{(n)}(t)   dt\geq \frac{3}{16}\cdot \sigma_i.
   \eeqnn
   By the mean value theorem, we have for any $|\lambda|\leq 1/T_0$,
   \beqnn
   \big|1-\hat{\phi}_{ii}^{(n)}(\lambda)\big|  \geq  \Big|\int_0^\infty \sin(\lambda t) \phi_{ii}^{(n)}(t) dt\Big|\geq \frac{3}{16} \cdot\sigma_i\cdot |\lambda|.
   \eeqnn
   Here we have proved the desired result for $|\lambda|\leq 1/T_0$.
   For $|\lambda|>1/T_0$, from Proposition~\ref{Thm.UpperBound}, there exists a constant $\lambda_0>0$ such that
   \beqnn
   \big|\hat\phi_{ii}^{(n)}(\lambda) \big|\leq \frac{1}{2}
   \quad\mbox{and hence}\quad
   \big|1-\hat\phi_{ii}^{(n)}(\lambda) \big|\geq \frac{1}{2}, 
   \eeqnn
  for any $n\geq 1$ and $|\lambda|\geq \lambda_0$.
   It is obvious that the desired result follows if  $1/T_0\geq \lambda_0$  and then the proof ends.
   If $1/T_0<\lambda_0$, it suffices to prove (\ref{LowerBound}) holds for $\lambda \in[1/T_0,\lambda_0]$.
   Notice that
   \beqnn
   \big|1-\hat\phi_{ii}^{(n)}(\lambda)\big|
   \geq  1-\int_0^\infty \cos(\lambda t)\phi_{ii}^{(n)}(t)dt=: 1- F^{(n)}(\lambda).
   \eeqnn
   The continuity of $F^{(n)}$ induces that $\lambda_n:=\arg\max_{|\lambda|\in[1/T_0,\lambda_0]} F^{(n)}(\lambda)$ is well defined.
   For any $T>0$, since $\cos(t)\leq 1$ we have
   \beqlb\label{UpperboundF} 
   F^{(n)}(\lambda_n)
   \ar\leq\ar \int_0^T \cos(\lambda_n t)\phi_{ii}^{(n)}(t)dt + \int_T^\infty \phi^{(n)}_{ii}(t)dt.
   \eeqlb
   Using the hypothesis \ref{H2} again, we can choose $T>0$ large enough such that
   \beqnn
   \sup_{n\geq 1}\int_T^\infty  \phi^{(n)}_{ii}(t)dt \leq \frac{1}{T}\int_T^\infty  t \cdot \bar\phi_{i}(t)dt \leq \frac{1}{2}
   \eeqnn
    and hence 
   \beqnn
   \inf_{n\geq 1}\int_0^T  \phi^{(n)}_{ii}(t)dt \geq \frac{1}{2}.
   \eeqnn
   By the periodicity of $\cos(\lambda_nt)$, we have
   \beqlb\label{Sum}
   \int_0^T \cos(\lambda_n t)\phi_{ii}^{(n)}(t)dt
   \ar\leq\ar \sum_{k=0}^{ [T\lambda_n /(2\pi)] } \int_{(2k\pi-\pi/2)/\lambda_n}^{{(2k\pi+\pi/2)/\lambda_n}} \cos(\lambda_n t) \phi_{ii}^{(n)}(t)dt.
   \eeqlb
   We now start to analyze the maximum of the sum above.
   Notice that $\cos(\lambda_nt)$ is unimodal on each interval $[(2k\pi-\pi/2)/\lambda_n, (2k\pi+\pi/2)/\lambda_n]$ for any $k\geq 0$ with the maximum arrived at the point $2k\pi /\lambda_n$.
   Thus the more weight of $\phi_{ii}^{(n)}$ is distributed around the local maximum points, the larger the sum above will be.
   To obtain the maximum of the summation in (\ref{Sum}) we should split the weight of $\int_0^T \phi^{(n)}_{ii}(t)dt$ uniformly around these maximum points.
   In precise, we choose $T>0$ large enough such that
   \beqnn
   R_{\lambda_n}:= \frac{\lambda_n\int_0^T \phi^{(n)}_{ii}(t)dt}{2\|\phi^{(n)}_{ii}\|_{\rm TV}\cdot ([T\lambda_n /(2\pi)]+1)}<1.
   \eeqnn
   From the previous observation and the fact that $\cos(\lambda_n t) \leq 1$, we have for any $k\geq 0$,
   \beqnn
   \int_{(2k\pi-\pi/2)/\lambda_n}^{{(2k\pi+\pi/2)/\lambda_n}} \cos(\lambda_n t) \phi_{ii}^{(n)}(t)dt
   \ar\leq \ar\|\phi^{(n)}_{ii}\|_{\rm TV}  \int_{(2k\pi-R_{\lambda_n})/\lambda_n}^{{(2k\pi+R_{\lambda_n})/\lambda_n}} \cos(\lambda_n t)dt\cr
   \ar=\ar
   \frac{\|\phi^{(n)}_{ii}\|_{\rm TV}}{\lambda_n}  \int_{ -R_{\lambda_n} }^{{ R_{\lambda_n} }} \cos(t)dt \cr
   \ar=\ar  \frac{ \int_0^T \phi^{(n)}_{ii}(t)dt}{ [T\lambda_n /(2\pi)]+1}\cdot\frac{\sin(R_{\lambda_n})}{R_{\lambda_n}}.
   \eeqnn
   Taking this back into (\ref{Sum}) and then (\ref{UpperboundF}), we have
   \beqnn
   F^{(n)}(\lambda_n) \leq \frac{\sin(R_{\lambda_n})}{R_{\lambda_n}} \int_0^T \phi^{(n)}_{ii}(t)dt  + \int_T^\infty \phi^{(n)}_{ii}(t)dt
   \eeqnn 
   and hence
   \beqnn
   \inf_{|\lambda|>1/T_0} \big|1-\hat\phi_{ii}^{(n)}(\lambda)\big|
   \ar\geq\ar 1-\|\phi^{(n)}_{ii}\|_{L^1}+  \Big(1-\frac{\sin(R_{\lambda_n})}{R_{\lambda_n}} \Big)\cdot \int_0^T \phi^{(n)}_{ii}(t)dt \cr
   \ar\geq\ar \frac{1}{2}\Big(1-\frac{\sin(R_{\lambda_n})}{R_{\lambda_n}} \Big) .
   \eeqnn
   From the fact that $\lambda_n\in(1/T_0,\lambda_0)$ for any $n\geq 1$ and $\sup_{n\geq 1} \|\phi^{(n)}_{ii}\|_{\rm TV}<\infty $, we  have 
   \beqnn
   \inf_{n\geq 1}R_{\lambda_n}>0
   \quad\mbox{and hence}\quad
   \sup_{n\geq 1}\frac{\sin(R_{\lambda_n})}{R_{\lambda_n}}<1 .
   \eeqnn 
   Consequently, there exists a constant $C_0>0$ such that for any \beqnn
   \inf_{n\geq n_0}	\inf_{|\lambda|>1/T_0}	\big|1-\hat\phi_{ii}^{(n)}(\lambda)\big| \geq C_0(|\lambda|\wedge 1).
   \eeqnn
   Putting all estimates above together, we can immediately get the desired result with $C_2:= (3\sigma_i/16) \wedge (1/2) \wedge C_0$.
   \qed
    
   \end{appendix}

 \begin{acks}[Acknowledgments]
  The author would like to thank the two professional referees and Yuchen Sun for their careful and insightful reading of the paper, and for comments, which led to many improvements.  
 \end{acks}

 \end{document}